\documentclass[a4paper,11pt]{article}

\usepackage[utf8]{inputenc}
\usepackage[T1]{fontenc}

\def\thebibliography#1{\section*{Literature Cited\markboth
{Literature Cited}{Literature Cited}}\list
{[\arabic{enumi}]}{\settowidth\labelwidth{[#1]}\leftmargin\labelwidth
\advance\leftmargin\labelsep 
\usecounter{enumi}}
\def\newblock{\hskip .11em plus .33em minus -.07em}
\sloppy
\sfcode`\.=1000\relax}

\usepackage{amssymb}
\usepackage{mathrsfs}
\usepackage{amsmath}
\usepackage{amsfonts}
\usepackage{amsthm}
 
\newtheorem{Th}{Theorem}
\newtheorem{Propn}[Th]{Proposition} 
 
\newtheorem{Lmm}[Th]{Lemma} 
\newtheorem{Hy}{Hypothesis} 
\newtheorem{Defn}[Th]{Definition}

\begin{document}
\begin{center}

{\LARGE 
On the Stochastic Burgers Equation and the Axiom of Choice}
\vspace{5mm} 

John M. Noble\\ 
Matematiska institutionen,\\
Linköpings universitet,\\
58183 LINKÖPING, Sweden
\end{center}

\begin{abstract}
This article considers the stochastic Burgers  equation

\[  \left\{ \begin{array}{l} \partial_t u^{(\epsilon)} =
(\frac{\epsilon}{2}u_{xx}^{(\epsilon)} -
\frac{1}{2}(u^{(\epsilon)2})_x)dt + \partial_t \zeta_x\\u_0 \equiv
0\end{array}\right.
\]

\noindent where $\zeta$ is a spatially homogeneous Gaussian random field, $2\pi$-periodic in the space variable, mean zero and with covariance $E_{\bf Q} \left \{\zeta(t,x)\zeta(s,y)\right\} = (s \wedge t) \sum_{n \geq 1} a_n^2 \cos (n(x-y)) = (s \wedge t)\Gamma(x-y)$ where $\Gamma$ is $8$ times differentiable. The main result is that 
\[ \lim_{\epsilon \rightarrow 0} \sup_{0 \leq t \leq T} \sup_{x \in [0,2\pi]} \left |E_{\bf Q} \left\{ u^{(\epsilon)2p}(t,x) \right\} -   (-\Gamma^{\prime \prime}(0))^p \prod_{j=1}^p (2j - 1)t^p \right | = 0.\]
\noindent for each positive integer $p$ and each $T < +\infty$. This result is of interest following the work of E, Khanin, Mazel and Sinai, proving existence of invariant measure for the stochastic inviscid Burgers equation, where the hypotheses on the Gaussian random field include those of this article. The method of E, Khanin, Mazel and Sinai is to construct a solution to the Stochastic Inviscid Burgers equation using the minimising trajectories of the associated action functional. This construction relies on the fact that the minimiser exists, which depends on the relative weak compactness of the unit ball in $L^2$ (Tychonov compactness). Kelley proved that Tychonov compactness is equivalent to the Axiom of Choice. This article therefore demonstrates that the Axiom of Choice leads to contradictory results in mathematical analysis. 

\end{abstract}

\section{Summary and Notations}
This article considers the stochastic Burgers equation

\begin{equation}\label{burger} \left\{ \begin{array}{l} \partial_t u^{(\epsilon)} =
(\frac{\epsilon}{2}u_{xx}^{(\epsilon)} -
\frac{1}{2}(u^{(\epsilon)2})_x)dt + \partial_t \zeta_x\\u_0 \equiv
0\end{array}\right.
\end{equation}

\noindent where $\partial_t$ denotes a stochastic differential with respect to
the variable $t$, subscripts denote derivative with respect to the argument labelled by the subscripted variable and
$\zeta$ is a space homogeneous, Gaussian, random field,   satisfying the following
hypotheses. 
\begin{Hy} \label{hyp}
 \begin{itemize}
\item  $(\Omega, {\cal G}, ({\cal G}_t)_{0 \leq  t < +\infty}, {\bf Q})$ denotes a filtered probability space, where \\ $(\beta^{jn})_{j=1,2; n \geq 1}$ are standard independent Wiener processes with respect to ${\bf Q}$, $\beta^{jn}(0) = 0$ for each $(j,n) \in \{1,2\} \times {\bf N}$  and ${\cal G}_ t$ is the $\sigma$-algebra generated by the increments \[\{(\beta^{jn}(v) - \beta^{jn}(u))_{0 \leq u \leq v \leq t}, \;(j,n) \in \{1,2\} \times {\bf N}  \}.\]

\item The Gaussian random field $\zeta$ is defined as 

\[ \zeta (t,x) = \sum_{n=1}^\infty a_n (\cos (nx) \beta^{1n}(t) + \sin
(nx)\beta^{2n}(t))\]

\noindent where $(a_n)_{n \geq 1}$ are real numbers satisfying
$\sum_{n=1}^\infty n^4 |a_n| < +\infty$ and
$\beta^{jn}$ are independent standard Brownian motions such that $\beta^{jn}(0) = 0$   with respect to   $(\Omega, {\cal G}, ({\cal G}_t)_{t \geq 0}, {\bf Q})$.
\end{itemize} 
\end{Hy}

\paragraph{Notation} $E_{\bf Q}\{.\}$ denote the expectation operator with respect to the measure
${\bf Q}$. \vspace{5mm}

\noindent $\Gamma$ denotes the function 
\begin{equation}\label{gamdef} \Gamma(x) = \sum_{n \geq 1} a_n^2 \cos (nx).
\end{equation}

\noindent note that $\Gamma \in C^{8}({\bf R})$ (eight times differentiable) and is $2\pi$ periodic.\vspace{5mm}

\noindent Under ${\bf Q}$, $\zeta(0,.) \equiv 0$ ${\bf Q}$-almost surely and $\zeta$ is Gaussian  satisfying $E_{\bf Q}\{\zeta\} = 0$,  with
covariance given by

\[ E_{\bf Q} \left \{\zeta(t,x)\zeta(s,y)\right\} = (s \wedge t) \sum_{n \geq 1} a_n^2 \cos (n(x-y)) = (s \wedge t)\Gamma(x-y);\]

\noindent The moment fields are considered. Firstly, a priori bounds are calculated for

\begin{equation}\label{momfield} m^{(\epsilon)}_p(t,x_1,\ldots, x_p) := E_{\bf Q} \left \{ u^{(\epsilon)}(t,x_1)\ldots
u^{(\epsilon)}(t,x_p)\right\}.\end{equation}

\noindent  for $t \in [0,T]$ where  $T < +\infty$. These bounds are 
independent of
$\epsilon$. Secondly, the moment fields are shown to be Lipschitz, with the
Lipschitz constant independent of $\epsilon$. Thirdly, the main results of the article are the following theorems: 

\begin{Th} \label{last} Let $u^{(\epsilon)}$ denote the solution to equation (\ref{burger}), where $\zeta$ satisfies hypothesis~\ref{hyp}. Then for all non negative integer $p$, 
\[ E_{\bf Q} \left\{ u^{(\epsilon)2p + 1}(t,x) \right\} = 0 \qquad \forall t \geq 0, \; \epsilon > 0, \; \forall x \in [0, 2\pi]\]
\noindent and for all positive integer $p$ and all $T < +\infty$,

\begin{equation}\label{momgrow} \lim_{\epsilon \rightarrow 0} \sup_{0 \leq t \leq T} \sup_{x \in [0,2\pi]} \left |E_{\bf Q} \left\{ u^{(\epsilon)2p}(t,x) \right\} -   (-\Gamma^{\prime \prime}(0))^p \prod_{j=1}^p (2j - 1)t^p \right | = 0.\end{equation}
\end{Th}

\noindent and 

\begin{Th}\label{secondlast2} There is an adapted function $u: \Omega \times {\bf R}_+ \times {\bf R} \rightarrow {\bf R}$ such that for each $p > 1$ and each $T < +\infty$
\[ \lim_{\epsilon \rightarrow 0} E_{\bf Q}\left\{   \int_0^T \int_0^{2\pi}\left |u^{(\epsilon)}(t,x) - u(t,x) \right |^p dx dt \right  \} = 0.\]
\noindent For each $p > 0$, this function solves

\[ \left\{ \begin{array}{l} \partial_t u = -\frac{1}{2} (u^2)_x dt + \partial_t \zeta_x \\ u(0,x) \equiv 0 \end{array} \right. \]
\end{Th}

\noindent These results ought to be
of interest, following the results of E, Khanin, Mazel and Sinai in~\cite{EKMS}, showing existence of
an invariant measure for this equation. The moments of the invariant
measure constructed by E, Khanin, Mazel and Sinai are discussed in section~\ref{sin}, where it is shown that $E_{\bf Q} \left\{ \sup_{0 \leq x \leq 2\pi} |u(x)|^{p} \right\} < +\infty$ for each $p > 1$, where the distribution of $u$ is the invariant measure for the equation.  

\paragraph{Brief Outline} From equation (\ref{burger}), it is proved (section~\ref{momsect})  that the moment fields defined by equation (\ref{momfield}) satisfy
\begin{eqnarray}\nonumber \lefteqn{ \frac{\partial}{\partial t}m_p^{(\epsilon)}(t;x_1,\ldots, x_p) = \frac{\epsilon}{2}\Delta_{\bf x} m_p^{(\epsilon)}(t;x_1,\ldots, x_p)}\\&& \label{momfeq} - \frac{1}{2} \sum_{j=1}^p \frac{\partial}{\partial x_j} m_{p+1}^{(\epsilon)}(t;x_1,\ldots, x_p,x_j) + \sum_{1 \leq j < k \leq p} (-\Gamma^{\prime \prime} (x_j - x_k))m_{p-2}(t;\hat{x}_j, \hat{x_k})\end{eqnarray}

\noindent where $\frac{\partial}{\partial x_j}$ means derivative with respect to $x_j$ (that is both appearances of $x_j$), $\hat{x}_j$ denotes that variable $x_j$ has been omitted. This requires a Fubini theorem and the use of It\^o's formula. 
 
The non-linearity in the Burgers equation means the $p$th equation depends on the $p+1$ moment field. To show that there is a well defined solution to this system of equations that gives the moment fields, several a-priori bounds on moments of the solution to equation (\ref{burger}) and its derivatives have to be computed. The  bounds on the moments and the first derivatives have to be uniform in $\epsilon$ to ensure that the limit, as $\epsilon$ tends to zero can be taken.

    Theorem~\ref{bound1} gives a bound (independent of $\epsilon$) on 

\[E_{\bf Q} \left\{   \sup_{0 \leq \epsilon \leq
1}\sup_{0 \leq t \leq T} \sup_{0 \leq x \leq 2\pi} \left |u^{(\epsilon)}(t,x) \right |^p \right \} \]

\noindent  and theorem~\ref{bound2} gives a bound (independent of $\epsilon$) on

\[\sup_{0
< \epsilon \leq 1}\sup _{0 \leq t \leq T}\sup_{x_1, \ldots, x_p} \left |\frac{\partial}{\partial x_1}E_{\bf Q} \left\{  u^{(\epsilon)}(t,x_1)
\ldots u^{(\epsilon)}(t,x_p) \right \} \right |.\]

\noindent After this, it is shown that the conditions are satisfied so that
It\^ o's formula may be applied and orders of integration exchanged to show that the moment fields given by equation (\ref{momfield})  satisfy equation (\ref{momfeq}). 

 The a priori upper bounds from theorems~\ref{bound1} and~\ref{bound2} are necessary for a-priori existence of solution to these moment equations; without them, there is no proof that the moment equations should have a solution. 
Next, an appropriate rescaling is carried out. Defining

\[ \phi_p^{(\epsilon)}(t;x_1,\ldots, x_p):=
\frac{m^{(\epsilon)}_p(t;\epsilon x_1,\ldots, \epsilon x_p)
- m_p^{(\epsilon)}(t; 0, \ldots, 0)}{\epsilon}\]

\noindent and

\[ \mu_p^{(\epsilon)}(t;x_1,\ldots,x_p) := m_p^{(\epsilon)} (t;\epsilon x_1,\ldots,
\epsilon x_p)\]

\noindent gives

\begin{eqnarray}\nonumber \lefteqn{ \frac{\partial}{\partial t} \mu_p^{(\epsilon)}(t;x_1,\ldots, x_p) = \frac{1}{2}\Delta_{{\bf x}}\phi^{(\epsilon)}_p(t;x_1,\ldots, x_p)}\\&& \label{rescaleeqn} - \frac{1}{2}\sum_{j=1}^p \frac{\partial}{\partial x_j} \phi^{(\epsilon)}_{p+1}(t;x_1,\ldots, x_p, x_j) + \sum_{1 \leq j < k \leq p} \mu^{(\epsilon)}_{p-2}(t;\hat{x}_j, \hat{x}_k) (-\Gamma^{\prime \prime}(\epsilon(x_j - x_k)).\end{eqnarray} 
 
\noindent The result in theorem~\ref{last} only requires the diagonal $m_p^{(\epsilon)}(t;0,\ldots,0) = E_{\bf Q} \left\{u^{(\epsilon)p}(t,x)\right\}$. This holds for all $x \in {\bf R}$ because the distribution of the random field and hence the distribution of $u^{(\epsilon)}(t,.)$ is spatially homogeneous. It is shown that for all $(x_1,\ldots, x_p) \in {\bf R}$ and all $T < +\infty$, 

\[ \lim_{\epsilon \rightarrow 0} \sup_{0 \leq t \leq T} |\mu_p^{(\epsilon)}(t;x_1,\ldots, x_p) - m_p^{(\epsilon)}(t;0,\ldots, 0)| = 0,\]

\noindent and that there is a function $M_p(t)$ such that for all $T < +\infty$ 
\[ \lim_{\epsilon \rightarrow 0} \sup_{0 \leq t \leq T} |m_p^{(\epsilon)}(t;0,\ldots, 0) - M_p(t) | = 0 \]
\noindent and such that for each $T < +\infty$, $M_p(t)$ is   Lipschitz for $t \in [0,T]$.  It is clear that 
\[ \lim_{\epsilon \rightarrow 0} |\Gamma^{\prime \prime}(\epsilon z) - \Gamma^{\prime \prime}(0)| = 0 \qquad \forall z \in {\bf R}.\]

\noindent For each $t \in {\bf R}$, $m_p^{(\epsilon)}(t;.)$ is Lipschitz in the space variable, uniformly in $\epsilon$, from which it follows that $\phi_p^{(\epsilon)}(t;.)$ is Lipschitz in the space variable, uniformly in $\epsilon$. The bounds on the growth, together with interpreting $\frac{1}{2}\Delta$ as the infinitesimal generator of the heat semigroup, enable the computations which show that    the terms containing $\phi^{(\epsilon)}$ do not yield a contribution as $\epsilon \rightarrow 0$. From this, the results  stated in theorem~\ref{last} are obtained. The uniform upper bound in $\epsilon$ from theorem~\ref{bound2} is necessary for the uniform bounds in $\epsilon$ on the Lipschitz constant for $\phi^{(\epsilon)}$ required to take the limit. By taking $\epsilon \rightarrow 0$ in equation (\ref{rescaleeqn}), the proof of theorem~\ref{last} may be completed. 

\section{Bounds on the Moments}    The following lemma is a necessary first step
towards proving  theorems~\ref{bound1} and~\ref{bound2}. 

\begin{Lmm}\label{prel}
Let $\beta^j$ be independent standard Brownian motions with $\beta^j(0) = 0$, let  \[ S^j(t) = \sup_{0
\leq s \leq t}|\beta^j(t)|\]

\noindent  and let $a_j$ be numbers such that $\sum_j |a_j| <
+\infty$. It follows that 

\begin{equation}\label{preleq} E_{\bf Q} \left \{ \exp\left\{\sum_j a_j S^j(t) \right \} \right \} \leq \exp\left \{ \frac{t}{2}\left(\sum_j a_j^2 \right) +
\sqrt{2  t}(\sqrt{\log 2} + 2\sqrt{\pi}) \left( \sum_j |a_j|\right)\right \}.\end{equation}
\end{Lmm}\vspace{5mm}

\noindent {\bf Proof of lemma~\ref{prel}} Firstly, from Revuz and Yor~\cite{RY}
page 55 proposition  (1.8), if $S(t) = \sup_{0 \leq s \leq t} \beta(s)$ where
$\beta$ is a standard Brownian motion with $\beta(0) = 0$, such that  for all $0 \leq s \leq u \leq v \leq t < +\infty$, $\beta(v) - \beta(u)$ is measurable  with respect to
$ {\cal G}_{t}$  then, using ${\bf
Q}\{.\}$ to denote the {\em probability} of an event with respect to ${\bf
Q}$, 

\[ {\bf Q}\{ S(1) \geq x \} \leq \exp \left \{ -\frac{x^2}{2} \right \}\]

\noindent and rescaling gives

\[ {\bf Q}\{ S(t) \geq x \} \leq \exp\left \{ -\frac{x^2}{2t}\right \}.\]

\noindent Let $\tilde{S}(t) = \sup_{0 \leq s \leq t}|\beta(s)| = \sup_{0 \leq s \leq
t}\beta(s) \vee \sup_{0 \leq s \leq t}(-\beta(s))$. Note that 

\begin{eqnarray}\nonumber {\bf Q}\left \{\tilde{S}(t) \geq x \right \} &=& {\bf Q}\left\{ \left \{\sup_{0
\leq s \leq t}\beta(s)
\geq x \right \}
\cup\left \{ \sup_{0 \leq s \leq t}(-\beta(s)) \geq x\right \}\right\}\\
 &\leq & \label{goodineq} 1 \wedge  2{\bf Q} \left \{S(t)
\geq x \right \} \leq 1 \wedge  2\exp \left \{-\frac{x^2}{2t} \right \}.
\end{eqnarray}

\noindent Note that 
\[2e^{-x^2/2t} = 1 \quad \Rightarrow \quad x = \sqrt{2t \log 2}\]
($\log$ means natural logarithm). 
It follows that for any $\alpha > 0$, 

\begin{eqnarray*} E_{\bf Q} \left\{ e^{\alpha \tilde{S}(t)}\right\} &=& \int_0^\infty {\bf
Q}\left \{e^{\alpha \tilde{S}(t)} \geq y\right \} dy\\ 
&=& \int_0^\infty {\bf Q} \left\{ \tilde{S}(t) \geq \frac{1}{\alpha} \log y\right \} dy\\
&\leq & 
 e^{\alpha \sqrt{2t \log 2}} + \int_{e^{\alpha \sqrt{2t \log 2}}}^\infty {\bf
Q}\left \{ \tilde{S}(t)  \geq \frac{1}{\alpha} \log y \right \} dy\\ 
\\ &\leq &
e^{\alpha \sqrt{2t \log 2}} + 2\int_{e^{\alpha\sqrt{2t \log 2}}}^\infty
e^{-\frac{1}{2\alpha^2 t} (\log y)^2} dy
\end{eqnarray*}

\noindent Substituting $x = \frac{\log y}{\alpha\sqrt{t}}$ so that $y = e^{\alpha \sqrt{t} x}$, 

\begin{eqnarray*} E_{\bf Q} \left\{ e^{\alpha \tilde{S}(t)}\right\}&=& e^{\alpha \sqrt{2t \log 2}} + 2 \alpha \sqrt{t}\int_{\sqrt{2\log 2}}^\infty e^{-\frac{x^2}{2} + \alpha \sqrt{t} x} dx \\ 
&=&
e^{\alpha\sqrt{2t\log 2}} + 2\alpha \sqrt{ t} e^{\alpha^2 t/2}
\int_{\sqrt{2\log 2} - \alpha \sqrt{t}}^\infty e^{-\frac{x^2}{2}  
 } dx 
\\ &\leq  &
\exp\left \{\alpha\sqrt{2t \log 2} \right\} + 2\alpha \sqrt{2\pi t} 
\exp\left \{\frac{\alpha^2 t}{2}\right \}.
\end{eqnarray*}

\noindent Now note that for any non negative numbers $a,b,c,\alpha$,

\[ e^{a\alpha} + b\alpha e^{c\alpha^2} \leq e^{a\alpha + c\alpha^2}(1+b\alpha)\leq
e^{(a+b)\alpha + c\alpha^2}.\]

\noindent It follows that, for any $\alpha \geq 0$, 

\[E_{\bf Q}\left \{ e^{\alpha \tilde{S}(t)}\right \} \leq \exp\left \{ \sqrt{2t}(\sqrt{\log 2} +
2\sqrt{\pi})\alpha + \frac{\alpha^2 t}{2}\right \}.\]

\noindent It follows that 

\begin{eqnarray*} E_{\bf Q}\left \{ e^{\sum_j a_j S^j(t)}\right \} &=& \prod_{j=1}^\infty
E_{\bf Q}\left  \{ e^{  a_j S^j(t)}\right \} \\
&\leq &
\exp\left \{ \frac{t}{2}\left(\sum_j a_j^2 \right) + (\sqrt{2t}(\sqrt{\log 2} +
2\sqrt{\pi}))\left( \sum_j
  |a_j|\right)  
\right \} 
\end{eqnarray*} 

\noindent and lemma~\ref{prel} is proved. \qed\vspace{5mm}

\noindent The following bound will also be useful in the sequel.\vspace{5mm}

\begin{Lmm}\label{prel2} Let $\beta^j$ be independent standard Brownian motions with $\beta^j(0) = 0$ 
and let

\[ S^j(t) = \sup_{0 \leq s \leq t}\left |\beta^j(s) \right |.\]

\noindent  Let $a_j$ be real numbers such that
$\sum_j |a_j| < +\infty$. Let 

\begin{equation}\label{thecons}
 G(p) = \int_0^\infty y^p \exp\left \{-\frac{y^2}{2} \right \}
dy = 2^{(p-1)/2} \Gamma_{Eu} \left( \frac{p+1}{2}\right),
\end{equation}

\noindent where $\Gamma_{Eu}$ denotes the Euler Gamma function $\Gamma_{Eu}(\alpha) = \int_0^\infty z^{\alpha - 1}e^{-z} dz$. It holds that 

\begin{equation}\label{prel2eq} E_{\bf Q} \left \{ \left |\sum_j a_j S^j(t)\right |^p \right \}  \leq \left (\sum_j |a_j| \right )^p t^{p/2}\left((2\log
2)^{p/2} + 2p G(p-1)
\right).\end{equation}

\end{Lmm}\vspace{5mm}

\noindent {\bf Proof} Using inequality (\ref{goodineq}), it follows that 

\[{\bf Q}\{S^j(t) \geq x\} \leq 1 \wedge 2\exp\left
\{-\frac{x^2}{2t}\right \}.\]

\noindent Provided $0 < \sum_j |a_j| < +\infty$, a standard application of Jensen's inequality gives, for $p \geq 1$,  

\begin{eqnarray*}  E_{\bf Q}\left  \{ \left |\sum_j a_j S^j(t) \right |^p \right \}  &\leq&  E_{\bf Q} \left \{ \left (\sum_j |a_j|S^j(t) \right )^p \right \} \\ & = &  \left (\sum_j |a_j| \right )^p E_{\bf Q} \left\{ \left( \sum_j \frac{|a_j|}{\sum_j |a_j|}  S_j(t) \right)^p \right\} \\ &\leq & 
\left (\sum_j |a_j|\right )^{p-1} \sum_j |a_j| E_{\bf Q} \left \{ S^j(t)^p \right \} \end{eqnarray*}

\noindent and

\begin{eqnarray*} E_{\bf Q} \left\{ S^j(t)^p \right\} &=& \int_0^\infty {\bf Q}\left \{S^j(t)^p \geq
x \right \}dx = \int_0^\infty {\bf Q}\left\{ S^j(t) \geq x^{1/p}\right\} dx\\
&\leq & \int_0^\infty \left (1 \wedge 2\exp\left \{-\frac{x^{2/p}}{2t} \right \} \right )dx \\
&=& (2t \log 2)^{p/2} + 2\int_{(2t\log 2)^{p/2}}^\infty e^{-x^{2/p}/2t} dx\\
&=& (2t \log 2)^{p/2} + 2pt^{p/2} \int_{(2\log 2)^{1/2}}^\infty z^{p-1}e^{-z^2/2}
dz\\
&\leq & t^{p/2} \left( (2\log 2)^{p/2} + 2p G(p-1)\right) 
\end{eqnarray*}

\noindent which gives 

\[ E_{\bf Q} \left \{ \left |\sum_j a_j S^j(t) \right |^p \right \} \leq \left (\sum_j |a_j|\right )^p  t^{p/2} \left( (2\log
2)^{p/2} + 2p G(p-1)\right) \]

\noindent as advertised. \qed \vspace{5mm}
 
\begin{Lmm} \label{heatlem} For $p > 2$, $T > 0$, let ${\cal S}_{p,T}$ denote the space of functions $f: {\bf R}_+ \times {\bf R} \times \Omega \rightarrow {\bf R}$, $2\pi$ periodic in the space (second) variable, such that for each $s \in {\bf R}$, $f : [0,s] \times {\bf R} \times \Omega \rightarrow {\bf R}$ is measurable with respect to ${\cal B}([0,s]) \otimes {\cal B}({\bf R}) \otimes {\cal G}_s$  (where ${\cal B}$ denotes the Borel $\sigma$-algebra) and such that 

\begin{equation}\label{normdef} |\|f|\|_{p,T} := \left(\sup_{0 \leq t \leq T} \frac{1}{2\pi}\int_0^{2\pi}  E_{\bf Q} \left \{ |f(t,x)|^p \right \} dx    \right)^{1/p} < +\infty  
\end{equation}

\noindent Let 
\[ {\cal S}_p = \cap_{T > 0} {\cal S}_{p,T} \]

\noindent and let 

\begin{equation}\label{eqspstar} {\cal S}_p^* = \{f \in {\cal S}_p\; | f(t,.) \in C^2({\bf R}) \; \forall t < +\infty \; {\bf Q}- a.s.\}\end{equation}

\noindent For each $\epsilon > 0$, $p \in (0,+\infty)$, there exists a unique solution such that  $U^{(\epsilon)} \in {\cal S}_{p}$   to the
equation 

\begin{equation}\label{heat}\left\{\begin{array}{l} \partial_t U^{(\epsilon)}  = 
\frac{\epsilon}{2}U^{(\epsilon)}_{xx}dt -
\frac{1}{\epsilon} U^{(\epsilon)} \circ \partial_t \zeta \\ U^{(\epsilon)}(0,x)
\equiv 1,\end{array}\right.
\end{equation}

\noindent where the $\circ$ denotes stochastic integration in the
Stratonovich sense. The solution $U^{(\epsilon)}$ satisfies the following regularity: for all $T < +\infty$, there is a version such that almost surely, $U^{(\epsilon)} \in C^{0,\alpha}([0,T] \times {\bf R})$ (Hölder continuous of order $\alpha$) for all $\alpha < \frac{1}{2}$ and for each $t \in [0, T]$, $U^{(\epsilon)}(t,.) \in C^{3,\gamma}({\bf R})$ (three times differentiable, third derivative Hölder continuous of order $\gamma$) for all $\gamma < 1$, where for each $x \in {\bf R}$, $U^{(\epsilon)}(.,x), U^{(\epsilon)}_x(.,x), U^{(\epsilon)}_{xx}(.,x), U^{(\epsilon)}_{xxx}(.,x)$ are Hölder continuous in the time variable of all orders less than $1/2$.   \end{Lmm}

\noindent {\bf Proof} Existence and uniqueness of solution to equation (\ref{heat}) is standard and may be  found (for example) in
Kunita~\cite{Kun}. Kunita also shows that if $\Gamma \in C^{2n}$ then, almost surely, $U^{(\epsilon)}(t,.) \in C^{n-1,\gamma}$ (that is $n-1$ times differentiable, $n-1$th derivative Hölder continuous of order $\gamma$) for all $\gamma < 1$ in the space variable  and if $n \geq 2$, then $U^{(\epsilon)}(.,x)$ and its first $n$ derivatives are  Hölder continous of all orders less than $\frac{1}{2}$ in the time variable. Kunita's results are more extensive; those stated above are the only ones needed here.\vspace{5mm}

\noindent To keep this article relatively self contained, an outline of the proof is sketched here.   The It\^o formulation of the mild form of equation (\ref{heat}) is

\begin{eqnarray}\nonumber \lefteqn{ U^{(\epsilon)}(t,x) = }\\&& \nonumber 1 + \frac{1}{\epsilon}\sum_{n \geq 1} a_n \left(\int_0^t P_{t-s}(U^{(\epsilon)}(s,.)\cos(n.))(x)d\beta^{1n}_s + \int_0^t P_{t-s}(U^{(\epsilon)}(s,.)\sin(n.))(x)d\beta^{2n}_s \right)\\&& \label{mild} + \frac{\Gamma(0)}{2\epsilon^2}\int_0^t P_{t-s} U^{(\epsilon)}(s,x) ds \end{eqnarray}

\noindent where for a bounded measurable function $f$, $P_t$ is defined such that 
\begin{equation}\label{eqhsg} P_t f(x) = \int_{-\infty}^\infty \frac{1}{\sqrt{2\pi \epsilon t}}e^{-y^2/2\epsilon t} f(x+y) dy. \end{equation}   

\noindent From Revuz and Yor page 137 definition (2.1) and theorem (2.2), it is sufficient (but not necessary) that $E_{\bf Q} \left\{ \int_0^t f_s^2 ds \right\} < +\infty$ for an adapted measurable function $f$ to ensure that the stochastic integral $\int_0^t f_s d\beta_s^{jn}$ is well defined. Set $U^{(\epsilon,0)} = U^{(\epsilon)}$,  $f_{1k}(x) = \frac{d^k}{dx^k} \cos(x)$ and   $f_{2k}(x) = \frac{d^k}{dx^k}\sin (x)$. Recall the standard result that for two functions $f$ and $g$, both $n$ times differentiable, the $n$th derivative of the product satisfies
\[ (fg)^{(n)} = \sum_{j=0}^n \left(\begin{array}{c} n \\ j \end{array}\right) f^{(j)} g^{(n-j)}.\]

\noindent For $a \geq 1$, let $U^{(\epsilon,a)}$ satisfy

\begin{eqnarray}\label{Usat} U^{(\epsilon,a)}(t,x) &
=& \frac{1}{\epsilon}\sum_{k=0}^a  \left(\begin{array}{c} a \\ k \end{array}\right) \sum_{n \geq 1} n^k a_n \left(\int_0^t P_{t-s}\left( U^{(\epsilon,a-k)}(s,.)f_{1k}(n.)\right)(x) d\beta^{1n}_s \right. \\&& \left. + \int_0^t P_{t-s} \left(U^{(\epsilon,a-k)}(s,.)f_{2k}(n.)\right)(x) d\beta^{2n}_s\right) +  \frac{\Gamma(0)}{2\epsilon^2} \int_0^t P_{t-s} U^{(\epsilon,a)}(s,x)ds \nonumber
\end{eqnarray}

\noindent Suppose that existence and uniqueness of solution in ${\cal S}_{p}$ for all $p \geq 2$ have been established for $b = 0,1,\ldots, a-1$. Consider the iterative sequence: $U_0^{(\epsilon,0)} \equiv 1$ and $U_0^{(\epsilon,a)} \equiv 0$ for $a = 1,2,3$  and 

\begin{eqnarray} \nonumber \lefteqn{U^{(\epsilon,a)}_{m+1}(t,x) = P_t U^{(\epsilon, a)} (0,x)}\\ && + \label{iterate} \frac{1}{\epsilon}\sum_{k=1}^a  \left(\begin{array}{c} a \\ k \end{array}\right) \sum_{n \geq 1} n^k a_n \left(\int_0^t P_{t-s}\left( U^{(\epsilon,a-k)}(s,.)f_{1k}(n.)\right)(x) d\beta^{1n}_s \right. \\&& \left. + \int_0^t P_{t-s} \left(U^{(\epsilon,a-k)}(s,.)f_{2k}(n.)\right)(x) d\beta^{2n}\right) +  \frac{\Gamma(0)}{2\epsilon^2} \int_0^t P_{t-s} U^{(\epsilon,a)}_m(s,x)ds \nonumber \\&& +
\frac{1}{\epsilon} \sum_{n \geq 1}   a_n \left(\int_0^t P_{t-s}\left( U^{(\epsilon,a)}_m(s,.)\cos(n.)\right)(x) d\beta^{1n}_s  + \int_0^t P_{t-s} \left(U^{(\epsilon,a)}_m (s,.)\sin (n.)\right)(x) d\beta^{2n}\right) 
\nonumber \end{eqnarray}

\noindent  Set $D_m^{(\epsilon)} = U_m^{(\epsilon)} - U_{m-1}^{(\epsilon)}$ for $m \geq 1$ and $D_m^{(\epsilon,a)} =  U_m^{(\epsilon,a)} - U_{m-1}^{(\epsilon,a)}$. Then, using $D_m^{(\epsilon)} = D_m^{(\epsilon,0)}$, it follows that for $a = 0,1,2,3$ and $m \geq 1$, 

\begin{eqnarray}\label{eqdee}
 \lefteqn{D_{m+1}^{(\epsilon,a)}(t,x) = \frac{\Gamma(0)}{2\epsilon^2}\int_0^t P_{t-s}D_m^{(\epsilon,a)}ds}\\&& \nonumber
+ \frac{1}{\epsilon}\sum_{n \geq 1}   a_n \left(\int_0^t P_{t-s}\left( D^{(\epsilon,a)}_m(s,.)\cos(n.)\right)(x) d\beta^{1n}_s  + \int_0^t P_{t-s} \left(D^{(\epsilon,a)}_m (s,.)\sin (n.)\right)(x) d\beta^{2n}\right).
\end{eqnarray}

\noindent while

\[D_1^{(\epsilon)}(t,x) = \frac{\Gamma(0)}{2\epsilon^2}t + \frac{1}{\epsilon}\sum_{n \geq 1} a_n \left(\int_0^t P_{t-s}\cos(nx)d\beta_s^{1n} +   P_{t-s}\sin(nx) d\beta_s^{2n}\right) \]

\noindent and, for $a = 1,2,3$,

\begin{eqnarray*}\lefteqn{D_1^{(\epsilon,a)}(t,x) = \frac{1}{\epsilon}\sum_{k=1}^a  \left(\begin{array}{c} a \\ k \end{array}\right) \sum_{n \geq 1} n^k a_n  \left(\int_0^t P_{t-s}\left( U^{(\epsilon,a-k)}(s,.)f_{1k}(n.)\right)(x) d\beta^{1n}_s\right.}\\&&  \hspace{50mm} \left. + \int_0^t P_{t-s} \left(U^{(\epsilon,a-k)}(s,.)f_{2k}(n.)\right)(x) d\beta^{2n}\right), \end{eqnarray*}

\noindent so that for each $p \geq 2$ there is a constant $C_p < +\infty$ such that (using $\Gamma(0) = \sum_{n \geq 1} a_n^2$)

\[ E_{\bf Q} \left\{ |D_1^{(\epsilon)}(t,x)|^p\right\} \leq C_p\left(\frac{\Gamma(0)}{2\epsilon^2} t\right)^p + \frac{C_p t^{p/2}}{\epsilon^p}\Gamma(0)^{p/2}\]

\noindent and 

\[ E_{\bf Q} \left\{ |D_1^{(\epsilon,a)}(t,x)|^p \right\} \leq \frac{ C_p}{\epsilon^p} \left(\sum_{n \geq 1} n^6 a_n^2 \right)^{p/2}t^{(p/2) - 1}\sum_{k=1}^a \int_0^t E_{\bf Q}\left\{   |U^{(\epsilon, a-k)}(s,x)|^p \right\}ds \quad a = 1,2,3\]

\noindent From equation (\ref{eqdee}), it follows that

\begin{eqnarray*} \lefteqn{ D_{m+1}^{(\epsilon, a)2p}(t,x)  \leq 2^{2p}\left( \left( \frac{1}{\epsilon}\sum_{n \geq 1} n^a a_n  \left(\int_0^t P_{t-s}(D_m^{(\epsilon,a)}(s,.)\cos(n.))(x)d\beta^{1n}_s\right. \right. \right.}\\&&  + \left. \left. \left.  \int_0^t P_{t-s}(D_m^{(\epsilon,a)}(s,.)\sin(n.))(x)d\beta^{2n}_s \right)\right)^{2p}       +   \left(\frac{\Gamma(0)}{2\epsilon^2} \int_0^t P_{t-s} D_m^{(\epsilon,a)} (s,.)(x) ds\right)^{2p}\right ).
 \end{eqnarray*} 

\noindent In the following, the constant $K$ may change from line to line. It denotes a positive finite value. Let

\[M_r = \frac{1}{\epsilon}\sum_{n \geq 1} n^a a_n\left(\int_0^r P_{t-s}(D_m^{(\epsilon,a)}(s,.)\cos(n.))(x)d\beta^{1n}_s +  \int_0^r P_{t-s}(D_m^{(\epsilon,a)}(s,.)\sin(n.))(x)d\beta^{2n}_s \right) \]

\noindent defined on the time interval $r \in [0,t]$. From equation (\ref{iterate}), it is clear (under the assumption that $U^{(\epsilon,b)} \in {\cal S}_{p}$ for each $b < a$) that $U^{(\epsilon,a)}_m \in {\cal S}_{p}$ for each $m < +\infty$, from which it follows directly that $D^{(\epsilon,a)}_m \in {\cal S}_{p}$ for each $m < +\infty$ and hence that $M$ is a local martingale. Note that 
\[ \langle M \rangle_r = \frac{1}{\epsilon^2} \sum_{n \geq 1} n^{2a}a_n^2  \int_0^r \left ( (P_{t-s}(D_m^{(\epsilon,a)}(s,.)\cos(n.))(x))^2 + (P_{t-s}(D_m^{(\epsilon,a)}(s,.)\sin(n.))(x))^2 \right ) ds.\]

\noindent It follows from the Burkholder-Davis-Gundy inequality that 

\begin{eqnarray*} \lefteqn{ E_{\bf Q}\left\{  D_{m+1}^{(\epsilon,a)2p}(t,x)\right\}   \leq K E_{\bf Q} \left \{ \left( \frac{1}{\epsilon^2 }\sum_{n \geq 1} n^{2a} a_n^2  \left(\int_0^t (P_{t-s} D_m^{(\epsilon,a)} (s,.)\cos(n.)(x) )^2   ds \right. \right. \right.  }\\&&  \hspace{5mm} + \left. \left. \left. \int_0^t  ( P_{t-s} D_m^{(\epsilon,a)} (s,.)\sin(n.)(x)) )^2  d s  \right)^{p} \right\}      +   K E_{\bf Q} \left \{ \left(\frac{\Gamma(0)}{2\epsilon^2}\int_0^t   P_{t-s} D_m^{(\epsilon,a)} (s,x)  ds\right)^{2p}\right \}\right ).
 \end{eqnarray*} 
 
\noindent Let $\|.\|_p$ denote the norm defined by
\begin{equation}\label{eqpnorm} \|f\|_p(t) = \sup_{0 \leq s \leq t} \left( \frac{1}{2\pi} \int_0^{2\pi} E_{\bf Q}\left\{ \left |f(s,x) \right |\right\} dx \right)^{1/p}
\end{equation}

\noindent Then  straightforward applications of Hölder's inequality give, for $p\geq 2$,  

\[ \|D_{m+1}^{(\epsilon,a)}\|_{p}^{p}(t) \leq    K \frac{  \left (\sum_{n \geq 1} n^{2a} a_n^2 \right )^{p/2}}{\epsilon^{p}} t^{(p/2) -1}\int_0^t \|D_m^{(\epsilon,a)}\|_{p}^{p}(s) ds + K \frac{\Gamma^{p}(0)}{\epsilon^{2p}}t^{p-1}\int_0^t \|D_m^{(\epsilon,a)} \|_{p}^{p}(s) ds.\]

\noindent Note that, since the field is spatially homogeneous, $\|D_m^{(\epsilon, a)} \|_p^p(t) = E_{\bf Q}\left\{|D_m^{(\epsilon,a)}(t,x)|^p \right\}$. Therefore, for all $t \in [0,T]$, for $T < +\infty$,

\[ \|D_{m+1}^{(\epsilon, a)}\|_{p}^{p}(t) \leq C(T)\int_0^t \|D_m^{(\epsilon, a)}\|_{p}^{p}(s) ds\]
\noindent where
\[C(T) = K \frac{ \left (\sum_{n \geq 1} n^{2a} a_n^2 \right )^{p/2}}{\epsilon^{p}}  T^{(p/2) -1} + K \frac{\Gamma(0)^{p}}{\epsilon^{2p}}T^{p-1}\int_0^t E_{\bf Q} \left\{ |D_1^{(\epsilon)}(s,x)|^p \right\} ds,\]

\noindent so that, for $0 \leq t \leq T$, 
\[ \|D_m^{(\epsilon,a)}\|_p(t) \leq \left(\frac{(C(T)  t)^{m-1}}{(m-1)!}\right)^{1/p} \left(\int_0^t \|D_1^{(\epsilon,a)}\|_p^p(s) ds \right)^{1/p}.\] 

\noindent Since $\sum_{n \geq 1}n^8a_n^2 < +\infty$, it follows that $\sum_{n \geq 1} n^{2a} a_n^2 < +\infty$ for $a = 0,1,2,3$, it follows inductively, that $\|D_m^{(\epsilon,a)}\|_p(t)$ is summable. This is clearly true for $a=0$, which implies that 

\[ \sup_{0 \leq s \leq T} E_{\bf Q} \left\{ |D_1^{(\epsilon)}(s,x)|^p \right\}  < +\infty\]

\noindent and hence is true inductively for $a = 1,2,3$. Existence and uniqueness for solutions to equation (\ref{mild}) in ${\cal S}_{p}$ for all $p \in [ 2, +\infty)$  now follows directly by standard Gronwall arguments. \vspace{5mm}

\noindent Now set 
\[ C_{a,p}(t) = \sup_{0 \leq s \leq t} E_{\bf Q} \left \{ |U^{(\epsilon,a)}(t,x)|^p\right\}\]

\noindent The preceeding gives that $\sup_{0 \leq t \leq T} C_{a,p}(t) < +\infty$ for $a = 0,1,2,3$. 

\noindent Now set

 \begin{equation}\label{veeeqn} V^{(a,h)}(t,x) = \frac{U^{(\epsilon, a)}(t, x+h) - U^{(\epsilon, a)}(t,x-h)}{2h},\end{equation}

\noindent $f_{1k}(x) = \frac{d^k}{dx^k}\cos(x)$ and $f_{2k}(x) = \frac{d^k}{dx^k}\sin(x)$. Then, since $U^{(\epsilon, 0)}(0,.) \equiv 1$ and $U^{(\epsilon,a)}(0,.) \equiv 0$ for $a = 1,2,3$, it follows that  $V^{(a,h)}$ satisfies

\begin{eqnarray}\label{mildder2} \lefteqn{V^{(a,h)} (t,x) =        \frac{1}{\epsilon}\sum_{k=0}^a \left(\begin{array}{c} a \\ k \end{array} \right)\sum_{n \geq 1}n^k a_n \left(\int_0^t P_{t-s}(V^{(a-k,h)}(s,.)f_{1k}(n(.+h)))(x)d\beta^{1n}_s \right.} \\&& \nonumber  \hspace{5mm} \left. + \int_0^t P_{t-s}(V^{(a-k,h)}(s,.)f_{2k}(n(.+h)))(x)d\beta^{2n}_s \right)  + \frac{\Gamma(0)}{2\epsilon^2}\int_0^t P_{t-s} V^{(a,h)}(s,.) ds\\&& \nonumber 
+ \frac{1}{\epsilon} \sum_{k=0}^a \left(\begin{array}{c} a \\ k \end{array} \right)\\&& \hspace{5mm} \times \nonumber \sum_{n \geq 1}  n^k a_n \left( \int_0^t P_{t-s}\left (U^{(\epsilon, a-k)}(s,.-h)\left(\frac{f_{1k}(n(.+h) - f_{1k}(n(.-h))}{2h}\right )\right) (x)d\beta^{1n}_s \right. \\&&  \hspace{5mm}  + \int_0^t P_{t-s}\left (U^{(\epsilon, a-k)}(s,.-h)\left ( \frac{f_{2k}(n(.+h)) - f_{2k}(n(.-h))}{2h}\right )\right)(x)d\beta^{2n}_s .\nonumber 
\end{eqnarray}

\noindent Set 
\begin{equation}\label{veebound} K_{a,p,h}(t) = \sup_{0 \leq s \leq t} E_{\bf Q}\left \{ |V^{(a,h)}(t,x)|^p \right \}.
\end{equation}
 
\noindent For all $2 \leq  p < +\infty$,  elementary arguments give the existence of a constant $c_1(p,a,T) < +\infty$ such that 

\begin{eqnarray*} K_{a, p,h}(t) &\leq&   \frac{c_1(p,a,T)}{\epsilon^{p}} \left(\sum_{n \geq 1} n^{6}a_n^2\right)^{p/2} t^{(p/2) - 1} \sum_{k=0}^a     \int_0^t K_{a-k, p,h}(s) ds\\&& +  \frac{c_1(p,a,T)}{\epsilon^{ p}}    \left( \sum_{n \geq 1} n^{8} a_n^2\right)^{p/2}t^{(p/2) - 1}\sum_{k=0}^a  \int_0^t C_{a-k, p}(s) ds \\&&  +  \frac{c_1(p,a,T)}{\epsilon^{2p}} \Gamma(0)^p t^{p-1}\int_0^t K_{a, p,h}(s) ds \end{eqnarray*}

\noindent from which it follows that for all $T < +\infty$,  $2 \leq p< +\infty$ and $a = 0,1,2,3,4$, $C_{a, p}(T)   < + \infty$  and for $a = 0,1,2,3$, $\tilde{K}_{a, p}(T):= \sup_{h > 0} K_{a, p, h}(T) < +\infty$. These bounds enable Kolmogorov's criterion to be applied to the space variable. \vspace{5mm}

\noindent The following computations enable the appropriate Hölder continuity to be proved for the time variable. 
Set
\begin{eqnarray*}I^{(\epsilon,a)}(r,t;x) &
=& \frac{1}{\epsilon}\sum_{k=0}^a  \left(\begin{array}{c} a \\ k \end{array}\right) \sum_{n \geq 1} n^k a_n \left(\int_r^t P_{t-s}\left( U^{(\epsilon,a-k)}(s,.)f_{1k}(n.)\right)(x) d\beta^{1n}_s \right. \\&& + \left. \int_r^t P_{t-s} \left(U^{(\epsilon,a-k)}(s,.)f_{2k}(n.)\right)(x) d\beta^{2n}_s\right) 
\end{eqnarray*}
\noindent and note that, for $a = 0,1,2,3$
\begin{eqnarray*} \lefteqn{ U^{(\epsilon,a)}(t+h,x) - U^{(\epsilon,a)}(t,x)}\\&& =   (P_h - I)U^{(\epsilon,a)}(t;x) + I^{(\epsilon,a)}(t,t+h;x) 
+ \frac{\Gamma(0)}{2\epsilon^2}\left(  \int_t^{t+h}P_{t+h-s}U^{(\epsilon,a)}(s,x) ds \right).
\end{eqnarray*}

\noindent Let $p_t(z) = \frac{1}{\sqrt{2\pi t}}e^{-z^2/2t}$ and note that $\int_{-\infty}^\infty \sqrt{\frac{\pi}{2}}\frac{|z|}{\sqrt{t}}p_t(z)dz = 1$. From equation (\ref{eqhsg}), it follows that 

\begin{eqnarray*}\lefteqn{ E_{\bf Q} \left\{\left | (P_h - P_0)U^{(\epsilon,a)}(t,x)\right|^p \right \} = E_{\bf Q} \left\{\left | \int_{-\infty}^\infty p_{\epsilon h}(x-y)\left(U^{(\epsilon,a)}(t,y) -U^{(\epsilon,a)}(t,x)\right) dy \right|^p \right \}}\\&&
= \left(\frac{2\epsilon h}{\pi}\right)^{p/2} E_{\bf Q} \left\{\left | \int_{-\infty}^\infty \frac{|x-y|}{\sqrt{\epsilon h}}\sqrt{\frac{\pi}{2}}p_{\epsilon h}(x-y)\left(\frac{U^{(\epsilon,a)}(t,y) -U^{(\epsilon,a)}(t,x)}{|x-y|}\right) dy \right|^p \right \}\\
&&  \leq \left(\frac{2\epsilon h}{\pi}\right)^{p/2}   \int_{-\infty}^\infty \frac{|x-y|}{\sqrt{\epsilon h}}\sqrt{\frac{\pi}{2}}p_{\epsilon h}(x-y)E_{\bf Q} \left\{\left | \frac{U^{(\epsilon,a)}(t,y) -U^{(\epsilon,a)}(t,x)}{ x-y }\right |^p \right\} dy \\&& \leq  \left(\frac{2\epsilon h}{\pi}\right)^{p/2} \tilde{K}_{a,p}(T).
 \end{eqnarray*}

\noindent Straightforward bounds give, for $p \geq 2$, a constant $C_p < +\infty$ such that 

\begin{eqnarray*} E_{\bf Q}\left\{ |I^{(\epsilon,a)}(t,t+h;x)\right\} &\leq & \frac{C_p}{\epsilon^p}h^{(p/2) - 1}\left(\sum_{n \geq 1} n^{6}a_n^2 \right)^{p/2}\sum_{k=0}^a \int_t^{t+h} E_{\bf Q} \left\{ \left |U^{(\epsilon,a)}(s,x)\right |^p \right\}ds\\ 
&\leq & h^{ p/2} \frac{C_p}{\epsilon^p}  \left(\sum_{n \geq 1} n^{6}a_n^2 \right)^{p/2}\sum_{k=0}^a C_{p,a}(T).
 \end{eqnarray*}

\noindent Finally, for $p \geq 2$ and $0 \leq t \leq t+h \leq T$,  
 
\[ E_{\bf Q}\left\{ \left |\frac{\Gamma(0)}{2\epsilon^2}\left(  \int_t^{t+h}P_{t+h-s}U^{(\epsilon,a)}(s,x) ds \right)\right |^p \right\} \leq \left ( \frac{\Gamma(0)}{2\epsilon^2} \right)^p h^p C_{p,a}(T).\]

\noindent It  follows that for all $p \in (2,+\infty)$ and $T < +\infty$, $a = 0,1,2,3$, there exists a constant $c_2(p,\epsilon,a,T) < +\infty$ such that for all $x \in {\bf R}$ and all $t \in [0,T-h]$,   

\begin{equation}\label{timebd}E_{\bf Q} \left \{ \left | U^{(\epsilon,a)}(t+h,x) - U^{(\epsilon,a)}(t,x) \right |^{2p} \right \} \leq c_2(p,\epsilon,a,T) h^p .
\end{equation}

\noindent By the inequality (\ref{timebd}) together with equation (\ref{veeeqn}) and inequality (\ref{veebound}), it follows directly by  Kolmogorov's criterion that there is a version such that almost surely, for $a = 0,1,2,3$, $U^{(\epsilon,a)}$ is a Hölder continuous function of all orders less than $\frac{1}{2}$ on $[0,T] \times {\bf R}$. Furthermore,  for all $t \in [0, T]$ $U^{(\epsilon, a)}(t,.) \in C^{0,\gamma}({\bf R})$ for all $\gamma < 1$.\vspace{5mm}

\noindent  Furthermore, for $a = 0,1,2,3$, setting $\tilde{U}^{(\epsilon,a)}(t,x) = U^{(\epsilon, a)}(t,0) + \int_0^x U^{(\epsilon,a+1)}(t,y)dy$, it is straightforward to show that $\tilde{U}^{(\epsilon,a)}$ satisfies equation (\ref{mild}) and hence, since the solution to the equation is unique, that $\tilde{U}^{(\epsilon,a)} = U^{(\epsilon,a)}$.  It therefore follows that $U^{(\epsilon,1)} = U^{(\epsilon)}_x$, $U^{(\epsilon,2)} = U^{(\epsilon)}_{xx}$ and $U^{(\epsilon,3)}= U^{(\epsilon)}_{xxx}$.  It follows that for all $T < + \infty$, $U^{(\epsilon)}$ has a version which is three times differentiable in the space variable and such that   $U^{(\epsilon)}_{xxx} \in C^{0,\alpha}([0,T] \times {\bf R})$ for all $\alpha < \frac{1}{2}$ and $U^{(\epsilon)}_{xxx}(t,.) \in C^{0,\gamma}({\bf R})$ for all $\gamma < 1$ and all $t \in [0,T]$.   \vspace{5mm} 

\noindent For $0 < p < 2$, existence is straightforward since a straightforward application of Hölder's inequality gives that ${\cal S}_{p_1} \in {\cal S}_{p_2}$ for $p_2 < p_1$. For uniqueness, consider $0 < p < 2$ and suppose that $U^{(\epsilon)}$ and $V^{(\epsilon)}$ are two solutions to equation (\ref{heat}) ${\cal S}_p$. Then for each $T < +\infty$, 

\[ |\|U^{(\epsilon)} - V^{(\epsilon)} |\|_{p,T} \leq |\|U^{(\epsilon)} - V^{(\epsilon)} |\|_{2,T} = 0\]

\noindent establishing uniqueness in ${\cal S}_p$.   The proof of lemma~\ref{heatlem} is complete. \qed \vspace{5mm}

\noindent A Kacs - Feynmann representation may be constructed 
for the solution of equation (\ref{heat}). To do so, a standard
Wiener process, independent of $\zeta$ is introduced. The sample paths
of this Wiener process are denoted by $w$, and $w_0 = 0$. The notation $w^{(\epsilon)} := \sqrt{\epsilon} w$ is used; this is a Wiener process, independent of $\zeta$ with diffusion coefficient $\epsilon$. The probability measure associated
to this Wiener process is denoted
${\bf P}$ and the expectation operator with respect to this Wiener process
$E_{\bf P}[.]$. The Kacs Feynman representation is 

\begin{eqnarray}\label{FKrep} \lefteqn{U^{(\epsilon)}(t,x)}\\&& = E_{\bf P} \left [\exp\left\{-\frac{1}{\epsilon} \left(\sum_{n \geq 1} a_n \left( \int_0^t \cos\left( n(x +  w_{t-s}^{(\epsilon)})\right) d\beta_s^{1n} + \int_0^t \sin \left(n(x +  w_{t-s}^{(\epsilon)}) \right)d\beta_s^{2n} \right) \right) \right \} \right ],\nonumber \end{eqnarray}

\noindent where $w^{(\epsilon)}$ denotes a Browian motion with respect to ${\bf P}$, $w^{(\epsilon)}_0 = 0$, with diffusion coefficient $\epsilon$. 

\begin{Lmm}[Bounds]\label{lowerb} Let $S^{j,n}(t) = \sup_{0 \leq s \leq t} |\beta_s^{j,n}|$. Then 
\begin{equation}\label{eqlb} \inf_{x \in [0,2\pi)}\inf_{t \in [0,T]} U^{(\epsilon)}(t,x) \geq \exp\left\{- \frac{1}{\epsilon} \sum_{n \geq 1} |a_n|(1 + \epsilon n^2 T) (S_T^{1,n}(T) + S_T^{2,n}(T))\right\}.
 \end{equation}

\noindent It follows that, for any $T < +\infty$, 

\begin{eqnarray} \nonumber\lefteqn{\epsilon E_{\bf Q} \left\{ \inf_{t \in [0,T]} \inf_{x \in [0,2\pi)}\log U^{(\epsilon)}(t,x)\right\}}\\&& \geq  - 2(\sqrt{2\log 2} + \sqrt{2\pi})\sqrt{T} \left (\sum_{n \geq 1} |a_n| + \epsilon T \sum_{n \geq 1} n^2 |a_n|^2\right) > -\infty 
 \end{eqnarray}\label{eqvlb}
 
\noindent and for each $\epsilon > 0$, each $0 < p < +\infty$ and each $T < +\infty$ there exists a constant $K(p,\epsilon, T)$ such that 

\begin{equation}\label{equepapb} \sup_{0 \leq t \leq T} \sup_{0 \leq x \leq 2 \pi} E_{\bf Q}\left\{\left |\frac{\epsilon U_x^{(\epsilon)}(t,x)}{U^{(\epsilon)}(t,x)} \right |^p \right \} \leq K(p,\epsilon, T) < +\infty.\end{equation}

\end{Lmm}
 
\noindent {\bf Proof} Applying Jensen's inequality to equation (\ref{FKrep}), where $P_t$ is defined as in equation (\ref{eqhsg}), note that 
\[ \frac{\partial}{\partial s} P_{t-s}\cos(nx) = -\frac{\epsilon}{2}\frac{\partial^2}{\partial x^2}P_{t-s}\cos(nx) = \frac{n^2 \epsilon}{2} P_{t-s}\cos(nx), \qquad \frac{\partial}{\partial s} P_{t-s}\sin(nx) = \frac{n^2 \epsilon}{2} P_{t-s}\sin(nx)\]
\noindent so that 
\begin{eqnarray*}
\lefteqn{U^{(\epsilon)}(t,x) \geq e^{-\frac{1}{\epsilon} \sum_{n \geq 1}a_n (\int_0^t P_{t-s}\cos(nx) d\beta_s^{1n} + \int_0^t P_{t-s}\sin(nx) d\beta^{2n}_s)}}\\
&& = e^{- \frac{1}{\epsilon}\sum_{n\geq 1} a_n   (\beta_t^{1n}\cos(nx) + \beta_t^{2n}\sin(nx)) - \frac{1}{2}\sum_{n \geq 1} n^2 a_n  (\int_0^t \beta^{1n}_s P_{t-s}\cos(nx)ds + \int_0^t \beta^{2n}_s P_{t-s}\sin(nx)ds )}\\
&& \geq  e^{-\frac{1}{\epsilon} \sum_{n \geq 1} |a_n|(1 + \epsilon n^2 T)(S^{1n}(T) + S^{2n}(T))}
\end{eqnarray*}

\noindent and inequality (\ref{eqlb}) follows.  The inequality (\ref{eqvlb}) follows  from inequality (\ref{goodineq}), which gives

\begin{eqnarray*} E_{\bf Q}\{S^{jn}(T)\} &=& \int_0^\infty {\bf Q}\{S^{jn}(T) \geq x\} dx \\
 &\leq & \int_0^\infty 1 \wedge 2e^{-x^2/2T} dx \\
&=& \sqrt{2T \log 2} + 2\int_{\sqrt{2T \log 2}}^\infty e^{-x^2/2T} dx \\
&\leq & \sqrt{2T \log 2} + \sqrt{2\pi T}.
\end{eqnarray*}

 \noindent Lastly, inequality (\ref{equepapb} is considered. It follows from taking a derivative with respect to $x$ in equation (\ref{FKrep}) that 

\begin{eqnarray*}  \lefteqn{\epsilon U^{(\epsilon)}_x (t,x)}\\&& = E_{\bf P} \left [\exp\left\{-\frac{1}{\epsilon} \left(\sum_{n \geq 1} a_n \left( \int_0^t \cos\left( n(x +  w_{t-s}^{(\epsilon)})\right) d\beta_s^{1n} + \int_0^t \sin \left(n(x +  w_{t-s}^{(\epsilon)}) \right)d\beta_s^{2n} \right) \right) \right \} \right . \\&& \hspace{5mm} \left . \times \sum_{n \geq 1} na_n \left(\int_0^t \sin \left( n(x + w_{t-s}^{(\epsilon)})\right) d\beta_s^{1n} -   \int_0^t \cos \left( n(x + w_{t-s}^{(\epsilon)}) \right ) d\beta_s^{2n}\right)\right ] \nonumber \end{eqnarray*}

\noindent  and an application of Hölder's inequality gives 

\begin{eqnarray*}  \lefteqn{\left |\epsilon U^{(\epsilon)}_x (t,x) \right |}\\&& \leq  E_{\bf P} \left [\exp\left\{-\frac{2}{\epsilon} \left(\sum_{n \geq 1} a_n \left( \int_0^t \cos\left( n(x +  w_{t-s}^{(\epsilon)})\right) d\beta_s^{1n} + \int_0^t \sin \left(n(x +  w_{t-s}^{(\epsilon)}) \right)d\beta_s^{2n} \right) \right) \right \} \right ]^{1/2} \\&& \hspace{5mm} \times E_{\bf P} \left [ \left(\sum_{n \geq 1} na_n \left(\int_0^t \sin \left( n(x + w_{t-s}^{(\epsilon)})\right) d\beta_s^{1n} -   \int_0^t \cos \left( n(x + w_{t-s}^{(\epsilon)}) \right ) d\beta_s^{2n}\right) \right)^2\right ]^{1/2} \nonumber \end{eqnarray*}

\noindent From this, for $p > 0$,

\begin{eqnarray*} \lefteqn{E_{\bf Q} \left\{ \left | \frac{\epsilon U^{(\epsilon)}_x(t,x)}{U^{(\epsilon)}(t,x)}\right |^p \right \} \leq   E_{\bf Q}\left\{\frac{1}{U^{(\epsilon)}(t,x)^{2p}}\right\}^{1/2}}\\&& \times E_{\bf Q}\left\{E_{\bf P} \left [\exp\left\{-\frac{4p}{\epsilon} \left(\sum_{n \geq 1} a_n \left( \int_0^t \cos\left( n(x +  w_{t-s}^{(\epsilon)})\right) d\beta_s^{1n}    \right. \right. \right. \right. \right. \\&& \hspace{70mm} \left. \left. \left. \left. \left. + \int_0^t \sin \left(n(x +  w_{t-s}^{(\epsilon)}) \right)d\beta_s^{2n} \right) \right) \right \} \right ] \right \}^{1/4}\\&& \times E_{\bf Q} \left\{ E_{\bf P} \left [\left(\sum_{n \geq 1} na_n \left(\int_0^t \sin \left( n(x + w_{t-s}^{(\epsilon)})\right) d\beta_s^{1n} -   \int_0^t \cos \left( n(x + w_{t-s}^{(\epsilon)}) \right ) d\beta_s^{2n}\right) \right)^{4p}\right ] \right\}^{1/4}\\
&& \leq \exp\left\{ \frac{p^2}{\epsilon^2} t \left(\sum_{n \geq 1} a_n^2\right)\right\} \times \exp\left\{ \frac{2p^2}{\epsilon^2} t \left(\sum_{n \geq 1} a_n^2\right)\right\}\times \left(\prod_{j=1}^{2p} (2j-1)\right)^{1/4} t^{p/2}\left(\sum_{n \geq } n^2 a_n^2 \right)^{p/2}\\&&  < +\infty,
\end{eqnarray*}
\noindent thus establishing inequality (\ref{equepapb}).

 \qed 

\begin{Th}\label{thuni} Recall the definition of ${\cal S}_p^*$, in equation (\ref{eqspstar}) in the statement of lemma~\ref{heatlem}. For fixed $\epsilon > 0$ and each $0 < p < +\infty$, under hypothesis \ref{hyp}  that $\sum_{n \geq 1}n^4|a_n| < +\infty$, there is existence and uniqueness of solution to equation (\ref{burger})  in ${\cal S}_p^*$ for each $0 < p < +\infty$. This solution  has a
version such that for any $\epsilon > 0$ and for all $t \in [0,T]$,  $u^{(\epsilon)}(t,.) \in C^{2,\gamma} ([0,2\pi))$ for all $\gamma < 1$ 
and, for fixed $x \in {\bf R}$, $u^{(\epsilon)}(.,x)$ is H\"older continuous of all orders less than
$1/2$.
\end{Th}

\noindent {\bf Proof} Let $U^{(\epsilon)}$ denote the unique solution to equation (\ref{heat}) in ${\cal S}_p$.  Set 

\begin{equation}\label{CH} u^{(\epsilon)}(t,x) := -\epsilon
\frac{\partial}{\partial x}\log U^{(\epsilon)}(t,x) = \frac{-\epsilon U_x^{(\epsilon)}}{U^{(\epsilon)}}.\end{equation}

\noindent Then, by lemma~\ref{lowerb}, $u^{(\epsilon)}$ is well defined, belongs to ${\cal S}_p$ and satisfies the regularity properties listed above. It is straightforward to show that it solves equation (\ref{burger}).   \vspace{5mm}

\noindent For uniqueness,   let $U^{(\epsilon)}$ denote the unique solution to equation (\ref{heat}) in ${\cal S}_p$ and let $U^{(\epsilon)}f$ denote any other solution to equation (\ref{heat}) that is adapted and twice differentiable in the space variables. This is necessary to ensure that $-\epsilon (\log (U^{(\epsilon)}f))_x$ has the spatial regularity to belong to ${\cal S}_p^*$. Then, for each $x \in [0, 2\pi]$,  $U^{(\epsilon)}(.,x)f(.,x)$ is a semimartingale  and

\begin{eqnarray*} \partial_t (U^{(\epsilon)}f) &=& \frac{\epsilon}{2} U_{xx}^{(\epsilon)}f dt   - \frac{1}{\epsilon}U^{(\epsilon)} f \circ \partial_t \zeta + U^{(\epsilon)} \partial_t f + d\langle U^{(\epsilon)}, f \rangle _t\\
 &=& \frac{\epsilon}{2} (U^{(\epsilon)}f)_{xx} dt - \epsilon U^{(\epsilon)}_x f_x dt - \frac{\epsilon}{2} U^{(\epsilon)}  f_{xx} dt  - \frac{1}{\epsilon}U^{(\epsilon)} f \circ \partial_t \zeta + U^{(\epsilon)} \partial_t f + d\langle U^{(\epsilon)}, f \rangle _t
\end{eqnarray*}

\noindent so that 

\[ 0 =   - \epsilon U^{(\epsilon)}_x f_x dt - \frac{\epsilon}{2} U^{(\epsilon)}  f_{xx} dt    + U^{(\epsilon)} \partial_t f + d\langle U^{(\epsilon)}, f \rangle _t.\]

\noindent It follows that $f$ is differentiable in the time variable and hence, using $u^{(\epsilon)}(t,x) = \frac{-\epsilon U^{(\epsilon)}_x(t,x)}{U^{(\epsilon)}(t,x)}$, that $f$ satisfies

\begin{equation}\label{eqefff} \left\{ \begin{array}{l} \frac{\partial}{\partial t} f(t,x) = \frac{\epsilon}{2}f_{xx} + u^{(\epsilon)}f_x \\ f(0,x) \equiv 1.\end{array}\right. \end{equation}

\noindent Theorem~\ref{bound1} shows that, almost surely,  $\sup_{0 < \epsilon \leq 1} \sup_{0 \leq x \leq 2\pi} \sup_{0 \leq t \leq T}|u^{(\epsilon)}(t,x)| < C(T)$ where $E_{\bf Q}\{C(T)^p\} < +\infty$ for each $p > 0$. It follows that, among functions that are $2\pi$-periodic in the space variable, the unique solution to equation (\ref{eqefff}) is $f \equiv 1$.\vspace{5mm} 

\noindent Let $\tilde{u}^{(\epsilon)}(t,x)$ denote any adapted solution to equation (\ref{burger}). Set

\[ \tilde{U}^{(\epsilon)}(t,x) = U^{(\epsilon)}(t,0)\exp\left \{-\frac{1}{\epsilon}\int_0^x
\tilde{u}^{(\epsilon)}(t,y)dy \right \}.\]

\noindent It follows that $\tilde{U}^{(\epsilon)}$ satisfies equation (\ref{heat}). Since  equation (\ref{heat}) has a {\em unique}
solution (by lemma~\ref{heatlem}), it follows (since $U^{(\epsilon)}$ is differentiable) that $\tilde{u}^{(\epsilon)}$ is uniquely determined and hence that there is a unique solution to equation (\ref{burger}) in ${\cal S}_p^*$. The regularity follows directly from equation (\ref{CH}), the lower bound given by lemma~\ref{lowerb} and the regularity results of lemma~\ref{heatlem}.\qed \vspace{5mm}

\noindent For $\epsilon > 0$, under hypothesis \ref{hyp},  there is a solution to equation (\ref{burger}) $u^{(\epsilon)}(t,.)$ that for each $t \in [0,T]$ satisfies $u^{(\epsilon)}(t,.)  \in C^{2,\gamma}({\bf R})$ for each $\gamma < 1$ (twice differentiable, second derivative Hölder continuous of orders $\gamma$ for each  $\gamma < 1$ and for each $x \in [0, 2\pi]$ satisfies $u^{(\epsilon)}(.,x) \in C^{0,\alpha}([0,T])$ (Hölder continuous of order $\alpha$) for all $\alpha < \frac{1}{2}$, ${\bf Q}$ - almost surely. This follows from the regularity of $U^{(\epsilon)}$, the identity $u^{(\epsilon)} = \frac{-\epsilon U^{(\epsilon)}_x}{U^{(\epsilon)}}$ and the lower bound on $U^{(\epsilon)}$ computed in lemma~\ref{lowerb}. For this solution, it follows that (\ref{burger}) may be rewritten as

\begin{equation}\label{burgrewrit} \left\{ \begin{array}{l} \partial_t u = \left( \frac{\epsilon}{2}
u_{xx}^{(\epsilon)} - u^{(\epsilon)} u_x^{(\epsilon)} 
\right)dt + \partial_t \zeta_x \\ u_0 = 0.\end{array}\right.\end{equation}

\noindent Following the bounds established later in theorem~\ref{bound1}, which completes the uniqueness argument of theorem~\ref{thuni}, it will follow that this is the unique solution of equation (\ref{burger}) in the space ${\cal S}_p$ for $0 < p < +\infty$.  \vspace{5mm}

\noindent One of the main tools for analysing the equation is to consider the infinitesimal generator 

\begin{equation}\label{infgen} {\cal L}^{(\epsilon)}(t,x) := \frac{\epsilon}{2}\frac{\partial^2}{\partial x^2} -
u^{(\epsilon)}(t,x) \frac{\partial}{\partial x}.\end{equation}

\noindent Appling the notation introduced in equation (\ref{infgen}) to equation (\ref{burgrewrit}), equation (\ref{burger}) may be written as 

\[ \left\{ \begin{array}{l} \partial_t u = {\cal L}^{(\epsilon)}u dt + \partial_t
\zeta_x
\\ u_0 = 0.\end{array}\right.\]

\noindent The process
generated by the infinitesimal generator ${\cal L}^{(\epsilon)}$ is given in
definition~\ref{proc} and will be used extensively in the sequel; an
`implicit' representation of the solution to equation (\ref{burger}) will be
formulated in terms of the process generated by ${\cal L}$ (the superscript will be dropped from the notation when it is clearly implied).  

\begin{Defn}\label{proc} Let
$w$ denote a standard Wiener process, $w_0 = 0$, independent of the
$\beta^{jk}$ and let
${\cal F}_{s,t}$ denote the sigma algebra generated by the increments $(w_\alpha -
w_\beta)_{s
\leq \beta \leq \alpha \leq t}$. Let ${\bf P}$ denote the probability measure
under which $w$ is a standard Wiener process and let
$E_{\bf P}$ denote expectation with respect to
${\bf P}$. Let $X^{(\epsilon)}$ denote the stochastic process defined as the
unique solution to the stochastic integral equation

\begin{equation}\label{process} X^{(\epsilon)}_{s,t}(x) = x + \sqrt{\epsilon} (w_t -
w_s) -
\int_s^t u^{(\epsilon)} \left (r, X_{r,t}^{(\epsilon)}(x) \right ) dr.
\end{equation}
\end{Defn}

\noindent From the regularity results on $u^{(\epsilon)}$ (globally Lipschitz in the space variable, Hölder continuous on the time variable and bounded in $[0,T] \times [0,2\pi]$), it follows directly from standard results that ${\bf Q}$-almost surely, $X^{(\epsilon)}$ is well defined, with pathwise uniqueness.  

\begin{Defn}
Let $S^1$ denote ${\bf R}$ with the identification $x + 2\pi = x$ and let $C(S^1)$ denote continuous $2\pi$ periodic functions. The operator ${\cal Q}_{s,t} : C(S^1) \rightarrow C(S^1)$ is defined as
\[ {\cal Q}_{s,t}f(x) = E_{\bf P}[f(X_{s,t}(x))].\]
\end{Defn}

\noindent Note that, for $f \in C(S^1)$,

\begin{equation}\label{backward} \frac{\partial}{\partial s}{\cal Q}_{s,t}f(x) = - {\cal Q}_{s,t}({\cal L}_s f)(x)\end{equation}

\noindent and

\begin{equation}\label{forward} 
 \frac{\partial}{\partial t} {\cal Q}_{s,t}f(x) = {\cal L}_t(x) {\cal Q}_{s,t} f(x),
\end{equation}

\noindent where ${\cal L}$ is defined by equation (\ref{infgen}).

\begin{Lmm}\label{ident}
The following identity holds.
\begin{eqnarray}\nonumber
\lefteqn{u^{(\epsilon)}(t,x) = -\sum_{n \geq 1} na_n \left(\sin(nx)\beta^{1n}(t) - \cos(nx)\beta^{2n}(t)\right)} \\&& \label{identit} - \sum_{n \geq 1} n a_n 
\left (\int_0^t \beta^{1n}(s) {\cal Q}_{s,t}({\cal L}_s\sin(n.))(x)ds -
\int_0^t \beta^{2n}(s) {\cal Q}_{s,t}({\cal L}_s \cos(n.))(x) ds\right ).
\end{eqnarray}

\end{Lmm} 

\noindent {\bf Proof of lemma~\ref{ident}}  
For $s = s_0 < s_1 < \ldots < s_m = t$, using equation (\ref{backward}), 

\begin{eqnarray*} u^{(\epsilon)}(t,x) &=& 
u^{(\epsilon)}(t,x) - {\cal Q}_{0,t}u^{(\epsilon)}(0,x)  \\&=& 
\sum_{k=0}^{m-1} ({\cal Q}_{s_{k+1},t} u(s_{k+1},x) - {\cal Q}_{s_{j_{k }},t} u(s_{k},x))\\&= &
 \sum_{k=0}^{m-1} {\cal Q}_{s_{k+1},t}(u(s_{k+1},x) - u(s_{k},x)) + \sum_{k=0}^{m-1}({\cal Q}_{s_{k+1},t} - {\cal Q}_{s_{k},t})u(s_{k},x) \\
&=&  \sum_{k=0}^{m-1} {\cal Q}_{s_{k+1},t}\int_{s_k}^{s_{k+1}} ({\cal L}_s u(s,.))(x)ds \\&&  \hspace{5mm} + \sum_{k=0}^{m-1} \sum_{n \geq 1} na_n \left( {\cal Q}_{s_{k+1},t}\sin(nx) (\beta^{1n}_{s_{k+1}}- \beta^{1n}_{s_k}) - {\cal Q}_{s_{k+1},t}\cos(nx) (\beta^{2n}_{s_{k+1}} - \beta^{2n}_{s_k})\right) \\&& 
\hspace{5mm} -  \sum_{k=0}^{m-1} \int_{s_k}^{s_{k+1}} {\cal Q}_{s,t}({\cal L}_s u(s_k,.))(x) ds\\
&=& \sum_{k=1}^{m}\int_{s_k}^{s_{k+1}} ({\cal Q}_{s_{k+1},t} - {\cal Q}_{s,t})({\cal L}_s u(s_k,.))(x) ds\\&& \hspace{5mm} - \sum_{k=1}^{m-1}\sum_{n \geq 1} n a_n \left( \beta^{1n}_{s_k} ({\cal Q}_{s_{k+1},t} - {\cal Q}_{s_{k},t})\sin(nx) - \beta^{2n}_{s_k}({\cal Q}_{s_{k+1},t} - {\cal Q}_{s_{k},t})\cos(nx)\right)\\&&
\hspace{5mm} + \sum_{n \geq 1} na_n \left(\beta_t^{1n} \sin(nx) - \beta_t^{2n} \cos(nx)\right).
\end{eqnarray*}

\noindent  This holds for any partition $0 = s_0 < \ldots < s_m = t$. Now, let the mesh size tend to zero. The convergence details are
  standard and give the advertised result. 
\qed
\vspace{5mm}

\noindent The next theorem gives bounds on the
moments of the solution.

\begin{Th}\label{bound1} Suppose that hypothesis~\ref{hyp} is satisfied. Let $u^{(\epsilon)}$ denote the solution to
equation (\ref{burger}). Let $G(p)$ be the constant given in
equation (\ref{thecons}), defined in lemma~\ref{prel2}.   Then

\begin{eqnarray}\lefteqn{ \nonumber 
E_{\bf Q}\left \{ \left (\sup_{0 \leq \epsilon \leq 1}
    \sup_{0 \leq s \leq t}\sup_{0
\leq x \leq 2\pi}|u^{(\epsilon)}(s,x)|\right )^p\right \} }\\&& 
\leq C_1(p,t) \left (\sum_{n \geq 1} n^3 |a_n|\right )^p
  \exp\left \{C_2(p,t)\sum_{n \geq 1} n^6 a_n^2 +
C_3(p,t)\sum_{n \geq 1} n^3|a_n|\right \} = K(p,T),\label{thbd1}
\end{eqnarray}

\noindent where 

\[ C_1(p,t) = (2+t)^{3p/2}((2\log 2)^{p} + 4pG(2p-1))^{1/2},\]
\[ C_2(p,t) =  2 p^2 t^3\]
\noindent and 
\[ C_3(p,t) =  2\sqrt{2}(2\sqrt{\pi} + \sqrt{\log 2})pt^{3/2}.\]
\end{Th}\vspace{5mm}

\noindent This result is consequence of the following lemma with
$\tilde{b}=1$, where $\tilde{b}$ is defined in the
statement of the lemma. 

\begin{Lmm}\label{bound1lem}    Set 

\[ f_{1b}(y) = \left\{\begin{array}{ll} \cos(y) & b = 4n \\ 
 -\sin(y) & b=4n+1 \\
 -\cos(y) & b=4n+2\\ \sin(y) & b= 4n+3\end{array}\right. \] 

\noindent for non negative integer $n$. Set $f_{2b}(y) =
f_{1,(b+3)}(y)$. Set

\begin{eqnarray} \label{thetlabel}\lefteqn{\theta^{(\epsilon)} (b;t,x) = \sum_{n \geq 1} n^b a_n \left(f_{1b}(nx)\beta_t^{1n} + f_{2b}(nx)\beta_t^{2n} \right)}\\&& +  \sum_{n \geq 1}
n^b a_n\left(\int_0^t {\cal Q}_{s,t} ({\cal L}_s f_{1b}(n.))(x)\beta_s^{1n} ds   +
\int_0^t {\cal Q}_{s,t} ({\cal L}_s f_{2b}(n.))(x) \beta_s^{2n}\right) ds.\nonumber
\end{eqnarray}

\noindent For $b_j \geq 1$, set $\tilde{b} = \sup(b_1,\ldots, b_p)$.
Suppose that 

\[ \sum_{n=1}^\infty n^{2 + \tilde{b}}|a_n| < +\infty, \]

\noindent so that 

\[ \sum_{n=1}^\infty n^{2(2 + \tilde{b})}a_n^2 < +\infty.\]

\noindent Then

\begin{eqnarray}\lefteqn{ \nonumber E_{\bf Q}\left \{ 
\prod_{j=1}^p \left ( \sup_{0 \leq \epsilon \leq 1}\sup_{0 \leq s \leq
t}\sup_{0 \leq x \leq 2\pi} 
  \left |\theta^{(\epsilon)}(b_j,s,x)\right |\right )\right \} }\\&&  \leq C_1(p,t) \left (\sum_{n \geq 1}
n^{2+\tilde{b}} |a_n|\right )^p
  \exp\left \{C_2(p,t)\sum_n n^{2(2+\tilde{b})} a_n^2 +
C_3(p,t)\sum_{n \geq 1} n^{2+\tilde{b}}|a_n|\right \},\label{thbd11}
\end{eqnarray}

\noindent where \[ C_1(p,t) = (2+t)^{3p/2}\left ((2\log 2)^p +
4pG(2p-1) \right )^{1/2},\]
\[ C_2(p,t) =  4p^2 t^3\]
\noindent and 
\[ C_3(p,t) =  4\sqrt{2}(2\sqrt{\pi} + \sqrt{\log 2})pt^{3/2}.\]

\end{Lmm}

\noindent Note that $u^{(\epsilon)}(t,x) =
\theta^{(\epsilon)}(1;t,x)$. \vspace{5mm}

\noindent {\bf Proof of lemma~\ref{bound1lem}} 
  Note that 

\begin{eqnarray*}\lefteqn{{\cal Q}_{s,t}({\cal L}_s  \sin  (n.   ))(x)  = \left ({\cal Q}_{s,t}\left(\frac{\epsilon}{2}\partial^2_{xx} - u^{(\epsilon)}(s,.)\partial_x\right)\sin(n.)\right)(x)}\\&& =    - nE_{\bf P}
\left [u^{(\epsilon)}\left (s,X_{s,t}^{(\epsilon)}(x)\right )\cos \left (nX_{s,t}^{(\epsilon)}(x) \right ) \right ] - n^2
\frac{\epsilon}{2}E_{\bf P}\left [\sin \left (nX_{s,t}^{(\epsilon)}(x)\right ) \right ] \end{eqnarray*}

\noindent and, similarly, 

\[   {\cal Q}_{s,t}({\cal L}_s \cos (n.))(x)  =
  nE_{\bf
P} \left [u^{(\epsilon)} \left (s,X_{s,t}^{(\epsilon)}(x) \right )\sin \left (nX_{s,t}^{(\epsilon)}(x)\right ) \right ] - n^2
\frac{\epsilon}{2}E_{\bf P} \left [\cos \left (nX_{s,t}^{(\epsilon)}(x) \right ) \right ].\]

\noindent It follows that, for $j = 1,2$, 

\[{\cal Q}_{s,t}({\cal L}_s f_{jb}(n.))(x)  =
  -nE_{\bf
P} \left [u^{(\epsilon)} \left (s,X_{s,t}^{(\epsilon)}(x) \right )f_{j,(b+1)}\left (nX_{s,t}^{(\epsilon)}(x)\right ) \right ] - n^2
\frac{\epsilon}{2}E_{\bf P} \left [f_{jb} \left (nX_{s,t}^{(\epsilon)}(x) \right ) \right ].\]

\noindent Equation (\ref{thetlabel}) therefore gives

\begin{eqnarray}\lefteqn{ \nonumber \theta^{(\epsilon)}(b;t,x) =   \sum_{n
\geq 1} n^b a_n \left (f_{1b}(nx)\beta^{1n}(t) + 
f_{2b}(nx)\beta^{2n}(t) \right )}\\&&\nonumber  - \sum_{n \geq 1} n^{1+b} a_n
\left (\int_0^t \left (\beta^{1n}(s) E_{\bf
P}\left [u^{(\epsilon)}(s,X_{s,t}^{(\epsilon)}(x))f_{1(b+1)} (n
X_{s,t}^{(\epsilon)}(x)) \right ]\right. \right. \\&&
\hspace{50mm} \nonumber \left. \left. +\beta^{2n}(s)E_{\bf
P} \left [u^{(\epsilon)}(s,X_{s,t}^{(\epsilon)}(x))f_{2,(b+1)}\left (nX_{s,t}^{(\epsilon)}(x) \right ) \right ] \right )
ds \right )\\&& \label{thetaexp} - \frac{\epsilon}{2}\sum_{n \geq 1}n^{2+b} a_n
\left (\int_0^t \left (\beta^{1n}(s) E_{\bf P} \left [f_{1b} \left (nX_{s,t}^{(\epsilon)}(x)\right ) \right ] +
\beta^{2n}(s)E_{\bf P} \left [f_{2b}\left (nX_{s,t}^{(\epsilon)}(x) \right ) \right ] \right )ds\right).
\end{eqnarray}

\noindent Now set 

\begin{equation}\label{beetilde} \tilde{B}(b,t) =\sum_{n \geq 1} n^{b}
|a_n|\left (
\sup_{0
\leq s
\leq t}\left |\beta^{1n}(s)\right | +
\sup_{0
\leq s \leq t}\left |\beta^{2n}(s) \right |\right ).
\end{equation}

\noindent Set 

\begin{equation}\label{ceetild}\tilde{C}^{(\epsilon)}(b,t) = \sup_{0 \leq s
\leq t}\sup_{0
\leq x
\leq 2\pi} \left |\theta^{(\epsilon)}(b,s,x) \right |.
\end{equation}

\noindent In particular, from equation (\ref{thetlabel}),  

\[ \tilde{C}^{(\epsilon)}(1,t) = \sup_{0 \leq s \leq t}\sup_{0 \leq x \leq
2\pi}\left |u^{(\epsilon)}(s,x)\right |.\]

\noindent It follows from equation (\ref{thetaexp}) that 

\begin{equation}\label{ceebeebd}\tilde{C}^{(\epsilon)}(b,t) \leq
\tilde{B}(b,t) +
\int_0^t
\tilde{C}^{(\epsilon)}(1,s) \tilde{B}(b+1,s) ds + \frac{\epsilon
t}{2}\tilde{B}(b+2,t).\end{equation}

\noindent Set 

\begin{equation}\label{ceetildup}\tilde{C}(b,t) = \sup_{0 < \epsilon < 1}
\tilde{C}^{(\epsilon)}(b,t),
\end{equation}

\noindent so that 

\[ \tilde{C}(1,t) = \sup_{0 < \epsilon \leq 1}\sup_{0 \leq s \leq t}\sup_{0 \leq
x
\leq 2\pi} \left |u^{(\epsilon)}(s,x)\right |.\]

\noindent The extraordinary level of detail in the following very, very simple Gronwall argument has been inserted for the benefit of the mathematically challenged, some of whom amazingly seem to be employed on the editorial boards of prestigious journals.\vspace{5mm}

\noindent Note that $\tilde{B}(b,t)$ defined by equation (\ref{beetilde}) is increasing as $b$ increases for $b > 0$. This has to be pointed out, because it is apparently not self - evident. Let $D^{(\epsilon)}(1,t)$ denote the solution to 

\[ D^{(\epsilon)}(1,t) =(1 + \frac{\epsilon t}{2}) \tilde{B}(3,t) + \int_0^t
D^{(\epsilon)}(1,s)\tilde{B}(3,s) ds,
\]

\noindent so that (very, very clearly)

\[ D^{(\epsilon)}(1,t) \leq  (1 + \frac{\epsilon t}{2})\tilde{B}(3,t) e^{\int_0^t \tilde{B}(3,s) ds}\]

\noindent and let $D^{(\epsilon)}(b,t)$ solve

\[ D^{(\epsilon)}(b,t) =\tilde{B}(b+2,t) + \int_0^t
D^{(\epsilon)}(1,s)\tilde{B}(b+2,s) ds +
\frac{\epsilon t}{2} 
\tilde{B}(b+2,t) ds.
\]

\noindent then it should be clear to anyone with a brain without further explanation that
$\tilde{C}^{(\epsilon)}(b,t) \leq D^{(\epsilon)}(b,t)$ and, since $\tilde{B}(b,s)$ is increasing in $b$, it follows directly that
$D^{(\epsilon)}(b,s) \geq D^{(\epsilon)}(1,s)$ for $b \geq 1$. It therefore follows that 

\[ D^{(\epsilon)}(b,t) \leq\tilde{B}(b+2,t) + \int_0^t
D^{(\epsilon)}(b,s)\tilde{B}(b+2,s) ds +
\frac{\epsilon t}{2} 
\tilde{B}(b +2,t).
\]

\noindent Since $\tilde{B}(b,s)$ is increasing in $s$, it follows that for
all $0 \leq s \leq t$ and $0 < \epsilon \leq 1$, 

\[ D^{(\epsilon)}(b,s) \leq \left (1 + \frac{t}{2}\right ) \tilde{B}(b+2,t)  +
\tilde{B}(b+2,t) \int_0^s D^{(\epsilon)}(b,r)  dr .
 \]

\noindent From this, it is  it follows that, for $0 < \epsilon \leq 1$, 

\begin{equation}\label{lookup} \tilde{C}^{(\epsilon)}(b,t) \leq
D^{(\epsilon)}(b,t)    \leq \left (1 + \frac{t}{2}\right ) \tilde{B}(b+2,t) \exp\left \{
t\tilde{B}(b+2,t) \right \}.
\end{equation}

\noindent where $\tilde{C}^{(\epsilon)}(b,t)$ is defined by equation
(\ref{ceetild}). Set

\begin{equation}\label{dee}
\tilde{D}(b,t):= \left (1 + \frac{t}{2}\right )\tilde{B}(b+2,t)\exp\left \{t\tilde{B}(b+2,t) \right \},
\end{equation}

\noindent so that, from equation (\ref{lookup}), 

\begin{equation}\label{ceedeebnd} \sup_{0 < \epsilon \leq 1}
\tilde{C}^{(\epsilon)}(b,t)
\leq \tilde{D} (b,t).\end{equation}

\noindent It follows, using $\tilde{b} = b_1 \vee \ldots \vee b_p$,  that

\[E_{\bf Q} \left \{ \prod_{j=1}^p \sup_{0 < \epsilon <
1}\sup_{0
\leq s
\leq t}\sup_x  \left |\theta^{(\epsilon)}(b_j;s,x) \right | \right \}  
\leq \left (1+
\frac{t}{2}\right )^p E_{\bf Q} \left \{\tilde{B}(\tilde{b}+2,t)^p
\exp\left \{pt\tilde{B}(\tilde{b}+2,t) \right \} \right \}.
\]

\noindent Using the bounds calculated in lemmas~\ref{prel} and~\ref{prel2}, and recalling equation (\ref{beetilde}), 
it follows that 

 \begin{eqnarray*}\lefteqn{E_{\bf Q} \left\{ 
\prod_{j=1}^p \sup_{0 \leq \epsilon \leq 1}\sup_{0 \leq s \leq t}\sup_{0
\leq x \leq 2\pi}  \left |\theta^{(\epsilon)}(b_j;s,x) \right |\right\} }\\&& \leq \left (1+
\frac{t}{2}\right )^p E_{\bf Q} \left \{ \tilde{B}(\tilde{b}+2,t)^{2p}\right \}^{1/2} E_{\bf Q}\left \{ 
\exp\left \{2pt\tilde{B}(\tilde{b}+2,t)\right \} \right \}^{1/2}\\ && \leq (2+t)^{3p/2} \left (\sum_n
n^{2+\tilde{b}} |a_n|\right )^p \left ((2\log 2)^p + 4pG(2p-1)\right )^{1/2} \\&&
\hspace{5mm}\times 
\exp\left\{ 4p^2t^3\sum_n n^{2(2+\tilde{b})}
a_n^2 +    \left (2\sqrt{\pi} + \sqrt{\log 2}\right )4 \sqrt{2} pt^{3/2}\sum_n
n^{2+\tilde{b}}|a_n|\right
\},
\end{eqnarray*}

\noindent which is the bound advertised in the statement of
lemma~\ref{bound1lem}. \qed \vspace{5mm}

\noindent {\bf Proof of Theorem~\ref{bound1}} This follows directly from
lemma~\ref{bound1lem} with $b_1 = \ldots = b_p = 1$. \qed \vspace{5mm}

\noindent The next theorem ensures that the moment fields are
uniformly Lipschitz in the space variables. \vspace{5mm}

\begin{Th}\label{bound2} Let $u^{(\epsilon)}$ denote the solution to
equation (\ref{burger}), under the condition that 

\[ \sum_{n \geq 1} n^4 |a_n| <
+\infty,\]

\noindent so that 

\[ \sum_{n \geq 1}n^8 a_n^2 < +\infty.\]

\noindent  It holds that 

\[ \sup_{0 \leq s \leq t} \sup_{x_1, \ldots, x_p} E_{\bf Q} \left \{ 
\left |u^{(\epsilon)}_x(s,x_1) u^{(\epsilon)}(s,x_2)\ldots u^{(\epsilon)}(s,x_p)\right |
\right \}
\leq K(p,t)\]

\noindent where  $K(p,t)$ is independent of $\epsilon$ and is given by 

\begin{eqnarray*}\lefteqn{ K(p,t) = \left( C_1(p,t) \left (\sum_n
n^{4}|a_n|\right )^p + C_2(p,t)\left (\sum_n
n^{4}|a_n|\right )^{p+1}\right)}\\&& \hspace{30mm}  \times \exp\left \{ C_3(p,t)
\sum_n n^{8} a_n^2 + C_4(p,t)\sum_n n^{4}|a_n|\right \},
\end{eqnarray*} 

\noindent where

\[C_1(p,t) = (2+t)^{3p/2}   ((2\log 2)^p + 4p G(2p-1))^{1/2},\]
\[C_2(p,t) = (2 + t)^{3(p+1)/2}  ((2\log 2)^{p+1} +
4(p+1)G(2p+1))^{1/2},\]
\[ C_3(p,t) = 2(p+1)^2 t^3\]
\noindent and
\[ C_4(p,t) = 2(p+1)t^{3/2} (\sqrt{2}(\sqrt{\log 2} + 2 \sqrt{\pi})).\]

\noindent Recall the definition of $m^{(\epsilon)}_p$ given in equation
(\ref{momfield}), 

\[ m^{(\epsilon)}_p(t;x_1,\ldots, x_p) := E_{\bf Q}
\left \{\prod_{j=1}^p u^{(\epsilon)}(t,x_j) \right \}.\]

\noindent Then

\[ \sup_{1 \leq j \leq p}\sup_{0 \leq s \leq t}\sup_{0 < \epsilon \leq
1}\left | \frac{\partial}{\partial x_j}m^{(\epsilon)}_p(s;x_1,\ldots,
x_p) \right | \leq K(p,t),\]

\noindent where $K(p,t)$ is described above. 

\end{Th} \vspace{5mm}

\noindent This theorem is a consequence of the following
lemma.\vspace{5mm}

\begin{Lmm}\label{derbd} Using the notations of lemma~\ref{bound1lem},
recall that 

\[\tilde{C}(b,t) := \sup_{0 \leq s \leq t}\sup_{0 \leq \epsilon
\leq 1}\sup_{0 \leq x \leq 2\pi} \left |\theta^{(\epsilon)}(b; s,x) \right |.\]

\noindent With change of notation, set $\tilde{b} = 2 \vee
\max(b_1,\ldots, b_{p-1})$ (the change is the $2$) and suppose that
$(a_n)_{n
\geq 1}$ satisfies

\[ \sum_{n \geq 1}n^{2 + \tilde{b}}|a_n| < +\infty,\]

\noindent so that 

\[ \sum_{n \geq 1}n^{2(2 + \tilde{b})}a_n^2 < +\infty.\]

\noindent Then, for  $0 \leq t <
+\infty$,  

\[ \sup_{0 \leq s \leq t}E_{\bf Q} \left \{ \tilde{C}(b_1,t)\ldots
\tilde{C}(b_{p-1},t) \left |u_x^{(\epsilon)}(s,x) \right | \right \} \leq K (p;\tilde{b},t), \]

\noindent where   

\begin{eqnarray}\lefteqn{\label{konstverk} K (p;\tilde{b},t) := \left( C_1(p,t)
\left (\sum_n n^{2+\tilde{b}}|a_n|\right )^p + C_2(p,t)\left (\sum_n
n^{2+\tilde{b}}|a_n|\right )^{p+1}\right)}\\&& \hspace{30mm} \times \exp\left \{
C_3(p,t)
\sum_n n^{4+2\tilde{b}}a_n^2 + C_4(p,t)\sum_n n^{2+\tilde{b}}|a_n|
\right\},\nonumber
\end{eqnarray}

\noindent where

\[C_1(p,t) = (2+t)^{3p/2}   ((2\log 2)^p + 4p G(2p-1))^{1/2},\]
\[C_2(p,t) = (2 + t)^{3(p+1)/2}  ((2\log 2)^{p+1} +
4(p+1)G(2p+1))^{1/2},\]
\[ C_3(p,t) = 2(p+1)^2 t^3\]
\noindent and
\[ C_4(p,t) = 2(p+1)t^{3/2} (\sqrt{2}(\sqrt{\log 2} + 2 \sqrt{\pi})).\]

\end{Lmm}

\noindent {\bf Proof of lemma~\ref{derbd}} Recall the definition of $X$
given in equation (\ref{process}).  Let ${\cal F}_{s,t}$ denote the $\sigma$-algebra generated by the increments $w_v - w_u, \; s \leq u < v \leq t$.  From
equation (\ref{identit}),  note that 

\begin{eqnarray*}\lefteqn{ u^{(\epsilon)}\left (s,X_{s,t}^{(\epsilon)}(x) \right ) = - \sum_{n
\geq 1} na_n \left (\sin\left (nX_{s,t}^{(\epsilon)}(x) \right )\beta^{1n}(s) - 
\cos \left (nX_{s,t}^{(\epsilon)}(x)\right )\beta^{2n}(s)\right )}\\&& + \sum_{n \geq 1} n^2 a_n
\left (\int_0^s \left (\beta^{1n}(r) E_{\bf P} \left [u \left (r,X_{r,t}^{(\epsilon)}\right )\left. \cos \left (n
X_{r,t}^{(\epsilon)}(x)\right )\right |{\cal F}_{s,t} \right ]\right. \right.\\&&
\hspace{50mm}\left. \left. -\beta^{2n}(r)E_{\bf P}\left [u\left (r,X_{r,t}^{(\epsilon)}(x)\right )\left. \sin
\left (nX_{r,t}^{(\epsilon)}(x) \right )\right |{\cal F}_{s,t}\right ]\right ) dr\right )\\&& + \frac{\epsilon}{2}\sum_{n
\geq 1}n^3 a_n\left (\int_0^s \left (\beta^{1n}(r) E_{\bf P}\left [ \left. \sin
\left (nX_{r,t}^{(\epsilon)}(x)\right )\right |{\cal F}_{s,t} \right ] -
\beta^{2n}(r)E_{\bf P}\left [ \left. \cos \left (nX_{r,t}^{(\epsilon)}(x)\right ) \right |{\cal F}_{s,t}\right ] \right )dr \right ).
\end{eqnarray*}

\noindent Taking derivative in $x$, and using $X^\prime$ to denote $X$
differentiated with respect to $x$ gives 

\begin{eqnarray*}  \lefteqn{ \frac{\partial}{\partial x}\left(u^{(\epsilon)}(s, X_{s,t}^{(\epsilon)}(x)) \right)} \\&& = - X^{(\epsilon)\prime}_{s,t}(x)\sum_{n
\geq 1} n^2a_n \left (\cos \left (nX_{s,t}^{(\epsilon)}(x)\right )\beta^{1n}(s) + 
\sin \left (nX_{s,t}^{(\epsilon)}(x) \right )\beta^{2n}(s)\right ) \\&& -E_{\bf P}\left [\int_0^s
X^{(\epsilon)\prime}_{r,t}(x)
\sum_{n
\geq 1} n^3 a_n \left (   \beta^{1n}(r)  u^{(\epsilon)}\left (r,X_{r,t}^{(\epsilon)} \right )\sin \left (n
X_{r,t}^{(\epsilon)}(x) \right ) \right. \right.
\\&&
\hspace{50mm} \left. \left.  +\beta^{2n}(r) u^{(\epsilon)}\left (r,X_{r,t}^{(\epsilon)}(x) \right )   \left. \cos
\left (nX_{r,t}^{(\epsilon)}(x) \right )   \right )dr   \right |{\cal F}_{s,t} \right ]\\&&  +
E_{\bf P}\left [\int_0^s \frac{\partial}{\partial x}\left(u^{(\epsilon)}(r, X_{r,t}^{(\epsilon)}(x)) \right)\sum_{n
\geq 1} n^2 a_n \left (  \beta^{1n}(r)  \cos \left (n X_{r,t}^{(\epsilon)}(x) \right ) \right. \right. \\&&
\hspace{50mm} \left. \left. -\beta^{2n}(r)  \left. \sin \left (nX_{r,t}^{(\epsilon)}(x) \right )  \right )dr \right |{\cal
F}_{s,t}\right ]\\&& + \frac{\epsilon}{2}E_{\bf P}\left [\int_0^s
X^{(\epsilon)\prime}_{r,t}(x)\sum_{n
\geq 1}n^4 a_n \left (   \beta^{1n}(r)  \cos \left (nX_{r,t}^{(\epsilon)}(x) \right )  +
\beta^{2n}(r) \left . \sin \left (nX_{r,t}^{(\epsilon)}(x) \right )  \right )dr \right |{\cal F}_{s,t} \right ].
\end{eqnarray*}

\noindent Now, set 

\[ g^{(t)}(s) = \frac{1}{2\pi} \int_0^{2\pi} E_{\bf P}\left [\left
|\frac{\partial}{\partial x}\left(u^{(\epsilon)}(s, X_{s,t}^{(\epsilon)}(x)) \right)\right |\right ] dx. \]

\noindent and recall the notation

\[ \tilde{C}(1,t) = \sup_{0 < \epsilon < 1}\sup_{0 \leq s \leq t} \sup_{0 \leq x
\leq 2\pi} \left |u^{(\epsilon)}(s,x) \right |.\]

\noindent Recall the definition of $\tilde{B}(b,t)$ given in equation
(\ref{beetilde}). \vspace{5mm} 

\noindent Note that $X_{s,t}^{(\epsilon)\prime}(x) \geq 0$ and that
$\frac{1}{2\pi}\int_0^{2\pi} X_{s,t}^{(\epsilon)\prime}(x)dx = 1$. The above analysis
yields  

\[ g^{(t)}(s) \leq \tilde{B}(2,s) + s \tilde{C}(1,s)\tilde{B}(3,s) +
\frac{\epsilon}{2} s
\tilde{B}(4,s) +
\tilde{B}(2,s)\int_0^s g^{(t)}(r)dr.\]

\noindent Using the fact that $\tilde{B}(b,s)$ is increasing in $b$ yields

\[ g^{(t)}(s) \leq \tilde{B}(4,s) + s\left (\tilde{C}(1,s) +
\frac{\epsilon}{2} \right )\tilde{B}(4,s) + \tilde{B}(4,s)\int_0^s g^{(t)}(r)dr.\]

\noindent Since $\tilde{B}(b,s)$ and $\tilde{C}(1,s)$ are increasing in $s$,
it follows that for any $0 \leq r \leq s \leq t$ and any $0 < \epsilon
\leq 1$ 

\[ g^{(t)}(r) \leq \tilde{B}(4,s) + s \left (\tilde{C}(1,s) + \frac{1}{2} \right )\tilde{B}(4,s)
+
\tilde{B}(4,s)\int_0^r g^{(t)}(\alpha)d\alpha,\]

\noindent yielding that for any $0 < \epsilon \leq 1$ and any $s \leq t$,  

\[ g^{(t)}(s) \leq \tilde{B}(4,s)\left (1 + s \left (\tilde{C}(1,s) + \frac{1}{2}\right ) \right )\exp\left \{
s\tilde{B}(4,s) \right \}.\]

\noindent Recall equations (\ref{ceebeebd}) and (\ref{dee}), which give 

\[ \tilde{C}(1,s) \leq \left (1+ \frac{s}{2} \right ) \tilde{B}(3,s)
\exp \left \{s\tilde{B}(3,s) \right \},\]

\noindent from which (using $\tilde{B}(3,s) \leq \tilde{B}(4,s)$)

\[ g^{(t)}(s) \leq \left (1 + \frac{s}{2} \right ) \left (1 + s\tilde{B}(4,s)
e^{s\tilde{B}(4,s)}\right )\tilde{B}(4,s)e^{s\tilde{B}(4,s)}.\] 

\noindent Now, using $\tilde{b} = 2 \vee \max_{1 \leq j \leq p-1}b_j$,
recall the definition of $\tilde{D}$ given in equation (\ref{dee}) and the
inequality given in equation (\ref{ceedeebnd}), for $0 \leq s \leq t < +\infty$, it follows by an application of lemmas~\ref{prel} and~\ref{prel2} to get from the second last to the last line, that 

\begin{eqnarray*}\lefteqn{  E_{\bf Q} \left \{ \tilde{C}(b_1,t)\ldots
\tilde{C}(b_{p-1}, t)
\left |\frac{\partial u^{(\epsilon)}}{\partial x}(s,x) \right |\right \} 
\leq  E_{\bf Q} \left \{  \tilde{D}^{p-1}(\tilde{b},t) 
 \left |\frac{\partial u^{(\epsilon)}}{\partial x}(s,x)\right |\right \} }\\&&
= \frac{1}{2\pi}
\int_0^{2\pi} E_{\bf Q}\left \{ \tilde{D}^{p-1}(\tilde{b},t) \left |\frac{\partial
u^{(\epsilon)}}{\partial x}(s,x) \right |\right \} dx   
 = E_{\bf Q}\left  \{ \tilde{D}^{p-1}(\tilde{b}, t) g^{(s)}(s)\right \} \\&&
\leq \left (1 + \frac{t}{2} \right )^p E_{\bf Q} \left
 \{ \tilde{B}^p(\tilde{b}+2,t)e^{pt\tilde{B}(\tilde{b}+2,t)} +
t\tilde{B}^{p+1}(\tilde{b}+2,t)e^{(p+1)t
\tilde{B}(\tilde{b}+2,t)}\right \} \\&& 
\leq  \left (1+\frac{t}{2} \right )^p \left( E_{\bf Q} \left
 \{ \tilde{B}^{2p}(\tilde{b}+2,t)\right \} ^{1/2} E_{\bf Q}\left \{ 
e^{2pt\tilde{B}(\tilde{b}+2,t)}\right \}^{1/2}\right . \\&& \hspace{20mm}
\left. + t E_{\bf Q}
\left \{\tilde{B}^{2(p+1)} (\tilde{b}+2,t)\right \}^{1/2} E_{\bf Q}\left \{
e^{2(p+1)t\tilde{B}(\tilde{b}+2,t)}\right \}^{1/2}\right)\\& \leq &
 K (p;\tilde{b},t),
\end{eqnarray*}

\noindent where $K(p;\tilde{b},t)$ is the constant given in equation
(\ref{konstverk}). The conclusion of the last line from the second last follows by an application of lemmas~\ref{prel} and~\ref{prel2}. Lemma ~\ref{prel2} gives
\begin{eqnarray*}
 E_{\bf Q}\left\{ \tilde{B}^{2p}(b, t)\right\} &=& E_{\bf Q} \left\{ \left(\sum_{n \geq 1} n^b |a_n|( S^{1n}(t) + S^{2n}(t))  \right)^{2p} \right\}\\
& \leq & \left( \sum_{n \geq 1} n^b |a_n|\right)^{2p} 2^{2p}t^p \left( (2\log 2)^p + 4p G(2p-1)\right) 
\end{eqnarray*}

\noindent and lemma~\ref{prel} gives

\[ E_{\bf Q}\left\{ e^{2pt \tilde{B}(b,t)}\right \} \leq e^{4p^2 t^3 \sum_{n \geq 1}n^{2b}a_n^2 + 4pt^{3/2} \sqrt{2}(\sqrt{\log 2} + 2\sqrt{\pi})}.\]

\noindent These bounds may be applied to give the desired result.\qed

\paragraph{Proof of Theorem~\ref{bound2}} Note that 
\begin{eqnarray*}
\sup_{0 \leq s \leq t} \left |\frac{\partial}{\partial x_j}
m^{(\epsilon)}(s;x_1,\ldots, x_p) \right | &=& \sup_{0 \leq s \leq t}
\left |E_{\bf Q} \left \{ \left (\prod_{k \neq j}
u^{(\epsilon)}(s,x_k) \right )u_x(s,x_j)\right \}\right |
\\ & \leq & \sup_{0 \leq s \leq t} E_{\bf Q} \{ 
\tilde{C}(1,s)^{p-1} |u_x(s,x_j)|\}
\end{eqnarray*}

\noindent and the result now follows by applying lemma~\ref{derbd} with
$b_1 = \ldots = b_{p-1} = 1$. \qed 

\begin{Lmm}\label{epsder} Suppose that $\sum_{n \geq 1} n^4 |a_n| <
+\infty$ so that $\sum_{n \geq 1}n^8 a_n^2 < +\infty$.  For any $T \geq 0$,
there exists a constant
$C(p,T)$  such that 

\[ \sup_{0 \leq t \leq T}\sup_{0 \leq \epsilon \leq 1}\frac{\partial}{\partial
\epsilon} m_{2p}^{(\epsilon)}(t,{\bf 0}) \leq C(p,T).
\]

\end{Lmm}\vspace{5mm}

\noindent Note that this lemma gives no information on a lower bound.
\vspace{5mm}

\noindent {\bf Proof} For $\epsilon > 0$, set $\tilde{u} = \frac{\partial u^{(\epsilon)}}{\partial
\epsilon}$. Then $\tilde{u}$ is differentiable in $t$ and satisfies

\begin{equation}\label{equeps}\left\{\begin{array}{l} \tilde{u}_t = \frac{\epsilon}{2} \tilde{u}_{xx} -
u\tilde{u}_x - u_x
\tilde{u} + \frac{1}{2}u_{xx}\\ \tilde{u}(0,x) \equiv 0. \end{array}\right. \end{equation}
\noindent Recall (suppressing the notation $\epsilon$ from the process $X$) that 

\[ X_{s,t}(x) = x + \sqrt{\epsilon}(w_t - w_s) - \int_s^t
u^{(\epsilon)}\left (r,X_{r,t}(x) \right )dr,\]

\noindent so that 

\[ \frac{\partial X_{s,t}(x)}{\partial x} = 1 - \int_s^t
u^{(\epsilon)}_x \left (r,X_{r,t}(x) \right ) \frac{\partial X_{r,t}(x)}{\partial x} dr\]

\noindent yielding

\begin{equation}\label{eqnicfor} \frac{\partial X_{s,t}(x)}{\partial x} = e^{-\int_s^t u_x(r,X_{r,t}(x)) dr}.\end{equation}

\noindent It follows from equation (\ref{equeps}), using equation (\ref{eqnicfor}) to go from the first line to the second, that 

\begin{eqnarray*} \tilde{u}(t,x) &=&   \frac{1}{2}\int_0^t E_{\bf
P}\left [u_{xx}(s,X_{s,t}(x))e^{-\int_s^t u_x(r,X_{r,t}(x))dr}\right ] ds \\ 
&=& 
\frac{1}{2}\int_0^t E_{\bf
P}\left [u_{xx}(s,X_{s,t}(x))\frac{\partial X_{s,t}(x)}{\partial x}\right ] ds \\
&=&   \frac{1}{2}\frac{\partial}{\partial x}\int_0^t E_{\bf
P}\left [ u_x(s,X_{s,t}(x))\right ] ds,
\end{eqnarray*}

\noindent yielding (for $\epsilon > 0$)

\begin{eqnarray*}  \frac{\partial}{\partial \epsilon}
m^{(\epsilon)}(t;x_1,\ldots, x_p) &=& \sum_{j=1}^p E_{\bf Q}\left\{ \left(
\prod_{k
\neq j} u(t,x_k)\right)\tilde{u}(t,x_j)\right\} \\
&=&  \frac{1}{2}\sum_{j=1}^p \frac{\partial}{\partial
x_j}\int_0^t E_{\bf Q}
\left\{ \left (\prod_{k \neq j} u(t,x_k) \right )E_{\bf P} \left [ u_x (s,X_{s,t}(x_j)) \right ]\right\} ds
\end{eqnarray*}

\noindent  so that 

\begin{eqnarray} \nonumber \frac{\partial}{\partial \epsilon}
m_p^{(\epsilon)}(t,{\bf 0}) &=& \frac{p}{2}\int_0^t
\frac{1}{2\pi}\int_0^{2\pi} E_{\bf Q} \left \{u^{p-1}(t,x)
\frac{\partial}{\partial x} E_{\bf P} \left [u_x(s,X_{s,t}(x)) \right ] \right \} dx ds \\
&=&  \label{lowernow} -\frac{p(p-1)}{2}\int_0^t E_{\bf Q} \left \{ u^{p-2}(t,x)
u_x(t,x) E_{\bf P}\left [u_x(r,X_{r,t}(x)) \right ] \right \} dr.
\end{eqnarray}

\noindent Now, set $v = u_x$ and note that 

\[ \left\{\begin{array}{l} \partial_t v = (\frac{\epsilon}{2} v_{xx} - v^2 -
uv_x)dt +
\partial_t
\zeta_{xx} \\ v(0,x) \equiv 0.\end{array}\right. \]

\noindent It follows that 

\begin{eqnarray*} v(t,x) &=& - \sum_{n \geq 1} n^2 a_n \left (\beta^{1n}(t)\cos(nx) + \beta^{2n}(t)\sin(nx)\right)\\ && - \sum_{n \geq 1} n^2 a_n \left(\int_0^t \beta^{1n}(s){\cal Q}_{s,t}({\cal L}_s\cos(n.))(x) ds + \int_0^t \beta^{2n}(s) {\cal Q}_{s,t} ({\cal L}_s \sin(n.))(x) ds \right) \\&&   - \int_0^t
E_{\bf P}[v^2(s,X_{s,t})] ds.\end{eqnarray*}

\noindent Recall the definition of $\theta$ given in equation (\ref{thetlabel}). For $b=2$, 

\begin{eqnarray*} \theta(2;t,x) &=&  -\sum_{n \geq 1}n^2 a_n \left (\beta^{1n}(t)\cos (nx) + \beta^{2n}(t) \sin(nx)\right)\\&& - \sum_{n \geq 1}n^2 a_n \left( \int_0^t \beta^{1n}(s) {\cal Q}_{s,t}({\cal L}_s \cos (n.))(x) ds + \int_0^t \beta^{2n}(s) {\cal Q}_{s,t}({\cal L}_s \sin(n.))(x)ds\right)  
\end{eqnarray*}

\noindent  so that 

\[u_x(t,x) = \theta(2;t,x) - \int_0^t E_{\bf P} [v^2(s,X_{s,t}(x))]ds.\]

\noindent Note that 

\[ E_{\bf P} [u_x (r,X_{r,t}(x))] = E_{\bf P}[\theta (2;r,X_{r,t}(x))] - \int_0^r
E_{\bf P}[v^2(s,X_{s,t}(x))] ds\]

\noindent It follows that 

\begin{eqnarray*} \lefteqn{ u_x(t,x) E_{\bf P}[u_x(r,X_{r,t}(x))] = 
\theta(2;t,x)E[\theta(2;r,X_{r,t})] }\\&& + \int_0^t E_{\bf P}
[v^2(s,X_{s,t}(x))]ds
\int_0^r E_{\bf P} [v^2(s,X_{s,t}(x))] ds \\ &&   
- \theta(2;t,x)\int_0^r E_{\bf P} [v^2(s,X_{s,t}(x))]ds -
E_{\bf P} [\theta(2;r,X_{r,t}(x))]\int_0^t E_{\bf P} [v^2(s,X_{s,t}(x))]ds.
\end{eqnarray*}

\noindent Putting this into equation (\ref{lowernow}) gives 

\begin{eqnarray*}  \frac{\partial}{\partial \epsilon}
m_{2p}^{(\epsilon)}(t, {\bf 0})
&\leq & - p(2p-1) \left (\int_0^t E_{\bf Q}\left \{ u^{2(p-1)}(t,x)\theta(2;t,x)E_{\bf
P}[\theta(2; r,X_{r,t}(x))]\right \}dr\right.  \\&&
-\int_0^t E_{\bf
Q} \left \{u^{2(p-1)}(t,x)\theta(2;t,x)\int_0^rE_{\bf P}[v^2(s,X_{s,t}(x)]ds \right \}dr\\&&  -
\left. \int_0^t E_{\bf Q} \left \{u^{2(p-1)}(t,x)E_{\bf
P}[\theta(2;r,X_{r,t}(x))]\int_0^t E_{\bf P} [v^2(s,X_{s,t}(x))]ds \right \}dr \right).
\end{eqnarray*} 

\noindent Note that 

\[ 0 \leq \int_0^rE_{\bf P}[v^2(s,X_{s,t}(x))]ds \leq \int_0^tE_{\bf
P}[v^2(s,X_{s,t}(x))]ds = \theta(2;t,x) - u_x(t,x),\]

\noindent from which it follows that 

\begin{eqnarray*}  \frac{\partial}{\partial \epsilon}
m_{2p}^{(\epsilon)}(t, {\bf 0})
&\leq & - p(2p-1) \left (\int_0^t E_{\bf Q}\left \{ u^{2(p-1)}(t,x)\theta(2;t,x)E_{\bf
P}[\theta(2; r,X_{r,t}(x))]\right \}dr\right.  \\&&
- t E_{\bf
Q} \left \{u^{2(p-1)}(t,x)\theta(2;t,x)(\theta(2;t,x) - u_x(t,x)) \right \}dr\\&&  -
\left. \int_0^t E_{\bf Q} \left \{u^{2(p-1)}(t,x)E_{\bf
P}[\theta(2;r,X_{r,t}(x))](\theta(2;t,x) - u_x(t,x)) \right \}dr \right).
\end{eqnarray*} 

\noindent Recall the definition of $\tilde{C}^{(\epsilon)}$ given in
equation  (\ref{ceetild}); namely, 

\[\tilde{C}^{(\epsilon)}(b,t) = \sup_{0
\leq s \leq t} \sup_{0 \leq x \leq 2\pi}\left |\theta^{(\epsilon)}(b,s,x) \right |.\]

\noindent In
particular,  

\[ \tilde{C}^{(\epsilon)} (1,t) =
 \sup_{0
\leq s
\leq t}\sup_{0
\leq x
\leq 2\pi}\left |u^{(\epsilon)}(s,x) \right |\]

\noindent  and 

\[ \tilde{C}^{(\epsilon)}(2,t)=   \sup_{0
\leq s
\leq t}\sup_{0
\leq x
\leq 2\pi}\left |\theta^{(\epsilon)}(2;s,x)\right |.\]

\noindent It follows that 

\begin{eqnarray*}\lefteqn{  \sup_{0 \leq s \leq t}
\frac{\partial}{\partial
\epsilon}m_{2p}^{(\epsilon)}(s,{\bf 0}) 
\leq p(2p-1)}\\&& \times
\left (3t E_{\bf Q} \left \{
\tilde{C}^{(\epsilon)}(1,t)^{2(p-1)}\tilde{C}^{(\epsilon)}(2,t)^{2}\right \}  + 2t
\sup_{0
\leq s
\leq t} E_{\bf
Q}\left \{\tilde{C}^{(\epsilon)}(1,t)^{2(p-1)}\tilde{C}^{(\epsilon)}(2,t)
 \left |u_x(s,x) \right |\right \}\right).
\end{eqnarray*}

\noindent The first term on the right hand side is bounded, independently of
$\epsilon$, by an application of lemma~\ref{bound1lem}. The second term is
bounded, independently of $\epsilon$ by an application of lemma~\ref{derbd}.  \qed 

\section{The Moment Equations} \label{momsect}
 Having constructed a priori bounds for the moments of solutions to
equation (\ref{burger}) and a priori bounds on the Lipschitz constant, which are 
independent of $\epsilon$, the system of equations for the moment fields
is now considered. \vspace{5mm}

\noindent Recall that   $u^{(\epsilon)}(t,.) \in C^{2,1}([0,2\pi])$. It follows from equation (\ref{burger}) that 
$u^{(\epsilon)}(.,x)$ is a semimartingale for each $x \in [0,2\pi]$ and therefore It\^o's formula may be applied to  $f(u^{(\epsilon)}(.,x_1),\ldots, u^{(\epsilon)}(.,x_p)) = u^{(\epsilon)}(.,x_1)\ldots
u^{(\epsilon)}(.,x_p)$.  Set $\Gamma(z) = \sum_{n \geq 1} a_n^2 \cos(nz)$. It\^o's
formula yields 

\begin{eqnarray} \nonumber \lefteqn{\prod_{j=1}^p u^{(\epsilon)}(t,x_j) =
\frac{\epsilon}{2}
\int_0^t
\sum_{k=1}^p u_{xx}^{(\epsilon)}(s,x_k)\prod_{j \neq k}
u^{(\epsilon)} (s,x_j) ds -
 \int_0^t
\sum_{k=1}^p  u^{(\epsilon)}_x(s,x_k)  \prod_{j =1}^p u^{(\epsilon)}(s,x_j)
ds}
\\&&
  + \label{itook} \sum_{k=1}^p \int_0^t \left (\prod_{j \neq k}
u^{(\epsilon)}(s,x_j) \right )
\partial_s
\zeta (s,x_k) + \sum_{j < k} \int_0^t \left (\prod_{l \neq j,k}
u^{(\epsilon)}(s,x_l) \right ) (-\Gamma^{\prime \prime}(x_j - x_k)) ds.
\end{eqnarray}

\noindent Recall the definition of $m_p^{(\epsilon)}$ given in equation (\ref{momfield}). To obtain an equation for $m_p$ from equation (\ref{itook}), it
is necessary to show that Fubini's theorem may be used on each term of
the right hand side of equation (\ref{itook}) and that the martingale term
is indeed a martingale.  Fubini's theorem may be applied to the second and
fourth terms, using theorems~\ref{bound2} and~\ref{bound1} 
respectively. For the (local) martingale term, note that by the
Burkholder Davis Gundry inequality (theorem 4.1 and corollary 4.2 on
pages 160 and 161 of Revuz and Yor \cite{RY}) there exist constants
$K(q) < +\infty$ such that 

\begin{eqnarray*}\lefteqn{ E_{\bf Q}\left \{ \sup_{0 \leq s \leq t}\left
|\left (
\int_0^s
\prod_{j
\neq k} u^{(\epsilon)} (r,x_j) \partial_r \zeta (r,x_j)\right)\right |^q
\right \}}\\ &&\leq  K(q)  E_{\bf Q}\left \{ \left ( \int_0^t (\prod_{j
\neq k} u^{(\epsilon)} (r,x_j))^2 (-\Gamma^{\prime \prime}(0))
dr\right)^{q/2}
\right \} \\
&& \leq K(q) (-\Gamma^{\prime
\prime}(0))^{q/2}t^{(q/2)-1} \int_0^t E_{\bf Q}\left \{
|u(r,x)|^{ q} \right \} dr
\end{eqnarray*}

\noindent which is bounded above by theorem \ref{bound1}. This
gives that the (local) martingale  is uniformly integrable and is therefore
a martingale. \vspace{5mm}

\noindent The only `problem' term is the first one on the right hand side.
Recall equation (\ref{CH}), where $U^{(\epsilon)}$ is given by equation
(\ref{heat}). This may be rewritten as

\[ u^{(\epsilon)} = -\epsilon\frac{U^{(\epsilon)}_x}{U^{(\epsilon)}}.\]

\noindent Then

\begin{equation}\label{secondder} u^{(\epsilon)}_{xx} = -\epsilon \left(
\frac{U_{xxx}^{(\epsilon)}}{U^{(\epsilon)}} -
\frac{3U_{xx}^{(\epsilon)}U_x^{(\epsilon)}}{U^{(\epsilon)2}} +
2\frac{U_x^{(\epsilon)3}}{U^{(\epsilon)3}}\right).
\end{equation}

\noindent The solution of equation (\ref{heat}), with initial condition
$U^{(\epsilon)}(0,x) \equiv 1$, may be expressed using the Feynman Kacs
representation in equation (\ref{FKrep}). Using  

\[ f(x) = \sum_{n \geq 1}a_n \left(\int_0^t \cos(n(x + w^{(\epsilon)}_{t-s}))d\beta^{1n}_s + \int_0^t \sin (n(x + w^{(\epsilon)}_{t-s}))d\beta^{2n}_s \right)\]
 to simplify notation, where under ${\bf P}$, $w^{(\epsilon)}$ is a Brownian motion, with $w^{(\epsilon)}_0 = 0$ and diffusion coefficient $\epsilon$, it follows that 

\begin{equation}\label{youone} U_x^{(\epsilon)} =- \frac{1}{\epsilon}E_{\bf
P} [ e^{-\frac{1}{\epsilon} f(x)}f_x(x)],\end{equation}

\begin{equation}\label{youtwo} U_{xx}^{(\epsilon)} =
\frac{1}{\epsilon^2}E_{\bf P} [e^{-\frac{1}{\epsilon} f(x)} (f_x(x))^2] -
\frac{1}{\epsilon} E_{\bf P} [e^{-\frac{1}{\epsilon} f(x)} f_{xx}(x)]
\end{equation}

\noindent and

\begin{equation}\label{youthree} U_{xxx}^{(\epsilon)} =
-\frac{1}{\epsilon^3} E_{\bf P}[e^{-\frac{1}{\epsilon} f(x)}(f_x(x))^3] +
\frac{3}{\epsilon^2}E_{\bf P}[e^{-\frac{1}{\epsilon} f(x)} f_x(x) 
f_{xx}(x)] - \frac{1}{\epsilon} E_{\bf P}[e^{-\frac{1}{\epsilon}
f(x)}f_{xxx}(x)].\end{equation}

\noindent To show that Fubini's theorem may be applied to the first term
on the right hand side of equation (\ref{itook}), note that 

\begin{eqnarray}\label{scnddrst} \lefteqn{E_{\bf Q} \left\{ 
 \left |u_{xx}^{(\epsilon)}(s,x_k) \right |
\prod_{j
\neq k} \left |u^{(\epsilon)}(s,x_j) \right |\right \}}\\&& \nonumber
\leq E_{\bf Q} \left\{ \left |u_{xx}^{(\epsilon)}(s,x) \right |^2 \right \} ^{1/2}E_{\bf
Q} \left \{ \left |u^{(\epsilon)}(s,x) \right |^{2(p-1)}\right \}^{1/2}
\leq C(s) E_{\bf Q} \left\{ \left |u_{xx}^{(\epsilon)}(s,x) \right |^2\right \}^{1/2},
\end{eqnarray}

\noindent where the constant $C(s) < +\infty$, increasing in $s$, is
obtained by an application of theorem~\ref{bound1} and is independent of
$\epsilon$.

\noindent Since the arguments for all the terms obtained by applying
equations (\ref{youone}), (\ref{youtwo}) and (\ref{youthree}) to equation
(\ref{secondder}) in estimating the right hand side of equation
(\ref{scnddrst}) are similar, only one will be sketched. Note that (for
example)

\[ E_{\bf Q}\left \{ \left |\frac{U_{xxx}^{(\epsilon)}}{U^{(\epsilon)}}\right
|^2\right \}
\leq E_{\bf Q} \left\{ U_{xxx}^{(\epsilon)4}  \right\} ^{1/2} E_{\bf Q}\left
\{ \frac{1}{U^{(\epsilon)4}}\right \}^{1/2}.\]

\noindent Note that, since $f(x)$ is Gaussian, for any $q \in {\bf R}$,  

\[ E_{\bf Q}\left \{  \left(e^{-\frac{1}{\epsilon}f(x)}\right)^q\right \} =
e^{\frac{q^2}{2\epsilon^2}\Gamma(0) t}\]

\noindent It follows (using Jensen's inequality on the function
$\frac{1}{x}$ for $x \in (0,+\infty)$, which is convex in that region) that 

\[ E_{\bf Q} \left \{ \frac{1}{U^{(\epsilon)4}}\right \} = E_{\bf Q}\left\{\frac{1}{E_{\bf P} \left [e^{-\frac{1}{\epsilon} f(x)}\right ]^4} \right\} \leq E_{\bf Q} \left \{ 
E_{\bf P}
\left [ e^{\frac{4}{\epsilon} f(x) }\right ] \right \} =
e^{\frac{16}{\epsilon^2}\Gamma(0) t}.\]

\noindent Let $f^{(n)} = \frac{\partial^n}{\partial x^n} f(x)$. Note that
${\bf P}$
almost surely $f^{(n)}(x)$,
is Gaussian with respect to ${\bf Q}$, with 
$E_{\bf Q} \left \{ f^{(n)2q-1}\right\}  = 0$ and
$E_{\bf Q} \left\{ f^{(n)2q} \right\}  = (-1)^{n}(\Gamma^{(2n)}(0))^{q}t^q\prod_{j=1}^q(2j -
1)$, for all integer $q \geq 1$, where $\Gamma^{(2n)}$ denotes the $2n$th derivative of $\Gamma$. It
is now easy to use Hölder's inequality to compute an upper bound for
$E_{\bf Q}\left\{ U_{xxx}^{(\epsilon)2k}\right\} $, $E_{\bf Q}\left\{ U_{xx}^{(\epsilon)2k}\right\}$ and
$E_{\bf Q} \left\{ U_x^{(\epsilon)2k}\right\}$ for all $k \geq 1$, since these will involve $\Gamma^{(2n)}(0)$ for $n \leq 3$, which is bounded by hypothesis~\ref{hyp}. These bounds depend on
$\epsilon$ and are increasing as $t \rightarrow +\infty$ and as $\epsilon
\rightarrow 0$. It now follows that there exists a non negative function
$C(\epsilon,t)$, increasing in $t$, such that for any
$t < +\infty$ and any $\epsilon > 0$, $C(\epsilon, t) < +\infty$ and such
that 

\[ \sup_{0 \leq s \leq t} \sup_{(x_1,\ldots, x_p) \in {\bf R}^p} E_{\bf
Q}\left\{ \left |u_{xx}^{(\epsilon)}(s,x_k)\right |\prod_{j \neq k} \left |u^{(\epsilon)}(s,x_j) \right | \right\} 
\leq C(\epsilon, t).\]

\noindent It follows that, for fixed $\epsilon > 0$ and $t < +\infty$, 
Fubini's theorem may be applied to equation (\ref{itook}). It has already
been seen that the martingale term is a martingale, starting from $0$ at $t = 0$ and therefore has
expected value
$0$. It follows that

 \begin{eqnarray}\label{moments}\lefteqn{\frac{\partial}{\partial
t}m_p^{(\epsilon)}(t;x_1,
\ldots, x_p) =
\frac{\epsilon}{2}\Delta_{\bf x} m_p^{(\epsilon)}(t;x_1,\ldots, x_p)}\\&&
\nonumber - \frac{1}{2}
\sum_{j=1}^p \frac{\partial}{\partial x_j} m_{p+1}^{(\epsilon)}(t;x_1,\ldots,
x_p, x_j) + \sum_{k < l} (-\Gamma^{\prime \prime}(x_k -
x_l))m_{p-2}^{(\epsilon)}(t;\hat{x_k}, \hat{x_l}),
\end{eqnarray}

\[ m_p^{(\epsilon)}(0;x_1,\ldots, x_p) \equiv 0\]

\noindent where $\frac{\partial}{\partial x_j}$ means differentiation with
respect to both appearances of $x_j$ and $m_{p-2}(\hat{x_k}, \hat{x_l})$ means
that the variables $ x_k$ and $x_l$ are excluded; $m_{p-2}$ is a function of the
other $p-2$ space variables.   \vspace{5mm}

\noindent $\tilde{u}(t,x) = - u(t, -x)$  and $\tilde{\zeta}(t,x)  = \zeta(t,-x)$. Then, it is straightforward to compute that

\[ \left\{\begin{array}{l} \partial_t \tilde{u} =
(\frac{\epsilon}{2}\tilde{u}_{xx} -
\frac{1}{2}(\tilde{u}^2)_x)dt + \partial_t \tilde{\zeta}_x\\ \tilde{u}(0,x)
\equiv 0.\end{array}\right.\]

\noindent From this,  and
noting that $\zeta$ and $\tilde{\zeta}$ are identically distributed, it
follows that 

\begin{eqnarray*} m_p(t;x_1,\ldots, x_p) &=& E_{\bf Q}\{u(t,x_1)\ldots u(t,x_p)\}\\ &=& E_{\bf Q}\{\tilde{u}(t,x_1)\ldots \tilde{u}(t,x_p)\} = (-1)^p m_p(t; -x_1, \ldots, -x_p).\end{eqnarray*}

\noindent It follows that $m_{2p}^{(\epsilon)}(t,.)$ is an even function and that
$m_{2p+1}^{(\epsilon)}(t,.)$ is an odd function for each integer $p \geq 0$. That is, for each integer $p \geq 1$, 

\begin{equation}\label{momeven}
 m_{2p}^{(\epsilon)}(t;x_1,\ldots, x_{2p}) = m^{(\epsilon)}_{2p}(t;-x_1,\ldots,-x_{2p})\qquad \forall \epsilon \geq 0, \quad {\bf x} \in {\bf R}^{2p}, \quad t \in {\bf R}_+
\end{equation}

\noindent and

\begin{equation}\label{momodd}  m_{2p+1}^{(\epsilon)}(t; x_1, \ldots, x_{2p+1}) = -m_{2p+1}(t; -x_1, \ldots, -x_{2p+1})\qquad \forall \epsilon \geq 0, \quad {\bf x} \in {\bf R}^{2p+1}, \quad t \in {\bf R}_+.\end{equation} 

\noindent In particular,  
$m^{(\epsilon)}_{2p+1}(t,{\bf 0})
\equiv 0$ for all $t \geq 0$ and all $p \geq 0$. This, together with the upper
and lower bounds on $m_p^{(\epsilon)}(t,{\bf 0})$ uniform in $\epsilon$ and
together with lemma~\ref{epsder} gives that for each integer $p \geq 0$
and all $T \geq 0$, 

\[ \limsup_{0 \leq \epsilon_1 \leq
\epsilon_2 
\rightarrow 0} \sup_{0 \leq t \leq T}|m^{(\epsilon_1)}_p(t,{\bf 0}) -
m^{(\epsilon_2)}_p(t,{\bf 0})| = 0.\]   

\noindent Set $M_p^{(\epsilon)}(t) := m^{(\epsilon)}_p(t;0,\ldots,
0)$. 

\begin{Lmm} \label{mexist} For all non negative integer $p$ and for all $t \in {\bf R}_+$, the limit $M_p(t):=
\lim_{\epsilon
\rightarrow 0} M_p^{(\epsilon)}(t)$ is well defined.
\end{Lmm}

\noindent {\bf Proof of Lemma~\ref{mexist}} Firstly, by theorem~\ref{bound1}, $\sup_{0 < \epsilon < 1} \sup_{0 \leq t \leq T} |M_p^{(\epsilon)}(t) | < K(p,T)$, where $K(p,T)$ is defined on the right hand side of inequality (\ref{thbd1}). Consider $p$ odd; that is $p = 2q + 1$ for non negative integer $q$. Then $M_{2q+1}^{(\epsilon)}(t) \equiv 0$, so that $M_{2q+1}(t) := \lim_{\epsilon \rightarrow 0} M_{2q+1}^{(\epsilon)}(t) = 0$. Secondly, consider $p$ even; that is $p = 2q$ for non negative integer $q$. Then $M_{2q}^{(\epsilon)}(t) \geq 0$ for all $\epsilon > 0$ and all $t \in {\bf R}$. Let
\[ \overline{M}_{2q}^{(\epsilon)} = \sup_{0 < \delta \leq \epsilon} M_{2q}^{(\delta)}(t), \qquad \underline{M}_{2q}^{(\epsilon)} = \inf_{0 < \delta \leq \epsilon} M_{2q}^{(\delta)}(t).\]

\noindent Then, from lemma~\ref{epsder}, which states that  for each $t > 0$, there is a constant $C(q,T) < +\infty$ such that $\sup_{\epsilon > 0} \frac{\partial}{\partial \epsilon} M_{2q}^{(\epsilon)}(t) < C(q,T)$ for all $t \in (0,T)$, it follows that 

\[ \underline{M}_{2q}^{(\epsilon)}  \leq \overline{M}_{2q}^{(\epsilon)} \leq \underline{M}_{2q}^{(\epsilon)} + \epsilon C(q,T),\]
\noindent from which
\[ \lim_{\epsilon \rightarrow 0} \underline{M}_{2q}^{(\epsilon)}(t) = \lim_{\epsilon \rightarrow 0} \overline{M}_{2q}^{(\epsilon)}(t) = \lim_{\epsilon \rightarrow 0} M_{2q}^{(\epsilon)}(t) .\]\qed  

\noindent Set 

\begin{equation}\label{mudefn} \mu_p^{(\epsilon)}(t;x_1, \ldots, x_p) :=
m^{(\epsilon)}_p(t;\epsilon x_1,
\ldots,
\epsilon x_p).\end{equation}

\noindent Let $K(p,T)$ denote the uniform Lipschitz constant for
$(m^{(\epsilon)}(t,.))_{0 \leq \epsilon \leq 1, 0 \leq t \leq T}$ found in
theorem~\ref{bound2}. That is, for the remainder of the article, $K(p,T)$ will denote a finite positive constant such that 

\begin{eqnarray}\nonumber \lefteqn{ \sup_{0 \leq t \leq T} \sup_{j \in \{1, \ldots, p\}} \sup_{(x_1,\ldots, x_p) \in {\bf R}^p}\left |\frac{\partial}{\partial x_j} m_p^{(\epsilon)}(t;x_1,\ldots, x_p) \right |}\\&& \label{Kdef} \leq \sup_{0 \leq t \leq T} \sup_{(x_1, \ldots, x_p) \in {\bf R}^p} E_{\bf Q} \left \{ 
\left |u^{(\epsilon)}_x(t,x_1) u^{(\epsilon)}(t,x_2)\ldots u^{(\epsilon)}(t,x_p)\right |
\right \}
\leq K(p,T)\end{eqnarray}

Note that, for any fixed
$x_1,
\ldots, x_p$,  

\begin{eqnarray*} \lefteqn{\lim_{\epsilon \rightarrow 0}
|\mu^{(\epsilon)}_p(t;x_1,
\ldots, x_p) - M_p(t)|} \\ &&   \leq  \lim_{\epsilon \rightarrow
0} |\mu^{(\epsilon)}_p(t;x_1,\ldots, x_p) - M^{(\epsilon)}_p(t)| +
\lim_{\epsilon \rightarrow 0}|M^{(\epsilon)}_p(t) - M_p(t)| \\
&& \leq \lim_{\epsilon \rightarrow 0} \epsilon (\sum_{j=1}^p |x_j|) K(p,T) +
0 = 0.
\end{eqnarray*}

\noindent  Set

\begin{equation}\label{phidefn} \phi^{(\epsilon)}_p(t;x_1, \ldots, x_p) :=
\frac{\mu^{(\epsilon)}_p (t;x_1,
\ldots, x_p) - M_p^{(\epsilon)}(t)}{\epsilon}.
\end{equation}

\noindent Note that, by equations (\ref{momodd}) and (\ref{momeven}),   $\phi^{(\epsilon)}_{2p + 1}$ is an odd function and $\phi^{(\epsilon)}_{2p}$ is an even function for all integer $p \geq 1$, all $\epsilon \geq 0$ all $t \in {\bf R}_+$. That is,

\begin{equation}\label{phiodd}  \phi_{2p+1}^{(\epsilon)}(t; x_1, \ldots, x_{2p+1}) = - \phi_{2p+1}(t; -x_1, \ldots, -x_{2p+1})\qquad \forall \epsilon \geq 0, \quad {\bf x} \in {\bf R}^{2p+1}, \quad t \in {\bf R}_+ \end{equation}

\noindent and 

\begin{equation}\label{phieven}
 \phi_{2p}^{(\epsilon)}(t;x_1,\ldots, x_{2p}) = \phi^{(\epsilon)}_{2p}(t;-x_1,\ldots,-x_{2p})\qquad \forall \epsilon \geq 0, \quad {\bf x} \in {\bf R}^{2p}, \quad t \in {\bf R}_+.
\end{equation}
 
\begin{Lmm}\label{LIP} It holds that 

\begin{equation}\label{lips}\sup_{0 \leq t \leq T} \sup_{0 < \epsilon \leq 1}
\sup_{x_1,\ldots, x_p}\frac{|\phi^{(\epsilon)}_p (t;x_1,
\ldots, x_p)|}{(\sum_{j=1}^p |x_j|)} \leq
 K(p,T)\end{equation}

\noindent and, for all $j \in \{1, \ldots, p\}$, 

\begin{equation}\label{lips2}
\sup_{0 \leq t \leq T}\sup_{0 < \epsilon \leq 1}\sup_{x_1,\ldots,
x_p}\limsup_{h \rightarrow 0}
\frac{|\phi^{(\epsilon)}_p (t;x_1,\ldots,x_j+h,\ldots, x_p) -
\phi^{(\epsilon)}_p (t;x_1,\ldots, x_j,\ldots, x_p)|}{|h|} \leq K(p,T)
\end{equation}

\noindent where the existence of a constant $K(p,T)$ independent of $\epsilon$, is guaranteed by theorem~\ref{bound2}.
\end{Lmm} 

\paragraph{Proof} This is an immediate consequence of Taylor's expansion
theorem, together with theorem~\ref{bound2}. In the second part, for
example,

\begin{eqnarray*}\lefteqn{\sup_{0 \leq t \leq T} \left
|\frac{\phi^{(\epsilon)}_p (t;x_1,\ldots, x_j+h,
\ldots, x_p) -
\phi^{(\epsilon)}_p (t;x_1,\ldots, x_p)}{h}\right |}\\ &&=  \sup_{0 \leq t
\leq T} \left |\frac{m^{(\epsilon)}_p (t;\epsilon x_1,\ldots, \epsilon x_j+
\epsilon h, \ldots,
\epsilon x_p) -
m^{(\epsilon)}_p (t; \epsilon x_1,\ldots, \epsilon x_p)}{\epsilon h}\right |\\&&
\leq  K(p,T).
\end{eqnarray*}
\noindent Equation (\ref{lips}) also follows directly from the definition, using 
\[|\phi^{(\epsilon)}_p(t;x_1, \ldots, x_p)| = |\frac{1}{\epsilon}(m^{(\epsilon)}(t; \epsilon x_1, \ldots, \epsilon x_p) - m^{(\epsilon)}(t;0,\ldots,0))| \leq K(p,T) \sum_{j=1}^p |x_j| \]
by Taylor's expansion theorem. 
 \qed.
\vspace{5mm}

\noindent  Let $\Phi^{(\epsilon)}_p : {\bf R}_+ \times {\bf R}^p \rightarrow {\bf R}$ be used to denote the function   

\begin{equation}\label{Phidef}  \Phi^{(\epsilon)}_p(t;.) =  \int_0^t  \phi_p^{(\epsilon)}(\alpha,.) d\alpha. 
\end{equation} 

\begin{Lmm} \label{bddapriori}$M_p$ is Lipschitz. That is, for each $T <
+\infty$, there exists a constant $C(p,T)$ such that 
\[ \sup_{0 \leq t \leq T}\left(\limsup_{h \rightarrow 0}
\frac{|M_p(t+h) - M_p(t)|}{h}\right)  \leq C(p,T).\]
\end{Lmm}
\noindent {\bf Proof of lemma~\ref{bddapriori}} Note that 

\begin{eqnarray*}\lefteqn{ \frac{\partial}{\partial t} \mu^{(\epsilon)}_p(t;x_1,
\ldots, x_p) =
\frac{1}{2}\Delta_{\bf x} \phi^{(\epsilon)}_p(t;x_1,\ldots, x_p)}\\&& -
\frac{1}{2}\sum_{j=1}^p \frac{\partial}{\partial
x_j}\phi^{(\epsilon)}_{p+1}(t;x_1,\ldots, x_p,x_j) + \sum_{j <
k}(-\Gamma^{\prime \prime}(\epsilon (x_j - x_k))\mu^{(\epsilon)}(t;\hat{x_j},
\hat{x_k}).\end{eqnarray*}

\noindent This may be rearranged as

\begin{eqnarray}\nonumber\lefteqn{ \frac{1}{2}\Delta
\phi_p^{(\epsilon)}(t;x_1,\ldots, x_p) =   \frac{1}{2}\sum_{j=1}^p
\frac{\partial}{\partial x_{j}}\phi^{(\epsilon)}_{p+1}(t;x_1,\ldots,
x_p,x_j)}\\&& \label{epsilprelim} +
  \left\{\frac{\partial}{\partial t} \mu^{(\epsilon)}_p(t;x_1,
\ldots, x_p)-\sum_{j <
k}(-\Gamma^{\prime \prime}(\epsilon(x_j -
x_k))\mu_{p-2}^{(\epsilon)}(t;\hat{x_j},
\hat{x_k})   \right\}.
\end{eqnarray}

\noindent With a change of notation from earlier, set 

\[ P_s f({\bf x}) = \int_{{\bf R}^p} \frac{1}{(2\pi
s)^{p/2}}\exp \left \{-\frac{|{\bf x}-{\bf y}|^2}{2s} \right \} f({\bf y}) d{\bf y} \]

\noindent and set
$p(r;{\bf z}) = \frac{1}{(2\pi r)^{p/2}}\exp\left \{-\frac{|{\bf z}|^2}{2r} \right \}$. By
integrating all terms of equation (\ref{epsilprelim}) against the test
function
$\frac{1}{s}\int_0^s p(r;{\bf x}-{\bf y}) dr$, it follows
 that, for all
$s
> 0$,

\begin{eqnarray*}\lefteqn{ \frac{1}{s}\int_{{\bf R}^p}\int_0^s p(r;{\bf x} -
{\bf y})
\frac{1}{2}\Delta
\phi_p^{(\epsilon)}(t;{\bf y}) dr d{\bf y} -   \frac{1}{2s}\sum_{j=1}^p
\int_{{\bf R}^p} \int_0^s \frac{\partial}{\partial
x_{j}}p(r;{\bf x} - {\bf y})\phi_{p+1}^{(\epsilon)}(t;y_1,\ldots, y_p,y_j) dr
d{\bf y}}\\&&
  \hspace{10mm}= \frac{\partial}{\partial t}
\left(\frac{1}{s}\int_0^s P_r \mu^{(\epsilon)}_p (t;x_1,\ldots, x_p)
dr\right)   -\sum_{j < k} \frac{1}{s}\int_0^s P_r(-\Gamma^{\prime
\prime}(\epsilon(x_j-x_k))\mu_{p-2}^{(\epsilon)}(t;\hat{x}_j,\hat{x}_k))dr.  
\end{eqnarray*} 

\noindent Note that $\lim_{\epsilon \rightarrow 0} \frac{1}{s}\int_0^s P_r \mu^{(\epsilon)}_p (t;x_1,\ldots, x_p)
dr = M_p(t)$ and that 

\[ \lim_{\epsilon \rightarrow 0}\sum_{j < k} \frac{1}{s}\int_0^s
P_r(-\Gamma^{\prime
\prime}(\epsilon(x_j-x_k))\mu_{p-2}^{(\epsilon)}(t;\hat{x}_j,\hat{x}_k))dr
= \frac{p(p-1)}{2}(- \Gamma^{\prime \prime}(0))M_{p-2}(t).\]

\noindent Since
$\frac{1}{2}\Delta$ is the infinitesimal generator of $P_s$, it follows that for all $s > 0$,  

\begin{equation} \label{journal} \frac{1}{s}\int_{{\bf R}^p} \int_0^s p(r;{\bf x} -
{\bf y})
\frac{1}{2}\Delta
\phi_p^{(\epsilon)}(t;{\bf y}) dr d{\bf y} = \frac{P_s
\phi_p^{(\epsilon)}(t,{\bf x}) - \phi_p^{(\epsilon)}(t,{\bf x})}{s}.\end{equation}

\noindent The observation that equation (\ref{journal}) holds for all $s > 0$ has to be made. When the article did not have this, it received a referee report stating that equation (\ref{journal}) did not hold for all $s > 0$; the referee stated that equation (\ref{journal}) only held in the limit as $s \rightarrow 0$. The referee therefore assumed that $\lim_{s \rightarrow 0}$ was intended and that the author had made a `flagrant error'. This was from a respectable journal and the editor sent the author an electronic mail assuring him that the referee was `an expert in the field' (`the field' was left undefined). \vspace{5mm}

\noindent Since $p(s; {\bf x}) = \frac{1}{(2\pi s)^{p/2}}\exp\left\{-\frac{|{\bf x}|^2}{2s}\right\}$, it follows that for any $s > 0$ and any bounded continuous function $f : {\bf R}^p \rightarrow {\bf R}$,  

\[\frac{\partial}{\partial x_j}\int_{{\bf R}^{p}}\frac{1}{(2\pi s)^{p/2}}e^{-|{\bf x} - {\bf y}|^2/2s} f(  y_1,\ldots,
y_{p} ) d{\bf y}  = -  \int_{{\bf R}^p} \frac{x_j - y_j}{s}p(s,{\bf x} - {\bf
y})  f ( y_1,\ldots, y_p ) d{\bf y}. \]

\noindent It follows that for all $h > 0$ and $0 < t < T-h$,

\begin{eqnarray}\lefteqn{  M_p (t+h) - M_p (t) =
\frac{p(p-1)}{2}(-\Gamma^{\prime
\prime}(0))\int_t^{t+h} M_{p-2} (\tau)d\tau} \nonumber \\&& +
\lim_{\epsilon \rightarrow 0}\left(\frac{1}{s}\int_t^{t+h} 
(P_s\phi_p^{(\epsilon)} (\tau,{\bf x}) -  
\phi_p^{(\epsilon)} (\tau,{\bf x}))  
d\tau
\right. 
\nonumber
\\&&\left. \label{intout} -
\frac{1}{2s}\sum_{j=1}^p
 \int_{{\bf R}^p}\int_0^s
\frac{(x_j - y_j)}{r}p_r({\bf x}-{\bf
y})\int_t^{t+h} \phi_{p+1}^{(\epsilon)} ( \tau; y_1,\ldots, y_p,y_j)
  d\tau d{\bf y} \right).
\end{eqnarray}

\noindent The estimate from lemma~\ref{LIP} (namely, that
$\sup_{0 \leq t \leq T} |\phi_p^{(\epsilon)}(t,{\bf x})| \leq K(p,T)\sum_j |x_j|$) gives

\begin{eqnarray}\sup_{0 < \epsilon \leq 1} |P_s\phi_p^{(\epsilon)} (t,{\bf
x})| &\leq &
\nonumber K(p,T)\sum_j \int_{{\bf R}^p} |y_j| p(s;{\bf x} - {\bf y})d{\bf y}
\\ &\leq & K(p,T) \sum_j \int_{{\bf R}^p} (|x_j - y_j| + |x_j|)p(s;{\bf x} -
{\bf y})d{\bf y} \nonumber \\ 
&=&  K(p,T)\left(\sum_j
|x_j|\right) +
\sqrt{\frac{2s}{\pi}}p K(p,T). \label{lip2}\end{eqnarray}

\noindent Similarly,  

\begin{eqnarray*}\lefteqn{ \sup_{0 < \epsilon \leq 1}\left |\sum_{j=1}^p
\int_{{\bf R}^p}\int_0^s
\frac{(x_j - y_j)}{r}p_r({\bf x}-{\bf y})\phi_{p+1}^{(\epsilon)} (y_1,\ldots,
y_p,y_j) d{\bf y}\right |}\\&& \leq 
K(p+1,T)\sum_{j=1}^p \int_{{\bf R}^q}\int_0^s \frac{|x_j - y_j|}{r}p(r;{\bf
x} - {\bf y}) \{ 2|y_j| + \sum_{k \neq j}|y_k|\}d{\bf y}\\&&
\leq K(p+1,T)\sum_{j=1}^p \left \{2|x_j| + \sum_{j \neq k}|x_k|\right \}\int_{{\bf
R}^p}\int_0^s \frac{|x_j - y_j|}{r}p(r;{\bf x} - {\bf y}) d{\bf
y}\\&&\hspace{5mm} + K(p+1,T)\sum_{j=1}^p \int_{{\bf R}^p}\int_0^s
\frac{|y_j|}{r}p(r;{\bf y})\left \{2|y_j| + \sum_{k \neq j} |y_k|\right \} d{\bf y}
\\&&
= \frac{2\sqrt{2}(p+1)}{\sqrt{\pi}}K(p+1,T) \left (\sum_j
|x_j|\right )\sqrt{s} + 2p  K(p+1,T)\left (1 + \frac{p-1}{\pi} \right ) s.
\end{eqnarray*}

\noindent Now, using the upper bound, uniform in $\epsilon$ given by
equation (\ref{thbd1}) in theorem~\ref{bound1}, it follows that for $t \in
[0,T]$, there is a constant
$c(p,T) < +\infty$ such that $ \sup_{0 \leq t \leq T}  |M_{2(p-1)}(t)| \leq
c(p,T)$.

\noindent Taking $s \rightarrow +\infty$, it follows from equation
(\ref{intout}) that 

\begin{eqnarray}\label{conclu}\lefteqn{\sup_{0 \leq t \leq T}\limsup_{h
\rightarrow 0}|
\frac{M_{2p}(t+h) - M_{2p}(t)}{h}|}\\&& \nonumber 
\leq p(2p-1)(-\Gamma^{\prime \prime}(0))c(p,T) +  4p  K(2p+1,T) \left (1 + \frac{4(2p-1)}{\pi}\right ).
\end{eqnarray}

\noindent   Lemma~\ref{bddapriori} follows. 
\qed \vspace{5mm}

\begin{Lmm} \label{secondlast}   Let  $\Phi_{p+1}^{(\epsilon)}$ be defined according to equation (\ref{Phidef}). Then for all
$p \geq 2$

\begin{eqnarray*}\lefteqn{ \lim_{\epsilon \rightarrow 0} \left | M_p(t) -
\frac{p(p-1)}{2}\int_0^t M_{p-2}(\alpha)d\alpha\right . }\\&& - \left. 
\lim_{s
\rightarrow +\infty}\frac{1}{2s}   \sum_{j=1}^p
\int_0^s \int_{{\bf R}^p} \frac{x_j - y_j}{r}p(r,{\bf x} - {\bf
y})\Phi_{p+1}^{(\epsilon)}(t;y_1,\ldots, y_p,y_j) d{\bf y} dr\right | = 0.
\end{eqnarray*}
\end{Lmm}

\paragraph{Proof}  From equation (\ref{intout}), it follows that 
 
\begin{eqnarray}\lefteqn{ \nonumber \lim_{\epsilon \rightarrow 0}\left |M_p(t) -
\frac{p(p-1)}{2}\int_0^t M_{p-2}(\alpha)d\alpha \right .}\\&& - \label{vanish}
 \left .
\frac{\int_0^t (P_s\phi^{(\epsilon)}_p(\alpha,{\bf x}) -
\phi^{(\epsilon)}_p(\alpha, {\bf x}))d\alpha}{s}\right. 
\\&& - \left. \frac{1}{2s}  
\sum_{j=1}^p
\int_0^s \int_{{\bf R}^p} \frac{x_j - y_j}{r}p(r,{\bf x} - {\bf
y})\left\{\int_0^t \phi_{p+1}^{(\epsilon)}(\alpha;y_1,\ldots, y_p,y_j) d\alpha \right\} d{\bf y}
dr
\right | = 0.\nonumber
\end{eqnarray} 

\noindent This holds for all $s > 0$.  Recall equation (\ref{Phidef}).  Using equations (\ref{lips}) and (\ref{lip2}), it follows that for all $0 < t < T < +\infty$

\[\left | \frac{\int_0^t (P_s\phi^{(\epsilon)}_p(\alpha,{\bf x}) -
\phi^{(\epsilon)}_p(\alpha, {\bf x}))d\alpha}{s} \right |\leq \frac{2}{s}TK(p,T) \sum_{j=1}^p |x_j| + \frac{1}{\sqrt{s}}\sqrt{\frac{2}{\pi}} p T K(p,T), \]

\noindent so that 

\begin{eqnarray}\nonumber \lefteqn{ \lim_{\epsilon \rightarrow 0} \left | M_p(t) -
\frac{p(p-1)}{2}\int_0^t M_{p-2}(\alpha)d\alpha \right. }\\&& \nonumber  \hspace{30mm} \left. - \frac{1}{2s}  
 \sum_{j=1}^p
\int_0^s \int_{{\bf R}^p} \frac{x_j - y_j}{r}p(r,{\bf x} - {\bf
y})  \Phi_{p+1}^{(\epsilon)}(t ;y_1,\ldots, y_p,y_j) d{\bf y}
dr
 \right | \\&&
\hspace{70mm}\label{thelimiteq} \leq \frac{2}{s}TK(p,T) \sum_{j=1}^p |x_j| + \frac{1}{\sqrt{s}}\sqrt{\frac{2}{\pi}} p T K(p,T).
\end{eqnarray}

\noindent Now, letting 
$s \rightarrow +\infty$, the result follows.  \qed
\vspace{5mm}

\begin{Propn}\label{help2}  

\noindent Let $p \in {\bf Z}_+$ (the non negative integers) let 
$\Phi^{(\epsilon)}_{p+1}$, be the function defined in equation (\ref{Phidef}). Then, for all $T \in (0,\infty)$ and any bounded domain $D \in {\bf R}^{p}$,

\[ \sup_{0 < \epsilon < 1} \lim_{s \rightarrow +\infty} \sup_{0 \leq t \leq T} \sup_{{\bf x} \in D }  \left |\sum_{j=1}^{p} \frac{\partial}{\partial x_j}\int_{{\bf R}^{p}}\frac{1}{(2\pi s)^{p/2}}e^{-|{\bf x} - {\bf y}|^2/2s} \Phi_{p+1}^{(\epsilon)}(t; y_1,\ldots,
y_{p}, y_j) d{\bf y}  \right |
= 0.\]
\end{Propn}

\noindent In the following, the $t$ in the notation will be suppressed; $\Phi_{p+1}^{(\epsilon)}(x_1,\ldots, x_{p+1})$ will be used to denote $\Phi_{p+1}^{(\epsilon)}(t ; x_1,\ldots, x_{p+1})$ and it will be assumed that $0 < t < T$ where $T < +\infty$.  

\paragraph{Proof of Proposition~\ref{help2}}  Set 

\[ \tilde{\Phi}_{p+1}^{(\epsilon)} (s; x_1, \ldots, x_{p+1})  = \int_{{\bf R}^{p-1}} \frac{1}{(2\pi s )^{(p - 1)/2}}e^{-\frac{1}{2s}\sum_{j=1}^{p - 1} y_j^2 } \Phi_{p+1}^{(\epsilon)} (  x_1 + y_1, \ldots , x_{p-1} + y_{p-1}, x_p, x_{p+1}  ) d{\bf y} \] 

\noindent  and set

\begin{eqnarray*} \psi^{(\epsilon)} (s; x_1,\ldots, x_{p-1},x_p, x_{p+1}) &=& \int_{-\infty}^\infty \frac{1}{(2\pi s)^{1/2}}e^{-z^2/2s} \tilde{\Phi}_{p+1}^{(\epsilon)}(s;x_1 + z, \ldots, x_{p-1} + z, x_p, x_{p+1})   dz \\
 &=&\int_{-\infty}^\infty \frac{1}{(2\pi s)^{1/2}}e^{-z^2/2s} \tilde{\Phi}_{p+1}^{(\epsilon)}(s;x_1 , \ldots, x_{p-1}  , x_p + z, x_{p+1} + z)   dz.
\end{eqnarray*}

\noindent Using the bound in equation (\ref{lips2}), together with equation (\ref{phidefn}), it is straightforward that 
\[ \sup_{0 \leq \epsilon < 1} \sup_{{\bf x} \in {\bf R}^{p+1}} \left |\frac{\partial}{\partial x_j} \Phi_{p+1}^{(\epsilon)}(x_1,\ldots, x_{p+1})\right | < TK(p+1,T).\]

\noindent From this, it follows directly that for all $(j,k) \in \{1, \ldots, p-1\}^2$, 

\begin{eqnarray}\nonumber \lefteqn{\sup_{0 < \epsilon < 1} \sup_{{\bf x} \in {\bf R}^{p+1}}\left | \frac{\partial^2}{\partial x_j \partial x_k} \psi^{(\epsilon)} (s; x_1, \ldots, x_{p -1},x_p,x_{p+1})\right |}\\&& \nonumber
 = \sup_{0 < \epsilon < 1}\sup_{{\bf x} \in {\bf R}^{p+1}}\left | \frac{\partial^2}{\partial x_j \partial x_k} \int_{{\bf R}}\frac{1}{(2\pi s)^{1/2}}e^{-z^2/2s} \int_{{\bf R}^{p-1}}\frac{1}{(2\pi s)^{(p-1)/2}}e^{-\frac{1}{2s}\sum_{j=1}^{p-1}(x_j - y_j - z)^2} \right. \\&& \nonumber \hspace{90mm} \times \left. \Phi_{p+1}^{(\epsilon)}(y_1,\ldots, y_{p-1},x_p,x_{p+1}) d{\bf y} dz \right | \\
&& \nonumber \leq \sup_{0 < \epsilon < 1}\sup_{{\bf x} \in {\bf R}^{p+1}}   \int_{{\bf R}}\frac{1}{(2\pi s)^{1/2}}e^{-z^2/2s} \int_{{\bf R}^{p-1}}\frac{1}{(2\pi s)^{(p-1)/2}}e^{-\frac{1}{2s}\sum_{j=1}^{p-1}(x_j - y_j - z)^2}   \\&& \nonumber \hspace{60mm} \times \left |\frac{x_k - y_k - z}{s} \right |\left|\frac{\partial}{\partial y_j} \Phi_{p+1}^{(\epsilon)}(y_1,\ldots, y_{p-1},x_p,x_{p+1})\right | d{\bf y} dz   \\
&& \nonumber  \leq  TK(p+1,T) \int_{-\infty}^\infty \frac{|y|}{s} \frac{1}{(2\pi s)^{1/2}}e^{-|y|^2/2s} dy \\
&&  \label{zerlim} =  \sqrt{\frac{2}{\pi s}} TK(p+1,T) \stackrel{s \rightarrow \infty}{\longrightarrow} 0.
\end{eqnarray}

\noindent Set 

\[ \gamma^{(\epsilon)}(s ;x_2,\ldots, x_{p-1},x_p,x_{p+1}) = \left .\frac{\partial}{\partial x_1} \psi^{(\epsilon)}(s ;x_1,x_2, \ldots, x_{p-1},x_p,x_{p+1})\right |_{x_1 = 0},\]

\noindent then it follows directly from (\ref{zerlim}) that 

\begin{equation}\label{somebound}\lim_{s \rightarrow +\infty} \sup_{(x_2,\ldots, x_{p+1}) \in {\bf R}^p}\max_{j \in \{2, \ldots, p-1\}}\sup_{0 < \epsilon < 1}\left | \frac{\partial}{\partial x_j} \gamma^{(\epsilon)} (s ;x_2,\ldots, x_{p-1},x_p,x_{p+1})  \right | = 0 \end{equation}

\noindent and hence that   for any bounded $D \subset {\bf R}^{p-1}$ and all $(x_p, x_{p+1}) \in {\bf R}^2$, 

\begin{equation}\label{somebound2} \lim_{s \rightarrow +\infty} \sup_{(x_2,\ldots, x_{p-1}) \in D} \sup_{0 < \epsilon < 1} \left | \gamma^{(\epsilon)} (s ;x_2,\ldots, x_{p-1},x_p,x_{p+1}) - \gamma^{(\epsilon)} (s ;0,\ldots,0, x_p, x_{p+1})\right | = 0. \end{equation}

\noindent From equation (\ref{somebound}), it follows by Taylor's expansion theorem that for any bounded subset  $D \subset {\bf R}^{p-1}$, 

\begin{eqnarray*} \lefteqn{\psi^{(\epsilon)}(s; x_1,\ldots, x_{p-1},x_p,x_{p+1}) = \psi^{(\epsilon)}(s;0,\ldots,0,x_p,x_{p+1})}\\&& = \gamma^{(\epsilon)}(s;0,\ldots,0,x_p,x_{p+1}) \sum_{j=1}^{p-1}x_j  + \frac{1}{2}\sum_{j,k = 1}^{p-1} \partial^2_{jk}\psi^{(\epsilon)}(s;x_1^*, \ldots, x_{p-1}^*, x_p, x_{p+1})
\end{eqnarray*}

\noindent where $|x_j^*| \leq |x_j|$ for $j = 1, \ldots, p-1$. It follows from (\ref{zerlim}) that for any bounded $D \subset {\bf R}^{p+1}$, 

\begin{equation} \label{somebound3} \lim_{s  \rightarrow +\infty}\sup_{{\bf x}\in D}  \sup_{0 < \epsilon < 1}   \left | \psi^{(\epsilon)} (s ;x_1,\ldots, x_{p-1}, x_p, x_{p+1}) - \gamma^{(\epsilon)}(s ; 0, \ldots, 0, x_p, x_{p+1}) \sum_{j=1}^{p-1} x_j\right | = 0. 
\end{equation}

\noindent Taylor's expansion theorem applied to the first $p-1$ variables of $\partial \psi^{(\epsilon)}(x_1,\ldots, x_{p-1}, x_p, x_{p+1})$ gives that there are points $x_1^*, \ldots, x_{p-1}^*$ such that $0 \leq |x_j^*| \leq |x_j|$ such that 

\[ \partial_j \psi^{(\epsilon)}(x_1,\ldots, x_{p-1}, x_p, x_{p+1}) = \gamma^{(\epsilon)}(0, \ldots, 0, x_p, x_{p+1}) + \sum_{k=1}^{p-1} x_k \partial^2_{jk} \psi^{(\epsilon)}(x_1^*, \ldots, x_{p-1}^*, x_p, x_{p+1})\]

\noindent from which it follows, using equations (\ref{zerlim}) and    (\ref{somebound2}), that for any bounded set $D \subset {\bf R}^{p+1}$,

\begin{eqnarray*}\lefteqn{ \lim_{s \rightarrow +\infty} \sup_{{\bf x} \in D } \sup_{0 < \epsilon < 1}  \left |\frac{\partial}{\partial x_j} \psi^{(\epsilon)} (x_1,\ldots, x_{p-1}, x_p, x_{p+1}) - \gamma^{(\epsilon)}(s ;0,\ldots,0,x_p,x_{p+1})  \right | = 0,} \\&& \hspace{110mm}  j \in \{2,\ldots, p-1\}.
\end{eqnarray*}

\noindent By construction, $\psi(s ; x_1,\ldots, x_{p-1},x_p, x_{p+1}) = \psi(s ; x_{\sigma(1)} + a, \ldots, x_{\sigma(p-1)} + a, x_{p+1} + a , x_p + a)$ for any permutation $\sigma$ of $\{1, \ldots, p-1\}$ and any $a \in {\bf R}$ and all $(x_1,\ldots, x_{p+1}) \in {\bf R}^{p+1}$. From equation (\ref{somebound3}), it follows that for any bounded $D \subset {\bf R}^{p+2}$, 

\begin{eqnarray}\label{somebound4}\lefteqn{\lim_{s  \rightarrow +\infty}  \sup_{0 < \epsilon < 1}\sup_{(x_1,\ldots, x_{p+1}, a) \in D}}\\&&\left | \psi^{(\epsilon)}(s ;x_1,\ldots, x_{p-1}, x_p, x_{p+1}) - \gamma^{(\epsilon)}(s ; 0, \ldots, 0, a + x_p,a +  x_{p+1}) \sum_{j=1}^{p-1} (a + x_j)\right | = 0. \nonumber
\end{eqnarray}

\noindent It follows that 

\begin{eqnarray*}\lefteqn{ \lim_{s \rightarrow +\infty}\sup_{0 < \epsilon < 1} \sup_{(x_1,\ldots, x_{p+1},a) \in D}\left |\left (\gamma^{(\epsilon)}(s ;0, \ldots, 0, x_p, x_{p+1}) - \gamma^{(\epsilon)}(s ;0, \ldots, 0, x_p+ a , x_{p+1}+ a) \right ) \sum_{j=1}^{p-1} x_j\right . }\\&& \hspace{75mm} \left. - (p-1) a \gamma^{(\epsilon)} (s ;0, \ldots, 0, x_p + a, x_{p+1}+a) \right | = 0.
\end{eqnarray*}

\noindent From this, it follows (by considering the coefficient of $\sum_{j=1}^{p-1} x_j$) that 

\begin{eqnarray*}\lefteqn{\lim_{s \rightarrow +\infty}\sup_{0 < \epsilon < 1} \left | \gamma^{(\epsilon)}(s ;0, \ldots, x_p,  x_{p+1}) - \gamma^{(\epsilon)}(s ;0, \ldots, 0, x_p+ a , x_{p+1}+ a)\right | = 0}\\&& \hspace{100mm} \forall (a, x_p, x_{p+1}) \in {\bf R}^3\end{eqnarray*}

\noindent and hence that

\[\lim_{s \rightarrow +\infty}\sup_{0 < \epsilon < 1} \left | \gamma^{(\epsilon)}(s ;0, \ldots, 0,0,z) - \gamma^{(\epsilon)}(s ;0, \ldots, 0, x , x+z)\right | = 0\qquad \forall (x,z) \in {\bf R}^2.
\]

\noindent Using this, it follows by considering the term $(p-1) a \gamma^{(\epsilon)}(s ;0, \ldots, 0, x_p + a, x_{p+1}+a)$, that 

\[ \lim_{s \rightarrow +\infty} \sup_{0 < \epsilon < 1} |\gamma^{(\epsilon)}(s ;0,\ldots, 0, x_p, x_{p+1}) | = 0 \qquad \forall (x_p, x_{p+1}) \in {\bf R}^2,\]

\noindent from which it follows that for any bounded domain $D \subset {\bf R}^{p+1}$ and all $j \in \{1, \ldots, p-1\}$

\[ \lim_{s \rightarrow +\infty} \sup_{{\bf x} \in D} \sup_{0 < \epsilon < 1} \left | \frac{\partial}{\partial x_j} \psi^{(\epsilon)}(s ;x_1,\ldots, x_{p-1}, x_p, x_{p+1}) \right | = 0.\] 

\noindent Now note that 

\[\frac{\partial}{\partial x_p}\int_{{\bf R}^{p}}\frac{1}{(2\pi s)^{p/2}}e^{-|{\bf x} - {\bf y}|^2/2s} \Phi_{p+1}^{(\epsilon)}( y_1,\ldots,
y_{p}, y_p) d{\bf y}  = \frac{d}{d x_p}\psi^{(\epsilon)} (s; x_1  , \ldots, x_{p-1} , x_p , x_p)\]

\noindent where the differential $\frac{d}{dx_p}$ refers to both appearances of the variable $x_p$ and, since 
\[ \psi^{(\epsilon)} (s; x_1  , \ldots, x_{p-1} , x_p , x_p)  = \psi^{(\epsilon)} (s; x_1-x_p  , \ldots, x_{p-1} - x_p , 0 , 0),\]

\noindent  it follows that 

\[\frac{\partial}{\partial x_p}\int_{{\bf R}^{p}}\frac{1}{(2\pi s)^{p/2}}e^{-|{\bf x} - {\bf y}|^2/2s} \Phi_{p+1}^{(\epsilon)}( y_1,\ldots,
y_{p}, y_p) d{\bf y} = -\sum_{j=1}^{p-1} \frac{\partial}{\partial x_j}\psi^{(\epsilon)} (s; x_1-x_p, \ldots, x_{p-1} - x_p,0,0),
\]

\noindent from which, for any bounded $D \subset {\bf R}^p$, 

\[\lim_{s \rightarrow +\infty} \sup_{0 < \epsilon < 1} \sup_{{\bf x} \in D} \left | \frac{\partial}{\partial x_p}\int_{{\bf R}^{p}}\frac{1}{(2\pi s)^{p/2}}e^{-|{\bf x} - {\bf y}|^2/2s} \Phi_{p+1}^{(\epsilon)}( y_1,\ldots,
y_{p}, y_p) d{\bf y} \right | = 0.\]

\noindent Let 

\[ \tau_j(k) = \left\{\begin{array}{ll} k & k \neq j \\ p & k = j \\ j & k = p \end{array}\right. \] 

\noindent The result now follows by noting that for $j \in \{1, \ldots, p\}$, 
\[ \Phi_{p+1}^{(\epsilon)}(x_1,\ldots, x_p, x_j) = \Phi_{p+1}^{(\epsilon)}(x_{\tau_j(1)}, \ldots, x_{\tau_j(p)}, x_j),\]

\noindent from which it follows that

\[ \sup_{0 < \epsilon < 1} \lim_{s \rightarrow +\infty} \sup_{0 \leq t \leq T} \sup_{{\bf x} \in D }  \left |  \frac{\partial}{\partial x_j}\int_{{\bf R}^{p}}\frac{1}{(2\pi s)^{p/2}}e^{-|{\bf x} - {\bf y}|^2/2s} \Phi_{p+1}^{(\epsilon)}(t; y_1,\ldots,
y_{p}, y_j) d{\bf y}  \right |
= 0 \]

\noindent for each $j = 1, \ldots, p$ and hence that 

\[ \sup_{0 < \epsilon < 1} \lim_{s \rightarrow +\infty} \sup_{0 \leq t \leq T} \sup_{{\bf x} \in D }  \left |\sum_{j=1}^{p} \frac{\partial}{\partial x_j}\int_{{\bf R}^{p}}\frac{1}{(2\pi s)^{p/2}}e^{-|{\bf x} - {\bf y}|^2/2s} \Phi_{p+1}^{(\epsilon)}(t; y_1,\ldots,
y_{p}, y_j) d{\bf y}  \right |
= 0.\]
 \qed

\paragraph{Proof of theorem~\ref{last}}  Lemma~\ref{secondlast} gave 

\begin{eqnarray*}\lefteqn{ \left | M_p(t) - 
\frac{p(p-1)}{2}(-\Gamma^{\prime \prime}(0))\int_0^t
M_{p-2}(\alpha) d\alpha\right |}\\&& \leq \limsup_{\epsilon \rightarrow 0} \left |  \lim_{s
\rightarrow +\infty} \frac{1}{2s}   \left(\sum_{j=1}^p\int_0^s
  P_r \left (\frac{d}{d x_j}
 \Phi_{p+1}^{(\epsilon)} (t; x_1,\ldots, x_p ,x_j)  \right ) 
dr\right ) \right |.
\end{eqnarray*}

\noindent By proposition~\ref{help2}, it follows that for each
$r > 0$ and each $T < + \infty$, ${\bf x} \in {\bf R}^p$, 

\[ \lim_{s \rightarrow +\infty} \sup_{0 \leq t \leq T} \sup_{0 < \epsilon < 1}\left | P_{sr}  \frac{d}{d x_j}
 \Phi_{p+1}^{(\epsilon)}(t; x_1,\ldots, x_p ,x_j)  \right | = 0.\]

\noindent This, together with the the uniform bound  
\[  \sup_{0 \leq t \leq T} \sup_{(x_1,\ldots, x_p) \in {\bf R}^p} \left | P_{sr}  \frac{d}{d x_j}
 \Phi_{p+1}^{(\epsilon)}(t ; x_1,\ldots, x_p ,x_j)  \right | \leq 2 T K(p+1,T)\]

\noindent and the bounded
convergence theorem, imply that
\begin{eqnarray*} 
\lefteqn{\lim_{s \rightarrow +\infty} \frac{1}{2s} \int_0^s
P_r \frac{d}{dx_j} \Phi_{p+1}^{(\epsilon)}( t  ; x_1,\ldots,  x_p, x_j) dr  }\\&&
=  \lim_{s \rightarrow +\infty}   \int_0^1
P_{sr} \frac{d}{dx_j} \Phi_{p+1}^{(\epsilon)} (t  ; x_1,\ldots, x_p, x_j) dr   \\&& = 0.
\end{eqnarray*}

\noindent It follows that  

\begin{equation}\label{diffequ}   M_p (t) = 
\frac{p(p-1)}{2}(-\Gamma^{\prime
\prime}(0)) \int_0^t M_{p-2}(\alpha) d\alpha.\end{equation}

\noindent Since
$M_1(t)
\equiv 0$ and
$M_0(t)
\equiv 1$, the result follows for $t \in [0,T]$, by solving the system of
equations (\ref{diffequ}). By taking $T$ arbitrarily large, the result holds
for all $0 < t < +\infty$. 
\qed \vspace{5mm}

\noindent Theorem~\ref{secondlast2} is now considered. Firstly, it is    shown that   the solutions to the
equation (\ref{burger}) converge in $L^p$ norm as $\epsilon \rightarrow
0$, for all $p \geq 2$.\vspace{5mm}

\noindent Let 
$\tilde{C} : {\bf R}_+ \times \Omega \rightarrow {\bf R}_+$ denote $\tilde{C}(1,.)$ from equation (\ref{ceetildup}) and recall that for each $t \in {\bf R}_+$ solutions to equation (\ref{burger}) satisfy 

\[ \sup_{0 \leq s \leq t} \sup_{x \in [0,2\pi)} \sup_{0 < \epsilon < 1}
|u^{(\epsilon)}(s,x)| \leq \tilde{C}(t)\]

\noindent and that, from the inequality (\ref{ceedeebnd}), the definition in equation (\ref{dee}) and the computation below equation (\ref{ceedeebnd}) that      $E_{\bf Q}\left\{\tilde{C}^p(t) \right \} < +\infty$. The bounds are given in
theorem~\ref{bound1}. Furthermore, the result of theorem~\ref{last} is that 

\[ \lim_{\epsilon \rightarrow 0} \frac{1}{2\pi} E_{\bf Q} \left \{ \int_0^{2\pi}
|u^{(\epsilon)}(t,x)|^{2p} dx \right \} = \left (\prod_{j=1}^p (2j - 1) \right )\left (-\Gamma^{\prime
\prime}(0) \right )^p t^p.\]

\noindent The following argument shows convergence of $u^{(\epsilon)}$ in norm as $\epsilon \rightarrow 0$ in the $L^p$ spaces. The following result is required.

\begin{Th} \label{thkolmkr} Let $v: [0,1] \times [0,2\pi] \rightarrow {\bf R}$ be a process for which there are three strictly positive constants
 $\gamma, c, \delta$ such that 
\[ E_{\bf Q}\left \{\left |v^{(\epsilon_1)}(x_1) - v^{(\epsilon_2)}(x_2)\right |^\gamma \right \} \leq c (|\epsilon_1 - \epsilon_2| + |x_1 - x_2|)^{2 + \delta}\]
\noindent then there is a modification $\hat{v}$ of $v$ such that
\[ E_{\bf Q}\left \{\left(\sup_{(\epsilon_1, x_1) \neq (\epsilon_2, x_2)} \frac{|\hat{v}^{(\epsilon_1)}(x_1) - \hat{v}^{(\epsilon_2)}(x_2)|}{(|x_1 - x_2| +  |\epsilon_1 - \epsilon_2|)^\alpha}\right)^\gamma\right \}< + \infty\]
for all $\alpha \in [0, \frac{\delta}{\gamma})$.
\end{Th}

\paragraph{Proof} This is a standard result and may be found, for example, as theorem (2.1) on page 26 of Revuz and Yor~\cite{RY}. There it is presented for $v: [0,1] \times [0,1] \rightarrow {\bf R}$; the rescaling is standard. \qed

\begin{Lmm}\label{eltwoconvergence} Consider the system of equations for $\epsilon > 0$: 
\begin{equation}\label{eqvsyst}\left\{\begin{array}{l}
\partial_t v^{(\epsilon)} = \frac{\epsilon}{2} v_{xx}^{(\epsilon)} dt - \frac{1}{2} u^{(\epsilon)2} dt + \partial_t \zeta \\ v^{(\epsilon)}(0,x) \equiv 0
                       \end{array} \right. \end{equation}
\begin{equation}\label{eqvsyst2}\left\{\begin{array}{l}
\partial_t u^{(\epsilon)} = \frac{\epsilon}{2} u_{xx}^{(\epsilon)} dt - \frac{1}{2} (u^{(\epsilon)2})_x dt + \partial_t \zeta_x \\ u^{(\epsilon)}(0,x) \equiv 0
                       \end{array} \right. \end{equation}
The solutions
$v^{(\epsilon)}$ and $u^{(\epsilon)}$ to equations (\ref{eqvsyst}) and (\ref{eqvsyst2}) respectively  in ${\cal S}_p^*$ (defined in equation (\ref{eqspstar}) in the statement of lemma~\ref{heatlem}) converge in
$L^p$ norm to functions $v$ and $u$ respectively, which respectively satisfy 

\begin{equation}\label{eqveee} \left\{ \begin{array}{l}\partial_t v +
\frac{1}{2}u^{ 2} dt =  \partial_t \zeta \\ v(0,x) \equiv 0 \end{array}\right.\end{equation}

\noindent and 
\begin{equation}\label{equuuuu}  \left\{\begin{array}{l} \partial_t u +
\frac{1}{2}(u^{ 2})_x dt =  \partial_t \zeta_x \\
 u(0,x) \equiv 0. \end{array} 
\right.
\end{equation}

\end{Lmm}

\paragraph{Proof} Firstly, theorem~\ref{thuni} gives that for $\epsilon > 0$, equation (\ref{eqvsyst2}) has a unique solution in ${\cal S}_p^*$ for each $p > 0$, hence equation (\ref{eqvsyst}) has a unique solution in ${\cal S}_p^*$ for each $p > 0$, since once $u^{(\epsilon)}$ is established, equation (\ref{eqvsyst}) is linear and existence an uniqueness follows directly in a straightforward manner. Set
$\tilde{v}^{(\epsilon)} = \frac{\partial}{\partial \epsilon} v^{(\epsilon)}$.
It follows, simply by taking the derivative with respect to $\epsilon$
in equation (\ref{eqvsyst}) and using $v_x^{(\epsilon)} = u^{(\epsilon)}$ and  $v_{xx}^{(\epsilon)} =
u_x^{(\epsilon)}$, that 

\[ \left\{\begin{array}{l} \tilde{v}_t^{(\epsilon)} = \frac{\epsilon}{2}\tilde{v}_{xx}^{(\epsilon)} +
\frac{1}{2}u^{(\epsilon)}_x - u^{(\epsilon)} \tilde{v}_x^{(\epsilon)} \\ \tilde{v}^{(\epsilon)}(0,x) \equiv 0.
\end{array}\right.\]

\noindent Since $\frac{\epsilon}{2}\frac{\partial^2}{\partial x^2} -
u^{(\epsilon)}(t,x)\frac{\partial}{\partial x}$ is the infinitesimal generator, given in equation (\ref{infgen}) of the
process
$X$ from definition~\ref{proc},  it follows directly that

\begin{equation}\label{eqvtil} \tilde{v}^{(\epsilon)}(t,x) = \frac{1}{2} \int_0^t E_{\bf P} \left [
u_x^{(\epsilon)}(s, X_{s,t}(x)) 
\right] ds.\end{equation}

\noindent Now recall that

\[X_{s,t}^{(\epsilon)}(x) = x + \left(w^{(\epsilon)}_t - w^{(\epsilon)}_s\right) - \int_s^t u^{(\epsilon)} \left (r,X_{r,t}^{(\epsilon)}(x) \right )dr\]

\noindent (equation (\ref{process})). Taking the derivative in $x$ and  suppressing some appearances of $\epsilon$ in the notation, 

\[   X^\prime_{s,t}(x) = 1 - \int_s^t u_x(r,X_{r,t}(x))X^\prime_{r,t}(x)dr  \qquad \forall x \in {\bf R},  \]

\noindent giving

\[ X_{s,t}^\prime (x) = \exp\left\{ - \int_s^t u_x (r, X_{r,t}(x)) dr \right \} \]

\noindent and hence 

\[ \log X_{0,t}^\prime (x)   = -\int_0^t
u_x^{(\epsilon)}(r,X_{r,t}(x))dr.\]

\noindent It now follows from equation (\ref{eqvtil}) that  

\[ \tilde{v}^{(\epsilon)}(t,x) = -\frac{1}{2}E_{\bf P}[ \log X_{0,t}^\prime
(x)].\]

\noindent It follows that 

\begin{eqnarray*} \lefteqn{E_{\bf Q}\{ |\tilde{v}^{(\epsilon)}(t,x) |^{2p}\} }\\&&\leq 
 \frac{1}{2^{2p}}\left( E_{\bf Q}\left \{  E_{\bf P} \left [(\log
X_{0,t}^\prime(x))^{2p}\chi_{X^\prime_{0,t} (x) > 1} \right ]\right \}  + E_{\bf Q} \left \{ E_{\bf
P}\left [(\log X_{0,t}^\prime (x))^{2p} \chi_{X_{0,t}^\prime (x) < 1} \right ] \right \}\right)\\
&& =  I + II.
\end{eqnarray*}

\noindent For part $I$, note that 
\[ \frac{d}{dx} (\log x)^{2p} = 2p \frac{(\log x)^{2p-1}}{x}\]
\noindent and 
\begin{eqnarray*} \frac{d^2}{dx^2} (\log x)^{2p} &=& \frac{2p (2p-1) (\log x)^{2p-2}}{x^2} - \frac{2p (\log x)^{2p - 1}}{x^2}.\\
 &=& \frac{2p (\log x)^{2(p-1)}}{x^2}\left((2p-1) - \log x \right).
\end{eqnarray*}

\noindent  The maximum of $\frac{d}{dx} (\log x)^{2p}$ in the range $x \in (0,+\infty)$ occurs at $e^{2p-1}$ and is $2p (2p-1)^{2p-1}e^{-(2p-1)}$. It follows that, for $x \in [1, +\infty)$,
\[ (\log x)^{2p} \leq 2p (2p-1)^{2p-1}e^{-(2p-1)}(x-1) \leq (2p)^{2p} x.\]

\noindent Since $\int_0^{2\pi} X^\prime (x) dx = 2\pi$,  $X^\prime \geq
0$ and $E_{\bf Q}\{|\tilde{v}^{(\epsilon)}(t,x) |^{2p}\} = \frac{1}{2\pi}
\int_0^{2\pi} E_{\bf Q}\{|\tilde{v}^{(\epsilon)}(t,x) |^{2p}\} dx$ and, for $x
\geq 1$, $(\log x)^{2p} \leq (2p)^{2p} x$, it follows that

\[ I \leq (2p)^{2p}.\] 

\noindent Set $u_x = \phi$, then $\phi$ satisfies

\[ \left\{\begin{array}{l} \partial_t \phi  = \frac{\epsilon}{2}\phi_{xx} dt - \phi^2 dt - u\phi_x dt
+ \partial_t \zeta_{xx}\\ \phi(0,x) \equiv 0. \end{array}\right.
\]

\noindent Note that $u_x = \phi \leq w$, where $w$ satisfies

\begin{equation}\label{eqwww}\left\{\begin{array}{l} \partial_t w = \frac{\epsilon}{2}w_{xx} dt - u w_x dt + \partial_t \zeta_{xx}\\ w(0,x) \equiv 0. \end{array}\right. 
 \end{equation}

\noindent The solution to equation (\ref{eqwww}) may be written as

\begin{eqnarray*} w(t,x) &=& -\sum_{n \geq 0} n^2 a_n \left(\int_0^t
E_{\bf P}[
\cos(nX_{s,t}(x))]d_s\beta^{1n}(s) + \int_0^t E_{\bf P}[
\sin(nX_{s,t}(x))]d_s\beta^{2n}(s) 
\right)\\
&=& \theta(2;t,x)
\end{eqnarray*}

\noindent where $\theta$ is defined in equation (\ref{thetlabel}). Since

\[ X^\prime_{0,t}(x) = e^{-\int_0^t u_x(r,X_{t,}(x))dr} \geq e^{-\int_0^t
\theta (2; r,X_{t,}(x))dr},\]

\noindent it follows that 

\[ E_{\bf Q}\{E_{\bf P}[(\log (X_{0,t}^\prime \wedge 1))^{2p}]\} \leq
t^{2p} E_{\bf Q}\left\{ \left( \sup_{0 \leq s \leq t}\sup_x|\theta(2;s ,x
)|\right)^{2p}\right\},\]

\noindent which is bounded above independently of $\epsilon$ by an application of lemma~\ref{bound1lem}, so that 

\[ II \leq K(2p,t) < +\infty\]

\noindent where $K(2p,.)$ is an increasing function such that $K(2p,t) < +\infty$ for each $t < +\infty$, which is independent of $\epsilon$. It follows that for each $T < +\infty$

\begin{equation}\label{eqvetildbd} \sup_{0 \leq t \leq T} \sup_{0 < \epsilon \leq 1} E_{\bf Q} \left\{ \left | \tilde{v}^{(\epsilon)}(t,x)\right |^{2p} \right \} \leq (2p)^{2p} + K(2p,T) < +\infty.\end{equation}

\noindent Recall the definition of $\|f\|_p(t)$ given in equation (\ref{eqpnorm}).
From equation (\ref{eqvsyst}),

\[ \partial_t v^{(\epsilon)p} = \left(\frac{\epsilon
p}{2}v^{(\epsilon)p-1}v_{xx}^{(\epsilon)} -
\frac{p}{2}v^{(\epsilon)p-1}u^{(\epsilon)2} + \frac{p(p-1)}{2}
v^{(\epsilon)p-2} \Gamma(0)\right) dt+ pv^{(\epsilon)p-1}
\partial_t\zeta.\]

\noindent Integration by
parts and applications of Hölder's inequality yield that for positive integer $p$, 

\begin{eqnarray*} \frac{d}{dt} \|v^{(\epsilon)}\|_{2p}^{2p}(t) &=& - \epsilon p(2p-1) E_{\bf Q} \left\{ \left( \int v^{(\epsilon)2(p-1)} v_x^2 dx\right) \right \} - p E_{\bf Q} \left\{ \left ( \int v^{(\epsilon)2p - 1} u^{(\epsilon)2} dx \right) \right \} \\&& +p(2p - 1) \Gamma(0) E_{\bf Q} \left\{\left( \int v^{(\epsilon) 2p - 2} dx \right) \right \}  \\
&\leq & 
p \|v^{(\epsilon)}\|_{2p}^{2p-1}
\|u^{(\epsilon)}\|_{4p}^2 +
 p(2p-1) \|v^{(\epsilon)}\|_{2p}^{2p-2}
\Gamma(0),
 \end{eqnarray*}

\noindent so that

\[ \left\{\begin{array}{l} \frac{d}{dt} \|v^{(\epsilon)}\|_{2p}^2(t) \leq \|v^{(\epsilon)}\|_{2p} \|u^{(\epsilon)}\|_{4p}^2 + (2p-1) \Gamma(0) \\ \|v^{(\epsilon)}\|_{2p}(0) = 0.\end{array}\right. \]

\noindent It has already been established that, for $T < +\infty$, there is a constant $K(p,T) < +\infty$ such that 
\[ \sup_{0 \leq t \leq T}\|u\|_{4p}^2(t) \leq E_{\bf Q}\{\tilde{C}^{4p}(1,T)\}^{1/2p}
< K(p,T) < +\infty,\]

\noindent where  $\tilde{C}(b,t)$ is given by equation
(\ref{ceetildup}) and the existence of a finite upper bound $K(p,T)$ follows
from lemma~\ref{bound1lem}. It follows that

\[ \|v^{(\epsilon)}\|_{2p}^2 (t) \leq (1 + (2p-1) \Gamma(0) t) \exp\{ K(p,T) t\}\]

\noindent and hence that for any
$T < +\infty$ and any $p >1$, 
$v
\in L^p([0,T]
\times [0,2\pi]
\times
\Omega)$. That is, for each there is a positive function $C(p,t)$, increasing in $t$, with $C(p,t) < +\infty$ if $t < +\infty$ such that

\[ |\|v^{(\epsilon)}|\|_{T,p}:= \left(\int_0^T \int_0^{2\pi} E_{\bf Q} \left \{| v^{(\epsilon)}(t,x)|^p \right\} dx dt \right)^{1/p} < C(p,T).\], 

\noindent   Let 
\[K_1(p,t) = E_{\bf Q}\left\{\sup_{0 \leq s \leq t} \sup_{0 \leq x \leq 2\pi} \sup_{0 \leq \epsilon \leq 1} |u^{(\epsilon)}(s,x)|^{p}\right\}\] and note that, since $v^{(\epsilon)}(t,x_2) - v^{(\epsilon)}(t,x_1) = \int_{x_1}^{x_2} u^{(\epsilon)}(t,y) dy$, it follows by a standard application of Hölder's inequality that 

\begin{equation}\label{eqbd1}  \sup_{0 \leq t \leq T} \sup_{0 < \epsilon \leq 1} E_{\bf Q} \left\{ |v^{(\epsilon)}(t,x_2) - v^{(\epsilon)}(t,x_1)|^{2p}\right \} \leq |x_2 - x_1|^{2p} K_1(2p,T). \end{equation}

\noindent Let $K(p,T)$ denote the same quantity as in equation (\ref{eqvetildbd}). From equation (\ref{eqvetildbd}), it follows that for all $(\epsilon_1, \epsilon_2) \in (0,1]^2$, 
\begin{equation}\label{eqbd2} \sup_{0 \leq t \leq T}\sup_{0 \leq x \leq 2\pi} E_{\bf Q} \left\{ |v^{(\epsilon_1)}(t,x) - v^{(\epsilon_2)}(t,x)|^{2p}\right \} \leq |\epsilon_1 - \epsilon_2|^{2p} \left((2p)^{2p} + K(p,t)\right)^{2p},\end{equation}

   From equations (\ref{eqbd1}) and (\ref{eqbd2}), it follows  
by an application of theorem~\ref{thkolmkr} that 

\[ \lim_{\epsilon \rightarrow 0}\sup_{0
\leq t \leq T} E_{\bf Q}\left\{
\sup_{0
\leq x
\leq 2\pi}|v^{(\epsilon)}(t,x) - v(t,x)|^p\right\} = 0.\] 

\noindent \noindent from which 
\[ |\|v^{(\epsilon)} - v|\|_{T,p} \stackrel{\epsilon \rightarrow 0}{\longrightarrow} 0\]
\noindent for each $T < +\infty$ and each $p > 0$. It follows that ${\bf Q}$ almost surely, for all $T > 0$, $(m,n) \in {\bf Z}^2$, there is a random variable $\lambda_T(m,n)$ such that $E_{\bf Q} \left\{ |\lambda_T(m,n)|^p\right \} < +\infty$ for all $0 < p < +\infty$ and such that 
\[ \int_0^T \int_0^{2\pi} e^{ism\frac{2\pi}{T} + ixn} v^{(\epsilon)}(s,x) dx dt \rightarrow \lambda_T(m,n).\]
\noindent  Recall that 
$u^{(\epsilon)} = v_x^{(\epsilon)}$. Also,

\[ \sup_{0 \leq \epsilon \leq 1}\sup_{0 \leq s \leq t}\sup_{0 \leq x \leq
2\pi}|u^{(\epsilon)}(s,x)| \leq \tilde{C}(1,t),\]

\noindent where $\tilde{C}(b,t)$ is given by equation (\ref{ceetildup}).

\noindent By lemma~\ref{bound1lem}, for each $p > 0$, there is an increasing non negative function $K(p,.)$ such that $K(p,t) < +\infty$ for $t < +\infty$ and $E\{\tilde{C}^p(1,t)\} \leq K(p,t)$. It follows that ${\bf Q}$ - almost surely, $u^{(\epsilon)}$ converges weakly in $L^2([0,T] \times [0,2\pi])$ to

\[ u(t,x) = \frac{1}{2\pi T}\sum_{mn} -ine^{-(itm\frac{2\pi}{T} + inx)} \lambda_T(m,n).\] 

\noindent It is standard that ${\bf Q}$ - almost surely, the weak limits of the solutions of equations (\ref{eqvsyst}) and (\ref{eqvsyst2}) solve equations (\ref{eqveee}) and (\ref{equuuuu}) respectively. On $[0,T] \times [0,2\pi]$, let

\[ u^{(\epsilon)}(t,x) = \frac{1}{2\pi T} \sum_{nm} -in e^{(itm\frac{2\pi}{T} + inx)} \lambda_T^{(\epsilon)}(m,n).\]

\noindent Then 

\[ u^{(\epsilon)2}(t,x) = - \frac{1}{(2\pi T)}\sum_{mn}\sum_{m_1n_1} n_1(n-n_1)e^{(itm\frac{2\pi}{T} + inx)} \lambda_T^{(\epsilon)}(m_1,n_1)\lambda_T^{(\epsilon)}(m - m_1,n - n_1) \]

\noindent and

\begin{eqnarray*}\lefteqn{ \left |\sum_{n_1m_1} n_1(n-n_1)  \left(\lambda_T^{(\epsilon)}(m_1,n_1)\lambda_T^{(\epsilon)}(m - m_1,n - n_1) - \lambda_T (m_1,n_1)\lambda_T (m - m_1,n - n_1) \right)\right |}\\&&
 \leq \sum_{m_1n_1} \left |n_1\lambda^{(\epsilon)}_T(m_1,n_1)\right | 
 \left |(n-n_1)\lambda^{(\epsilon)}_T(m-m_1,n-n_1) - (n-n_1)\lambda_T (m-m_1,n-n_1) \right|\\&& \hspace{10mm} + \sum_{m_1n_1} \left |(n - n_1)\lambda_T (m- m_1,n - n_1)\right | \left |n_1\lambda^{(\epsilon)}_T ( m_1, n_1)   -n_1 \lambda_T ( m_1, n_1) \right|\\&& \leq  \left(\left (\sum_{m ,n } n^2 |\lambda^{(\epsilon)}_T (m , n )|^2\right)^{1/2} +  \left (\sum_{m ,n } n^2 |\lambda_T (m , n )|^2\right)^{1/2}\right)\\&& \hspace{50mm} \times  \left(\sum_{mn} n^2 \left |\lambda^{(\epsilon)}_T (m,n) - \lambda_T (m,n) \right|^2\right)^{1/2}.
\end{eqnarray*}

\noindent Firstly, 

\[ \sum_{m ,n } n^2 |\lambda^{(\epsilon)}_T (m , n )|^2 = \frac{1}{2\pi T} \int_0^T \int_0^{2\pi} |u^{(\epsilon)}(t,x)|^2 dx dt < \tilde{C}^2(T) \qquad \forall 0 \leq \epsilon < 1\]

\noindent where $\tilde{C}(T) = \sup_{0 \leq \epsilon < 1} \sup_{0 \leq t \leq T} \sup_{0 \leq x \leq 2\pi} |u^{(\epsilon)}(t,x)|$ and  $E_{\bf Q} \left\{\tilde{C}(T)^p\right\} < +\infty$ for all $0 < p < +\infty$. This implies that the dominated convergence theorem may be used on 

\[ \sum_{mn} n^2 \left |\lambda^{(\epsilon)}_T (m,n) - \lambda_T (m,n) \right|^2.\] 

\noindent Since $|\lambda_T^{(\epsilon)}(m,n) - \lambda_t(m,n)| \stackrel{\epsilon \rightarrow 0}{\longrightarrow} 0$ for eacn $(m,n)$, it follows that $u^{(\epsilon)2}$ converges to $u^2$ and hence that weak limits of the solutions of equations (\ref{eqvsyst}) and (\ref{eqvsyst2}) solve equations (\ref{eqveee}) and (\ref{equuuuu}) respectively. By integrating equations (\ref{eqvsyst}) and (\ref{eqveee}), it follows that 

\[ \int_0^T \int_0^{2\pi} u^{(\epsilon)2}(t,x) dx dt \stackrel{\epsilon \rightarrow 0}{\longrightarrow} \int_0^T \int_0^{2\pi} u^2(t,x) dx dt.\]

\noindent Furthermore, since ${\bf Q}$ - almost surely, $u$ is the weak limit of $u^{(\epsilon)}$ and $u \in L^2 ([0,T] \times [0, 2\pi])$, it follows that 

\[ \int_0^T \int_0^{2\pi} u^{(\epsilon)}(t,x) u(t,x) dx dt \stackrel{\epsilon \rightarrow 0}{\longrightarrow} \int_0^T \int_0^{2\pi} u^2(t,x) dx dt.\]

\noindent It follows that, ${\bf Q}$-almost surely, for all $T < +\infty$,

\begin{eqnarray*} \int_0^T \int_0^{2\pi} \left(u^{(\epsilon)}(t,x) - u(t,x) \right)^2 dx dt &=& \int_0^T \int_0^{2\pi} u^{(\epsilon)2}(t,x) dx dt + \int_0^T \int_0^{2\pi} u^{ 2}(t,x) dx dt\\&& - 2 \int_0^T \int_0^{2\pi} u^{(\epsilon)}(t,x)u(t,x) dx dt \stackrel{\epsilon \rightarrow 0}{\longrightarrow} 0. 
\end{eqnarray*}

\noindent It follows that ${\bf Q}$ - almost surely, $u^{(\epsilon)}(t,x) \stackrel{\epsilon \rightarrow 0}{\longrightarrow} u(t,x)$ for Lebesque almost all $(t,x) \in [0,T]\times [0,2\pi]$ for all $T < +\infty$. The dominated convergence theorem therefore gives 

\begin{eqnarray*} \lim_{\epsilon \rightarrow 0}E_{\bf
Q}\left \{\frac{1}{2\pi}\int_0^T \int_0^{2\pi} u^{(\epsilon)2p}(t,x) dx
dt\right \} &=&
 E_{\bf Q}\left \{\frac{1}{2\pi}\int_0^T
\int_0^{2\pi}
\lim_{\epsilon \rightarrow +\infty} u^{(\epsilon)2p}(t,x) dx dt \right \} \\
&=& E_{\bf
Q} \left \{\frac{1}{2\pi}\int_0^T
\int_0^{2\pi}
  u^{ 2p}(t,x) dx dt \right \},
\end{eqnarray*}

\noindent  From the bounds on $u^{(\epsilon)}(t,x)$ uniform in $(\epsilon, t,x) \in [0,1]\times [0,T] \times [0,2\pi]$,  convergence of $u^{(\epsilon)}$ to $u$ ${\bf Q}$ almost surely for Lebesgue almost all $(t,x) \in [0,T] \times [0,2\pi]$ and convergence of the   $L^p$ norms, it follows that $u^{(\epsilon)}$ converges to
$u$ in the $L^p$ {\em norm} topology for each $1 < p < +\infty$. Lemma~\ref{eltwoconvergence} is proved. 
   \qed 

\paragraph{Proof of theorem~\ref{secondlast2}} This follows directly from lemma~\ref{eltwoconvergence}. \qed

\section{The Invariant Measure for the Stochastic Burgers Equation}\label{sin}

The result in this article given by equation (\ref{momgrow}), concerning the 
growth of the moments for equation (\ref{burger}) is of interest, following
the results found in the article~\cite{EKMS}.  These results show
existence of an {\em invariant measure} for the viscosity solution of the inviscid Burgers equation

\[ \partial_t u + \frac{1}{2}(u^2)_x dt= \partial_t \zeta_x,\]

\noindent where the hypotheses on $\zeta$ in that article include the hypotheses
stated in hypothesis~\ref{hyp} of this article. The viscosity solution is the solution obtained by letting $\epsilon \rightarrow 0$ in equation (\ref{burger}). All moments
of the invariant measure  considered in that article exist,
as  outlined below. The argument presented by E, Khanin, Mazel and Sinai
in~\cite{EKMS} is based on Varadhan's theorem from large deviations. Starting from equation (\ref{heat}), solutions to equation 
(\ref{burger}) are given by the Cole Hopf transformation, equation
(\ref{CH}).   Consider the {\em action functional}

\begin{equation}\label{thaction} {\cal A}(\xi;0,t) = \frac{1}{2}\int_0^t
\dot{\xi}^2(s) ds +
\sum_{n=1}^\infty  a_n \left (\int_0^t \cos(n \xi(s))d_s\beta ^{1n}(s) + \int_0^t
\sin (n \xi(s))d_s\beta ^{2n}(s) \right ).
\end{equation}

\noindent It is a relatively straight forward application of Varadhan's theorem from Large Deviation theory to show
that 

\[ \lim_{\epsilon \rightarrow 0} - \epsilon \log U^{(\epsilon)}(t,x) = \inf_{\xi : \xi(t) =
x} {\cal A}(\xi; 0,t).\]

\noindent It is relatively standard, and is shown in the article~\cite{EKMS},  
that there is a trajectory  
$(\eta^{(t,x)}(s))_{s
\in [0,t]}$ that minimises ${\cal A}(.;0,t)$ subject to the constraint that
$\eta^{(t,x)}(t) = x$ and that, furthermore, this minimiser satisfies

\begin{equation}\label{theula}  \frac{\partial \eta^{(t,x)}}{\partial s}(s) = 
-
\sum_{n=1}^\infty na_n \left (\int_0^s \sin (n\eta^{(t,x)}(r))d_r \beta^{1n}(r) -
\int_0^s \cos (n\eta^{(t,x)}(r))d_r\beta^{2n}(r) \right ).
\end{equation} 

\noindent Set

\[ u(t,x) = \frac{\partial}{\partial x} {\cal A} (\eta^{(t,x)}(s); 0,t).\] 

\noindent By taking the derivative, integrating by parts and using equation 
(\ref{theula}),
it follows that 

\[ u(t,x) = \left. \frac{\partial \eta^{(t,x)}(s)}{\partial s}\right |_{s=t}.\]

\noindent By lemma~\ref{eltwoconvergence}, $u$ is the limit in $L^p$ norm of $u^{(\epsilon)}$, for any $p < +\infty$. Therefore $u$ satisfies 

\begin{equation}\label{invstbur}
\left\{\begin{array}{l} \partial_t u + \frac{1}{2}(u^2)_x dt = \partial_t \zeta_x \\u_0
\equiv 0.\end{array}\right.
\end{equation}

\noindent Using this, an upper bound is given in the article~\cite{EKMS}  for $\sup_{x
\in [0,2\pi]}|u(t,x)|$ and the distribution of this upper bound is shown to be
independent of $t$. It is also shown in~\cite{EKMS} that all the moments
of this distribution exist. An outline of the proof is reproduced here.
 
\begin{Th} \label{allmom} Let $u$ denote the solution to
equation (\ref{invstbur}). There exist
constants
$C(p) < +\infty$, independent of $t$, such that for all $t \in {\bf R}_+$, 

\[ E_{\bf Q} \left\{ \left (\sup_{0 \leq x \leq 2\pi}  | u(t,x) | \right )^{2p} \right\} \leq C(p).\]
\end{Th}  

\noindent {\bf Proof} The analysis follows that given in~\cite{EKMS}.  It is assumed,
following the discussion in~\cite{EKMS}, that the representation $u(t,x) = 
\frac{\partial
\eta^{(t,x)}}{\partial s}(s)|_{s=t}$ holds, where $u$ is the solution to
equation (\ref{invstbur}), $\eta$ minimises the action functional (\ref{thaction}) and
satisfies equation (\ref{theula}). Suppressing the superscripts, let
$0
\leq t_1
\leq t_2
\leq t$, then 

\begin{eqnarray*}\lefteqn{ \dot{\eta}(t_2) - \dot{\eta}(t_1) = - \sum_{n=1}^\infty
na_n
\int_{t_1}^{t_2} (\sin (n\eta(s)) d_s\beta^{1n}(s) -
\cos(n\eta(s))d_s\beta^{2n}(s))}\\ && = -\sum_{n=1}^\infty na_n
((\beta^{1n}(t_2)-\beta^{1n}(t_1))\sin (n\eta(t_2)) - 
(\beta^{2n}_{t_2}-\beta^{2n}_{t_1})\cos (n\eta(t_2)))\\&&
\hspace{5mm} + \sum_{n=1}^\infty n^2 a_n \int_{t_1}^{t_2} \dot{\eta}(s)
(\cos(n\eta(s))(\beta^{1n}(s) - \beta^{1n}(t_1)) +
\sin(n\eta(s))(\beta^{2n}(s) -
\beta^{2n}(t_1))) ds,
\end{eqnarray*}

\noindent so that, setting 

\[ C(s,t) = \sum_{n=1}^\infty n^2 |a_n|\left( \sup_{s\leq r_1 \leq r_2 \leq t}
|\beta^{1n}(r_2) - \beta^{1n}(r_1)|  + \sup_{s\leq r_1 \leq r_2 \leq t}
|\beta^{2n}(r_2) - \beta^{2n}(r_1)|\right ),\] 

\noindent  it follows that

\[ ||\dot{\eta}|(t_2) - |\dot{\eta}|(t_1)| \leq C(t_1, t_2) +
\int_{t_1}^{t_2}|\dot{\eta}| (s) C(t_1,s) ds.\] 

\noindent From this, it is straightforward to see that  for $s \in [t-1,t]$,

\begin{equation}\label{gineq1} \inf_{t-1 \leq s \leq t} |\dot{\eta}|(s) \geq
|\dot{\eta}|(t)e^{- C(t-1,t)} - C(t-1,t)
\end{equation}

\noindent and 

\begin{equation}\label{gineq2}
 \sup_{t-1 \leq s \leq t}|\dot{\eta}|(s) \leq
(|\dot{\eta}|(t) + C(t-1,t))e^{C(t-1,t)}.
\end{equation}

\noindent Now, the minimising trajectory minimises the action functional 

\begin{eqnarray*}\lefteqn{
{\cal A}(0,t;\xi) = \frac{1}{2}\int_0^t \dot{\xi}^2(s) ds + \sum_{n \geq 1} a_n
\int_0^t (\cos(n\xi(s))d_s\beta^{1n}(s) + \sin
(n\xi(s))d_s\beta^{2n}(s))}\\&& =    {\cal A}(0,t-1;\xi) +
\frac{1}{2}\int_{t-1}^t \dot{\xi}^2(s)ds\\&& \hspace{5mm} +
\sum_{n=1}^\infty  a_n ((\beta^{1n}(t) - \beta^{1n}(t-1))\cos (n\xi(t)) +
(\beta^{2n}(t)-\beta^{2n}(t-1))\sin (n \xi(t)))\\&& \hspace{5mm}
+ \sum_{n=1}^\infty na_n \int_{t-1}^t \dot{\xi}(s)\left((\beta^{1n}(s) -
\beta^{1n}(t-1)) \sin (n\xi(s)) - (\beta^{2n}(s) -
\beta^{2n}(t-1)) \cos (n\xi(s))\right) ds  
\end{eqnarray*}

\noindent subject to the constraint that $\xi(t) = x$.  Setting $\sup_{0 \leq x \leq
2\pi} |u(t,x)| = K $ and
$C = C(t-1,t)$ and using the inequalities (\ref{gineq1}) and (\ref{gineq2}), 

\[ {\cal A}(0,t;\xi) - {\cal A}(0,t-1;\xi) \geq \frac{1}{2}((|\dot{\xi}(t)|^2e^{-C} - C) \vee 0)^2 
-C -C||\dot{\xi}(t)|e^{-C} - C| . \]

\noindent To get an upper bound on ${\cal A}(0,t;\xi) - {\cal A}(0,t-1;\xi)$, where
$\xi$ is the minimiser of ${\cal A}(0,t; .)$ subject to $\xi(t) = x$, note that
$u(t,x)$ is $2\pi$ periodic in $x$ and  consider the constant velocity
trajectory $\eta$ such that $\eta(t) = x$ and
$\eta(t-1) =
\xi(t-1)$  to obtain

\[ {\cal A}(0,t;\xi) - {\cal A}(0,t-1; \xi) \leq 2\pi^2 + C + 2\pi C.\]

\noindent Since this holds for all $x \in [0,2\pi)$, it follows that 

\[ \frac{1}{2}((Ke^{-C} - C) \vee 0)^2 
-C -C|Ke^{-C} - C| \leq 2\pi^2 + C + 2\pi C.\]

\noindent It follows that either $K < Ce^C$, or

\[ \frac{1}{2}K^2 e^{-2C} - 2KCe^{-C} + \frac{3}{2}C^2  - C \leq 2\pi^2 + C(1 + 2\pi),\]
\noindent so that 
\[ (Ke^{-C})^2 -  4C (Ke^{-C}) + (3C^2 - 4(1+\pi)C - 4\pi^2) < 0\]

\noindent giving

\[ \sup_{0 \leq x \leq 2\pi} |u(t,x)| \leq K < 2Ce^C + e^C \sqrt{C^2 + 4(1 + \pi)C + 4\pi^2} \leq (3C  + 2(1+\pi))e^C < 10 e^{2C}.\]

\noindent To obtain estimates on $E_{\bf Q} \{ K^{p}\}$, first note that 

\[\sup_{t-1 \leq r \leq s
\leq t} |\beta^{1n}(s) - \beta^{1n}(r)| \leq 2 \sup_{t-1 \leq s \leq
t}|\beta^{1n}(s) -
\beta^{1n}(t-1)|,\]

\noindent  so that, setting $S^{jn} = \sup_{t-1 \leq s \leq t}|\beta^{jn}(s) -
\beta^{jn}(t-1)|$,

\[ E_{\bf Q}\{|K|^p\} \leq 10^p E_{\bf Q} \left\{ \exp\{4p\sum_{n=1}^\infty n^2 |a_n| (S^{1n} + S^{2n})\} \right \}
 = 10^p \prod_{n=1}^\infty E_{\bf Q} \left\{ \exp\{4p n^2 |a_n| S^{1n}\}\right\}^2.\]

\noindent An application of lemma~\ref{prel} gives 

\[ E_{\bf Q} \left\{|K|^p \right\} \leq 10^p\exp\left\{32 p^2 \sum_{n \geq 1}n^4 a_n^2 +
8p (\sqrt{2 \log 2} + 2 \sqrt{2\pi})\sum_{n \geq 1} n^2 |a_n|\right\},\]

\noindent concluding the proof of theorem~\ref{allmom}. This is the line of
the proof found in ~\cite{EKMS}.
\qed

\section{Conclusion} There is a striking dichotomy here. The steps that were sketched in section~\ref{sin} are justified elsewhere in the literature. This is 
discussed in ~\cite{EKMS}. The first of these is the application of Varadhan's theorem
in the stochastic case. This is straightforward; an integration by parts of
the potential term removes the stochastic integral. Secondly, the fact that
the action functional, in the stochastic case, has a minimiser and the fact
that this minimiser, in the stochastic case, satisfies the Euler - Lagrange
equations. These steps are relatively straight forward; after the `stochastic' integral has been removed by an integration by parts, the proof
depends on the $\dot{\xi}^2$ term and the fact that $|\dot{\xi}|$ is raised
to a power strictly greater than $1$, using standard arguments that date back to
Tonelli.  Once these steps are justified, the conclusion is that the Choice Axiom
leads to inconsistent results. The results from classical dynamics,
stating that a minimising trajectory for the action functional exists use
crucially the   relative weak compactness of the
unit ball in
$L^2$. By Tychonoff's theorem, the Choice Axiom implies relative weak
compactness of the unit ball in $L^2$. Kelley~\cite{Ke} showed that
relative weak compactness of the unit ball in $L^2$ implied the Choice
Axiom. It is the Choice Axiom, at the countable level, which is employed
in the arguments in this article. The conclusion is therefore that this article has
provided an example that demonstrates that employing the 
Choice Axiom leads to contradictory results in analysis. \vspace{5mm}

\noindent Electronic mail address for correspondence: {\tt jonob@mai.liu.se}\vspace{5mm}

\setlength{\baselineskip}{2ex}


\begin{thebibliography}{99} 

 
\bibitem{EKMS} Weinan E, K.Khanin, Mazel, Ya.G. Sinai [2000] {\em Invariant
Measure for Burgers Equation with Random Forcing}   Annals of Mathematics,
vol. 151 , pp 877 - 960

 
\bibitem{Ke}J.L. Kelley [1950]{\em The Tychonoff Product Theorem Implies the
Axiom of Choice} Fund. Math. 37 pp 75-76


 \bibitem{Kun} H. Kunita [1990]{\em Stochastic Flows and Stochastic
Differential Equations}\\Cambridge University Press
 
\bibitem{RY} D. Revuz and M. Yor [1999] {\em Continuous Martingales and
Brownian Motion} (3rd edition) \\ Springer

 
\end{thebibliography}
\end{document}